\documentclass[journal]{IEEEtran}
\usepackage{indentfirst}
\usepackage{mathrsfs}
\usepackage{dsfont}
\usepackage{graphicx}
\usepackage{epstopdf}
\usepackage{amsmath, amssymb, amsthm, bm}
\usepackage[pdftex,colorlinks,linkcolor=blue,anchorcolor=blue,citecolor=blue]{hyperref}
\usepackage{leftidx}
\usepackage{extarrows}
\usepackage{stfloats}
\usepackage{picinpar}
\usepackage{enumerate}
\usepackage{algorithm}
\usepackage{algorithmicx}
\usepackage{algpseudocode}
\usepackage{epsfig}
\usepackage{latexsym}
\usepackage{amsfonts}
\usepackage{enumerate}
\usepackage{graphics}
\usepackage{float}
\usepackage{pict2e}
\usepackage{tikz}
\usepackage{lineno,hyperref}
\usepackage{moreverb}
\usepackage{color}
\usepackage[nosort]{cite}
\usepackage{subfigure}
\usepackage{multirow}
\usetikzlibrary{arrows,automata}
\usepackage{booktabs}
\usepackage{graphicx}
\usepackage{colortbl}
\usepackage{diagbox}
\usepackage{framed}
\usepackage{units}
\usepackage{colortbl}
\usepackage{subfigure,graphicx}

\definecolor{mygray}{gray}{.87}

\newfont{\bb}{msbm10}

\DeclareMathOperator{\FP}{\mathds{F}}
\DeclareMathOperator{\GP}{\mathds{G}}
\DeclareMathOperator{\TP}{\mathds{T}}

\newcommand{\tr}{^{\sf T}}

\newcommand{\vb}{\bm{b}}

\newcommand{\vg}{\bm{g}}

\newcommand{\vu}{\bm{u}}
\newcommand{\vv}{\bm{v}}
\newcommand{\vw}{\bm{w}}
\newcommand{\vx}{\bm{x}}
\newcommand{\vy}{\bm{y}}
\newcommand{\vz}{\bm{z}}
\newcommand{\mA}{\bm{A}}
\newcommand{\mB}{\bm{B}}
\newcommand{\mC}{\bm{C}}
\newcommand{\mD}{\bm{D}}
\newcommand{\mE}{\bm{E}}

\newcommand{\mH}{\bm{H}}
\newcommand{\mI}{\bm{I}}
\newcommand{\mJ}{\bm{J}}
\newcommand{\mK}{\bm{K}}
\newcommand{\mL}{\bm{L}}

\newcommand{\mO}{\bm{O}}
\newcommand{\mP}{\bm{P}}
\newcommand{\mQ}{\bm{Q}}

\newcommand{\mS}{\bm{S}}

\newcommand{\mU}{\bm{U}}

\newcommand{\mW}{\bm{W}}

\newcommand{\mPi}{\bm{\Pi}}
\newcommand{\mGamma}{\bm{\Gamma}}

\newcommand{\sE}{\mathcal{E}}

\newcommand{\sG}{\mathcal{G}}

\newcommand{\sS}{\mathcal{S}}

\newcommand{\sV}{\mathcal{V}}

\newcommand{\RR}{\mathbb{R}}

\newcommand{\EE}{\mathbb{E}}

\newcommand{\vvarepsilon}{\bm{\varepsilon}}

\newcommand{\vlambda}{\bm{\lambda}}

\newcommand{\veta}{\bm{\eta}}

\newcommand{\one}{\mathbf{1}}
\DeclareMathOperator{\col}{col}
\DeclareMathOperator{\prox}{prox}
\DeclareMathOperator{\diag}{diag}
\DeclareMathOperator{\dist}{dist}

\allowdisplaybreaks[4]

\newtheorem{thm}{Theorem}
\newtheorem{ass}{Assumption}

\newtheorem{lem}{Lemma}

\newtheorem{rmk}{Remark}
\newtheorem{defn}{Definition}

\newtheorem{prop}{Proposition}

\ifCLASSINFOpdf

\else

\fi

\begin{document}
%
\title{Differentially Private Decentralized Optimization with Relay Communication}
%
%
%

\author{Luqing~Wang,
Luyao~Guo,
Shaofu~Yang,~\IEEEmembership{Member,~IEEE,}
Xinli~Shi,~\IEEEmembership{Senior Member,~IEEE.}
\thanks{This work was supported in part by the National Natural Science Foundation of China under Grant Nos. 62176056, 62003084, 62473098, and Young Elite Scientists Sponsorship Program by CAST, 2021QNRC001. (\emph{Corresponding author: Shaofu Yang.})}
\thanks{Luqing Wang is with the School of Computer Science and Engineering, Southeast University, Nanjing 210096, China
(e-mail: \href{mailto:luqing_wang@seu.edu.cn}{{luqing\_wang@seu.edu.cn}}).}
\thanks{Luyao Guo is with the School of Mathematics, Southeast University, Nanjing 210096, China
(e-mail: \href{mailto:ly_guo@seu.edu.cn}{{ly\_guo@seu.edu.cn}}).}
\thanks{Shaofu Yang is with the School of Computer Science and Engineering, Southeast University, Nanjing 210096, China
(e-mail: \href{mailto:sfyang@seu.edu.cn}{{sfyang@seu.edu.cn}}).}
\thanks{Xinli Shi is with the School of Cyber Science \& Engineering, Southeast University, Nanjing 210096, China
(e-mail: \href{mailto:xinli_shi@seu.edu.cn}{{xinli\_shi@seu.edu.cn}}).}

}

\maketitle


\begin{abstract}
Security concerns in large-scale networked environments are becoming increasingly critical.  To further improve the algorithm security from the design perspective of decentralized optimization algorithms, we introduce a new measure: Privacy Leakage Frequency (PLF), which reveals the relationship between communication and privacy leakage of algorithms, showing that lower PLF corresponds to lower privacy budgets. Based on such assertion, a novel differentially private decentralized primal--dual algorithm named DP-RECAL is proposed to take advantage of operator splitting method and relay communication mechanism to experience less PLF so as to reduce the overall privacy budget. To the best of our knowledge, compared with existing differentially private algorithms, DP-RECAL presents superior privacy performance and communication complexity. In addition, with uncoordinated network-independent stepsizes, we prove the convergence of DP-RECAL for general convex problems and establish a linear convergence rate under the metric subregularity. Evaluation analysis on least squares problem and numerical experiments on real-world datasets verify our theoretical results and demonstrate that DP-RECAL can defend some classical gradient leakage attacks.

\end{abstract}
%
\IEEEpeerreviewmaketitle

\section{Introduction}
Decentralized optimization, as a collaborative model entailing local computing and information exchanging, which is imperative to reduce storage costs, computation burdens and improve the robustness and scalability of network, has been widely adopted and deployed in various applications  \cite{Nedic2015}. However, the iterative process of decentralized optimization algorithms requires large, representative datasets owned by local agents, and then the communication among agents may lead to the leakage of sensitive information about the datasets, so that the adversaries can gain privacy from the shared model. Recent studies \cite{Melis2019,Zhu2019} reveal the unimaginable truths that not only the characteristics of data and their correlations, but even the original data can be precisely reversely inferred. Thus, privacy-preserving decentralized optimization algorithms have received a lot of attention.

A common class of approaches to guarantee privacy is to employ trusted third-party computation such as homomorphic encryption \cite{C. Zhang2019} and garbled circuit \cite{A. Yao1986}, which, nevertheless, will generate large computational burdens and incur exponential runtime overheads. Other results take advantage of the structural properties of decentralized network to ensure privacy. For example, the work \cite{Y.Ye2020} shows that time-varying parameters can be added to the projection or stepsizes, and noise that cover gradient can be constructed through the network structure. Although  such mechanism also guarantees the security of the algorithm, it is only applicable to specific network topology such as Hamiltonian graph or the requirement that the agents do not share information with adversaries. Therefore, we focus on the popular researches about differential privacy \cite{Dwork2014}.

Differential privacy adds inexact terms to shared message in order to interfere with the information obtained by the adversaries to achieve security in a probabilistic sense. \cite{H.-P. Cheng2019} and \cite{Lin2022} start to introduce differential privacy and client-level differential privacy into distributed learning systems for algorithm and agent privacy protection. However, there exists a contradiction that a small amount of differential privacy noise is insufficient to defend against strong privacy attacks, while a large amount of noise would impact the convergence rate and accuracy \cite{Cao2021}. Some existing works \cite{Abadi2016,Zhang2017,Huang2020} try to narrow the privacy budget of algorithms with new analysis of differential privacy theory. Other works \cite{ACheng2022,Xiao2023,Zhang2018_1} add additional restrictions or clipping to algorithms to minimize the effect of noise on accuracy. Only a few works attempt to reduce the noise coupled into the iterations while ensuring privacy preservation from the perspective of algorithmic framework design. \cite{Zhang2018} applies a linear approximation to modify the framework for even-numbered iterations on the basis of \cite{Zhang2018_1}, compromising privacy in only half of the iterations.

However, the work \cite{Zhang2018} retains a original ADMM method during its odd-numbered iterations to eliminate the information leakage occurring in even-numbered iterations, which imposes a significant computational burden when the internal optimization problem lacks an analytical solution. There is no neat and intuitive theoretical result revealing that the privacy budget of \cite{Zhang2018} has an order reduction compared to that of \cite{Zhang2018_1}. Additionally, existing research still lacks discussion on the relationship between communication and privacy leakage in algorithm design frameworks, failing to provide guidance on constructing a decentralized optimization algorithm with enhanced security and accuracy. Besides, some of the above algorithms are only applicable to convex and smooth functions or centralized optimization problems rather than general problems with nonsmooth terms in a decentralized form.

 In this paper, to shed new light on the connection between the communication strategy and the privacy performance of differentially private decentralized optimization algorithms, we define a new measure Privacy Leakage Frequency (PLF), which denotes the maximum number of times experienced by an agent to query local gradients and publish the messages during iterations of the algorithm. We point out a quantitative relationship between PLF and privacy budget, revealing that when converging to the fixed accuracy, the less PLF an algorithm has, the less privacy is compromised. Based on this assertion, we propose a novel privacy-preserving decentralized optimization algorithm, called the differentially private relayed and communication-efficient algorithm (DP-RECAL) to experience less PLF to achieve superior privacy preservation.

To the algorithm development, we first reformulate the considered problem as a composite optimization problem with equation constraints. Inspired by AFBS, which is a state-of-the-art operator splitting method for solving the zeros of monotone operators, we provide a new (centralized) primal-dual proximal algorithm to this reformulation. Then, we introduce the relay communication strategy, which not only allows decentralized implementation of our proposed algorithm, but also reduces the PLF of the algorithm. Furthermore, we add Gaussian noise as a perturbation term to the passed message and finally obtain the DP-RECAL. We highlight the following aspects of our proposed algorithm:

Compared to existing differentially private optimization algorithms \cite{Zhang2017,Huang2020,Zhang2018_1,Zhang2018}, we argue that DP-RECAL exhibits superior privacy-preserving performance and communication complexity. Besides, DP-RECAL adopts uncoordinated network-independent stepsizes and the private convergence stepsize condition of each agent is $\alpha_i\in(0,\frac{2}{L_i+1})$, which implies that the stepsizes can be decided locally by agents. In terms of its convergence analysis, we employ a new analytical technique to prove the convergence of DP-RECAL for general convex problems and establish a linear convergence rate under the metric subregularity.

 \emph{Organization:} The rest of our paper is organized as follows. Section \ref{Problem Formulation} introduces the optimization problem and generalizes the threat model. Then, we describe the concept of PLF through differential privacy theory and the motivation for constructing DP-RECAL. Section \ref{DP-RECAL} presents DP-RECAL in detail and make sufficient contrasts with related works. Section \ref{Convergence analysis} establishes the convergence analysis of DP-RECAL under general convexity and its linear convergence rate under metric subregularity, respectively. Section \ref{Privacy Analysis} studies the privacy performance of DP-RECAL. Moreover, Section \ref{Communication and privacy budget} compares the performance of RECAL with other related algorithms by listing the computational complexity, communication complexity and PLF in the least squares problem. Finally, Section \ref{Performance evaluation} implements two numerical experiments on real-world datasets and Section \ref{Conclusion} concludes this paper.

\emph{Notations:} $\RR^n$ denotes the $n$-dimensional vector space with inner-product $\langle \cdot,\cdot \rangle$. The $1$-norm and Euclidean norm are denoted as $\|\cdot\|_1$ and $\|\cdot\|$. $0$ and $\mI$ denote the null matrix and the identity matrix. For a given set $\{\vx_1,\cdots,\vx_n\}$, $\col\{\vx_1,\cdots,\vx_n\}=[\vx_1\tr,\cdots,\vx_n\tr]\tr$ and $\diag\{\vx_1,\cdots,\vx_n\}$ are denoted as a diagonal matrix consisting of blocks $\vx_1,\cdots,\vx_n$. For a matrix $\mA\in \RR^{n \times m}$, $\|\mA\| := \max_{\|\vx\|=1}\|\mA\vx\|$. For a symmetric matrix $\mA$, $\mA\succ 0$ means that $\mA$ is positive definite. If matrix $\mA \succ 0$, $\|\vx\|^2_{\mA}:=\langle \vx,\mA\vx \rangle$, $\mathcal{P}_{\mathcal{X}}^{\mA}(\vx):=\arg\min_{\vx'\in\mathcal{X}}\|\vx-\vx'\|_{\mA}$. If $f:\mathbb{R}^n\rightarrow \mathbb{R}\cup\{\infty\}$ is a proper, closed and convex function, its subdifferential $\partial f: \RR^n \rightrightarrows \RR^n$ is defined as
 $\partial f := \{\vy |~ f(\vx)+\langle \vy, \vz-\vx \rangle \leq f(\vz) , \forall \vz \in\RR^n\}$. Then, the resolvent of $\partial f$ called the proximal operator is denoted as $\prox_f^{\mA}(\vx)=\arg\min_{\vy \in \RR^n} \{f(\vy)+\frac{1}{2}\|\vx-\vy\|_{\mA}^2\}$, where $\mA$ is a positive definite matrix.

\section{Problem formulation}\label{Problem Formulation}
\subsection{Optimization model}
Throughout the paper, we consider optimization problems on an undirected graph $\sG=(\sV,\sE)$, where $\sV=\{1,2,\cdots,n\}$ is the set of $n$ agents, and $\sE$ is the set of $m$ edges. Let $\mD$ and $\mL$ denote the degree matrix and the Laplace matrix of graph $\sG$, respectively. The form of the optimization problem is as follows:
\begin{align}\label{pro1}
\min_{\vx \in \RR^q} ~~ \sum_{i=1}^nf_i(\vx;\mathcal{D}_i)+nr(\vx),
\end{align}
where $f_i:\mathbb{R}^q\rightarrow \mathbb{R}$ is a convex local loss function, $\mathcal{D}_i$ denotes a local dataset and $r: \mathbb{R}^{q}\rightarrow \mathbb{R}\cup\{\infty\}$ is a regular term. When there is no confusion, we will omit the $\mathcal{D}_i$ in $f_i$, i.e., $f_i(\vx):=f_i(\vx;\mathcal{D}_i)$. Then, we make the following assumption on local loss functions.
\begin{ass}\label{Ass1}
Each $f_i(\vx)$ is differentiable and $L_i$-smooth, i.e., for some $L_i>0$,
$$
\|\nabla f_i(\vx_1)-\nabla f_i(\vx_2)\|\leq L_i \|\vx_1-\vx_2\|, \forall i\in \mathcal{V},\forall \vx_1,\vx_2\in \mathbb{R}^q.
$$
The Lipschitz matrix $\mQ$ is stacked by $L_i$ of each agent, \textcolor[rgb]{0.00,0.00,0.63}{i.e.,} $\mQ=\diag\{L_1\mI,\cdots,L_n\mI\}$.
\end{ass}

\subsection{Threat model}\label{Attack model}
The adversary can be either a honest-but-curious agent of the network or an external eavesdropper with the goal of uncovering local private datasets. Such an adversary can be commonly seen in \cite{Abadi2016,Zhang2017,Huang2020,ACheng2022}. It is allowed to access and store the information passed between agents as well as the structure of optimization model, but must not interfere with the collaborative algorithm, such as modifying the model architecture and sending malicious parameters. The information passed between the agents will be gradients or state variables, and based on the state variables, the adversary can also indirectly infer the gradients from the iteration sequence.

The captured gradients can be exploited by some well-known gradient leakage attack algorithms \cite{Zhu2019,Geiping2020,Yin2021} to expose local private datasets, which often contain confidential information, such as the sensors' position in sensor network localization problem \cite{C. Zhang2019} and the patients' physical condition in disease classification problem \cite{Shruthi2019}. The core idea of these attack algorithms is to iterate the real dataset from the initialized dummy dataset by minimizing the distance between the dummy gradient and the real gradient based on real dataset.

This threat model with the attack algorithms will be adopted as an experimental simulation to test the security of our proposed method.
\subsection{Differential privacy}
We employ differential privacy as the base method for our analysis of algorithmic privacy. Differential privacy (DP) mechanism defined in \cite{Dwork2014}, a well-known standard approach for mathematical privacy-preserving, is gaining in importance that makes it nearly impossible for an adversary to identify a valid data record. DP guarantees that the addition
or removal of a single data sample does not substantially affect the outcome.
\begin{defn}\label{DP} 
A randomized algorithm $\mathcal{M}$ satisfies $(\epsilon,\delta)$-differential privacy if for any two neighbouring datasets $\mathcal{D}$ and $\hat{\mathcal{D}}$ differing in at
most one record and any output $E \in range(\mathcal{M})$, it holds that
\begin{align*}
\mathbf{P}(\mathcal{M}(\mathcal{D})=E) \leq e^{\epsilon} \mathbf{P}(\mathcal{M}(\hat{\mathcal{D}})=E)+\delta,
\end{align*}
where $\mathbf{P}(\cdot)$ is the probability measure and $\epsilon>0$.
\end{defn}
Intuitively, with a smaller $\epsilon$ the adversary can't tell whether the output $E$ is obtained from $\mathcal{D}$ or $\hat{\mathcal{D}}$, and hence almost impossibly infer the existence of any valid data record. In fact, $\epsilon$ of the algorithm resorting to DP mechanism is named as ``privacy budget '' in \cite{Dwork2014} and will be set in advance. During the implementation of the algorithm, some privacy budget is deducted for each access to datasets, and when the budget is used up, datasets can no longer be accessed. Thereby, $\delta$ and $\epsilon$ show the effectiveness of privacy protection by $(\epsilon,\delta)$-DP.

 In order to examine the privacy budget of an algorithm from a local to a holistic view, we focus on the following propositions of $(\epsilon,\delta)$-differential privacy.
\begin{prop}[Sequential combinations]\cite{Dwork2014}\label{sumCDP} \\
If randomized mechanisms $\mathcal{M}_1, \mathcal{M}_2, \cdots, \mathcal{M}_{K}$ are respectively $(\epsilon_1,\delta)$-DP, $(\epsilon_2,\delta)$-DP, $\cdots$, $(\epsilon_{\mathrm{K}},\delta)$-DP, the algorithm
\begin{align*}
M(\mathcal{D})=\mathcal{M}_1\circ \mathcal{M}_2 \circ \cdots \circ \mathcal{M}_{\mathrm{K}}(\mathcal{D})
\end{align*}
satisfies $(\sum_{k=1}^{\mathrm{K}}\epsilon_k,\delta)$-$DP$.
\end{prop}
\begin{prop}[Parallel combinations]\cite{Dwork2014}\label{sumCDP1} \\
Let a dataset $D=\bigcup_{h=1,2,\cdots,m} D_i$. If randomized mechanisms $\mathcal{M}_1, \mathcal{M}_2, \cdots, \mathcal{M}_m$ are respectively $(\epsilon_1,\delta)$-DP, $(\epsilon_2,\delta)$-DP, $\cdots$, $(\epsilon_m,\delta)$-DP, the algorithm
\begin{align*}
M(\mathcal{D})=(\mathcal{M}_1(\mathcal{D}_1), \mathcal{M}_2(\mathcal{D}_2), \cdots, \mathcal{M}_m(\mathcal{D}_m))
\end{align*}
satisfies $(\max_{i=1,\cdots,m}\epsilon_i,\delta)$-$ DP$.
\end{prop}
Proposition \ref{sumCDP} reveals that for an agent the privacy budget depends on the number of iterations leaking privacy (i.e., involved in the communication) and Proposition \ref{sumCDP1} reveals that the privacy budget of the algorithm, named as overall privacy budget, depends on the maximum privacy budget of the agents.

\subsection{Proposal of PLF}
Following the above properties in this article, we consider the privacy performance of algorithms from a new perspective of communication and define a new measure called privacy leakage frequency (PLF).

Taking the existing decentralized algorithms EXTRA \cite{PGEXTR} and TriPD-Dist \cite{Latafat2019} as examples, the PLF of EXTRA will degenerate to iteration number $K$ given that at each iteration all agents need to communicate and publish messages, while the PLF of TriPD-Dist is equivalent to $\max_{i \in\mathcal{V}} p_i K$ given that each agent $i$ is independently activated with probability $p_i$.
\vspace{-3mm}
\begin{flushleft}
\centering\fbox{
 \parbox{0.47\textwidth}{
\textbf{Privacy Leakage Frequency--PLF:}\\
The PLF $\Theta_K$ of a decentralized optimization algorithm after $K$ iterations is defined as
$$\Theta_K = \max_{i \in\mathcal{V}} \sum_{k=1}^{K} \omega_i(k),\  \omega_i(k)=\left\{\begin{array}{cl}
                                                                                           1, & \text{agent $i$ is activated}, \\
                                                                                           0, & \text{o.w.}
                                                                                         \end{array}\right.
$$
Agent $i$ is activated means the agent queries local gradients and publishes the treated information in the $k$-iteration.}
}\end{flushleft}

It is clear from Proposition \ref{sumCDP} and \ref{sumCDP1} that the less PLF the differentially private algorithm has, the less privacy is compromised. To show the quantitative relationship between PLF and privacy budget for differentially private algorithms, we give the following theorem.
\begin{thm}\label{LCI-pb}
Let $\texttt{ALG}$ be an optimization algorithm and $\tau_i(K)=\sum_{k=1}^{K} \omega_i(k)$ denote the number of times agent $i$ has been activated after $K$ iteration. If the privacy budget of agent $i$ activated at $\tau_i(k)$-th time during $k$-iteration is $\epsilon_{\tau_i(k)}$, then it holds that
\begin{align*}
\big [\text{Overall Privacy Budget of } \texttt{ALG}\big ]=\sum_{j=1}^{\Theta_K}\epsilon_j,
\end{align*}
where $\Theta_K$ is the PLF after $K$ iterations.
\end{thm}
\begin{IEEEproof}
Let $\sS_{i,k},\  k\in\{1,\ldots,K\}:=\mathcal{K}$ be the state of agent $i$ generated by an optimization algorithm in $k$-th round of iteration, and $\mathcal{A}_i\subset \mathcal{K}$ be the subsequence that agent $i$ is activated. When agent $i$ is activated, it will result in a privacy breach, whereas when agent $i$ remains inactive, as there is no message exchange, there is no risk of privacy leakage Therefore, from Proposition \ref{sumCDP}, the privacy budget of $\sS_{i,1}\circ\sS_{i,2}\circ\cdots\circ\sS_{i,K}$ for agent $i$ is $\sum_{k \in \mathcal{A}_i} \epsilon_{\tau_i(k)}+\sum_{k \in \mathcal{K}\backslash\mathcal{A}_i}0$, which, recalling the definition of $\omega_i(k)$, is equivalent to $\sum_{k=1}^K \omega_i(k)\epsilon_{\tau_i(k)}$. Note that $\Theta_K = \max_{i \in\mathcal{V}} \sum_{k=1}^{K} \omega_i(k)$ and $\tau_i(K)=\sum_{k=1}^{K} \omega_i(k)$. From Proposition \ref{sumCDP1}, we conclude that  the overall privacy budget of the optimization algorithm is
\begin{align*}
&\max_{i\in\mathcal{V}}\left\{\sum_{k=1}^K \omega_i(k)\epsilon_{\tau_i(k)}\right\}
=\max_{i\in\mathcal{V}} \sum_{j=1}^{\sum_{j=1}^{K}\omega_i(k)}\epsilon_j\\
&=\sum_{j=1}^{\max_{i\in\mathcal{V}}\sum_{k=1}^{K}\omega_i(k)}\epsilon_j=\sum_{j=1}^{\Theta_K}\epsilon_j.
\end{align*}
Thus, we complete the proof.
\end{IEEEproof}
  From this result, we can infer that the impact of PLF on the privacy budget is crucial. Moreover, it is worth pointing out PLF is related to both computation and communication, but is not equivalent to computational and communication complexity. This is because that if the node participates in the computation but does not share state information, there is no privacy leakage; if the node passes the same information in multiple rounds of communication, this corresponds to the privacy leakage from a single communication.

To our knowledge, existing differentially private algorithms only focus on the privacy budget per iteration and do not take a holistic view to reduce the overall privacy budget, much less design a new algorithm by reducing the PLF of iterations. Thereby, we desire to design a differentially private and decentralized optimization algorithm with less PLF to achieve better privacy performance for problem \eqref{pro1}.

\section{DP-RECAL}\label{DP-RECAL}
In this section, we propose a differentially private algorithm with efficient communication and draw comparisons with the existing related algorithms.
\subsection{Algorithm construction}
To illustrate the activity of each agent, we introduce variables $\vy_i \in \RR^q$ indicating the local copies of $\vx$ . By defining
$$
\vy=\col \{\vy_1,\vy_2,\cdots,\vy_n\},~ f(\vy)=\sum_{i=1}^nf_i(\vy_i),
$$
where the operation $\col(\cdot)$ represents stacking vectors column by column, problem \eqref{pro1} can be rewritten in a matrix form as
\begin{equation}\label{pro2}
\min ~f(\vy)+nr(\vx),
\text{ s.t. }\mU \vx-\vy=0,
\end{equation}
where $\mU = \one \otimes \mI$. 
Then, the Lagrangian of problem \eqref{pro2} is
\begin{align}\label{Lagrangian}
&\mathcal{L}(\vx,\vy,\vlambda)=nr(\vx)+f(\vy)+\langle \vlambda, \mU \vx-\vy\rangle,
\end{align}
where $\vlambda=\col\{\vlambda_1,\cdots,\vlambda_n\}$ is the Lagrange multiplier. For brevity, define $\vw=\col\{\vlambda,\vx,\vy\}$, and the following operators:
\begin{align*}
\mA: (\vlambda,\vx, \vy) &\mapsto (-\mU\vx+\vy,\mU\tr\vlambda, -\vlambda),\\
\FP: (\vlambda,\vx, \vy) &\mapsto (0,0,\nabla f(\vy)),\\
\GP: (\vlambda,\vx, \vy) &\mapsto (0,n\partial r(\vx),0 ).
\end{align*}
Then the saddle subdifferential is denoted as
$
\partial \mathcal{L}(\vw)= \mA\vw +\FP\vw+\GP\vw.
$
From the Karush-Kuhn-Tucker (KKT) conditions, we have that $\vw^*$ is a saddle point of $\mathcal{L}(\vw)$ if and only if $0 \in \partial\mathcal{L}(\vw^*)$. Define $\mathcal{W}^*=\mathcal{X}^* \times \mathcal{Y}^*\times {\Lambda}^*$ as the set of saddle-points, where $\mathcal{X}^* \times \mathcal{Y}^*$ is the primal optimal solution set and ${\Lambda}^*$ is the dual optimal solution set. Thus, problem \eqref{pro2} is equivalent to the following fundamental zero point problem:
\begin{align}\label{zero}
\mathrm{Find}~ \vw^*\in \mathbb{R}^{q} \times \mathbb{R}^{nq}\times \mathbb{R}^{nq}~ \mathrm{such~that} ~ 0 \in \partial\mathcal{L}(\vw^*).
\end{align}
According to the structure of $\partial\mathcal{L}$,
inspired by the prediction-correction structure of AFBS \cite{Latafat2017}, we propose the following variable metric fixed point iteration:
\begin{subequations}\label{eq:AFBS}
\begin{align}
&\vw^{k+\frac{1}{2}}=(\mH_1+\mA+\GP)^{-1}(\mH_1-\FP)\vw^k,\label{AFBS1}\\
&\vw^{k+1}=\vw^k+\mH_2(\vw^{k+\frac{1}{2}}-\vw^k)\label{AFBS2},
\end{align}
\end{subequations}
where
\begin{align*}
&\mH_1: (\vlambda,\vx,\vy)\mapsto(\nicefrac{1}{\beta}\vlambda+\mU\vx-\vx,\vx,\mGamma^{-1}\vy),\\
&\mH_2: (\vlambda,\vx,\vy)\mapsto(\vlambda+\beta\mU\vx-\beta\vy,\vx,\vy),
\end{align*}
and $\beta=\nicefrac{1}{2(n+1)}$ and $\mGamma = \diag\{\alpha_1\mI, \alpha_2\mI,\cdots, \alpha_n\mI\}\in \RR^{nq\times nq}$ denoting the stepsize matrix. The equivalent form of \eqref{eq:AFBS} is given as
\begin{subequations}\label{Alg1}
\begin{align}
&\vlambda^{k+\frac{1}{2}}=\vlambda^k+\beta(\mU\vx^k-\vy^k), \label{Alg1a}\\
&\vx^{k+1}=\prox_{ nr}(\vx^k-\mU\tr\vlambda^{k+\frac{1}{2}}).  \label{Alg1b} \\
&\vy^{k+1}=\vy^{k}-{\mGamma}(\nabla f(\vy^k)-\vlambda^{k+\frac{1}{2}}), \label{Alg1c}\\
&\vlambda^{k+1}=\vlambda^{k+\frac{1}{2}}+\beta[\mU(\vx^{k+1}-\vx^k)-(\vy^{k+1}-\vy^k)].\label{Alg1d}
\end{align}
\end{subequations}
\begin{rmk}
In \cite{Latafat2017}, to solve monotone inclusion problems of the form
$0\in Ax+Mx+Cx$,
where $A$ is a maximally monotone operator, $M$ is a bounded linear operator and $C$ is cocoercive, a general operator splitting scheme AFBS has been given.
\begin{align*}
\text{AFBS: }\left\{
\begin{array}{rll}
\bar{z}^k&=&(H+A)^{-1}(H-M-C)z^k,\\
\alpha_k&=&\frac{\lambda_k\|\bar{z}^k-z^k\|^2_P}{\|(H+M\tr)(\bar{z}^k-z^k)\|^2_{S^{-1}}},\\
z^{k+1}&=&z^k+\alpha_kS^{-1}(H+M\tr)(\bar{z}^k-z^k),
\end{array}\right.
\end{align*}
where $S$ is a strongly positive and self-adjoint operator, $K$ is a skew-adjoint operator, and $H=P+K$, with $P$ being a strongly positive and self-adjoint operator. Inspired by the prediction-correction structure of AFBS, we provide a new splitting method \eqref{eq:AFBS} with a simpler stepsize condition and easier implementation.
\end{rmk}

The iteration \eqref{Alg1} provides us a basic form for designing decentralized algorithm, which will be further elaborated in the next section.

\subsection{Decentralized implementation via relay communication }
Referring to \eqref{Alg1}, we observe that
\begin{align*}
\mU\tr\vlambda^{k+\frac{1}{2}}=\sum_{i=1}^n\vlambda_{i}^{k+\frac{1}{2}},
\end{align*}
 which is required to know the global information on $\vlambda$ at each iteration, is involved in the implement of update \eqref{Alg1b}, which implies that the proposed algorithm \eqref{Alg1} needs $\mathcal{O}(n)$ communications per iteration, and it can not be performed in a decentralized manner. It inspires us that if only one agent $i_k$ is activated and passes information in the $k$-iteration $(k\geq1)$ while the other agents remain unchanged, the coupled term $\mU \tr \vlambda^{k+\frac{1}{2}}$ can be disassembled and computed locally by agents. To be precise, we can reformulate $\mU \tr \vlambda^{k+\frac{1}{2}}$ as
\begin{align*}
 \mU \tr \vlambda^{k} +\vlambda_{i_k}^{k+\frac{1}{2}} -\vlambda_{i_k}^{k},
\end{align*}
 since $\vlambda_{i}^{k+\frac{1}{2}}$ is updated from $\vlambda_{i}^k$ only within agent $i_k$.
Notably, we then introduce a transition variable $\vu^k = \mU \tr \vlambda^{k}$ to record and recursively compute $\mU \tr \vlambda^{k}$, as the alteration of $\vlambda_{i}^{k+1}$ also exclusively occurs in one agent. This leads us to deduce
\begin{align*}
&\mU \tr \vlambda^{k+\frac{1}{2}}=\vu^k+\vlambda_{i_k}^{k+\frac{1}{2}} -\vlambda_{i_k}^{k}\\
&\mU \tr \vlambda^{k+1}=\vu^{k+1}=\vu^{k}+\vlambda_{i_k}^{k+1} -\vlambda_{i_k}^k.
\end{align*}
Such a similar technique can also be found in \cite{Mao2020}, thereby enhancing the facilitation of decentralized implementation.

Furthermore, following this approach, each iteration only needs one information interaction and takes $\mathcal{O}(1)$ communications, which implies that only one agent leaks privacy information and increases the privacy budget per iteration. By appropriately setting the initial iteration condition $\sum_{i=1}^{n}\vlambda_{i}^{1}=0$, \eqref{Alg1} is converted into the relayed and communication-efficient algorithm (RECAL):
\begin{subequations}\label{OurAlg2}
\begin{align}
&\vlambda_{i}^{k+\frac{1}{2}}=
\begin{cases}
\vlambda_{i}^{k}+\beta(\vx^{k}-\vy_i^{k}), ~i=i_{k},  \\
\vlambda_{i}^{k}, ~\text{ o.w.},
\end{cases} \label{OurAlg2a}\\
&\vx^{k+1}=\prox_{nr}\left(\vx^{k}-(\vu^{k}+\vlambda_{i_k}^{k+\frac{1}{2}}-\vlambda_{i_k}^{k})\right), \label{OurAlg2b} \\
&\vy_i^{k+1}=
\begin{cases}
\vy_i^k-\alpha_i(\nabla f_i(\vy_i^k)-\vlambda_i^{k+\frac{1}{2}}), ~i=i_{k}, \\
\vy_{i}^{k}, ~\text { o.w.},
\end{cases}  \label{OurAlg2c}\\
&\vlambda_{i}^{k+1}=
\begin{cases}
\vlambda_{i}^{k+\frac{1}{2}}+\beta[(\vx^{k+1}-\vx^{k})-(\vy_i^{k+1}-\vy_i^k)], ~i=i_{k},  \\
\vlambda_{i}^{k}, ~\text{ o.w.},
\end{cases} \label{OurAlg2d}\\
&\vu^{k+1}=\vu^{k}+\vlambda_{i_{k}}^{k+1}-\vlambda_{i_{k}}^{k}.\label{OurAlg2e}
\end{align}
\end{subequations}

RECAL is named because its implementation can be seen as a ``relay'' with $(\vu,\vx)$ considered as the ``baton''. For example, after completing the computation of the $k$-iteration, the agent $i_k$ will pass its computed baton $(\vu^{k+1},\vx^{k+1})$ with probability $P_{i_k,j}$ to some agent $j$, where $j \in \mathcal{N}_{i_k}$ and $\mathcal{N}_{i_k}$ denotes the set of neighbors of $i_k$, and then the neighbor of $i_k$ who receives the baton will be called $i_{k+1}$. In addition, this process can be considered as a Markov chain, since the baton passing of each iteration is controlled by only one agent and independent of past information.
The transition probability matrix $\mP=[P_{i,j}]_{i,j \in \mathcal{V}}$ of the Markov chain satisfies that $0<P_{ij}<1$ if $j \in \mathcal{N}_{i_k}$; otherwise, $P_{ij}=0$. It follows that $\sum_{j=1}^n P_{ij}=1$.

Given that each agent randomly selects a neighbourhood in its respective invariant set of neighbours and activates its corresponding neighbourhood edges, non-uniform selection implies that each agent needs to cache the diverse sampling probabilities for each neighbourhood edge, which is impractical as described in \cite{Chen2015}. Thus, we can give the following assumption.
\begin{ass}\label{ass:choose probability}
We assume $P_{ij}=\frac{1}{d_i}$, if $(i,j)\in \mathcal{E} $, where $d_i$ is the degree of the agent $i$.
\end{ass}
It implies that the probability of an agent choosing its neighbors to communicate with should be equal and then we have $\mP=\mI-\mD^{-1}\mL$.

\subsection{Discussion}
From the view point of computation, to solve problem \eqref{pro1} in a decentralized way, another alternative problem transformation is
\begin{equation}\label{pro3}
\begin{aligned}
\min ~~f(\vy)+r(\vy),~s.t.~~(\mW\otimes I_q)\vy=0,
\end{aligned}
\end{equation}
where $r(\vy)=\sum_{i=1}^{n}r(\vy_i)$, and $\mW\succeq 0$, if $(i,j)\notin \mathcal{E}$ and $i\neq j$, $\mW_{ij}=\mW_{ji}=0$, and $\mathrm{null}(\mW)=\mathrm{span}(\one_q)$. Based on the transformation, several state-of-the-art decentralized primal-dual proximal algorithms such as IC-ADMM \cite{Chang2015}, PG-ADMM \cite{Aybat2017}, PG-EXTRA \cite{PGEXTR} and NIDS \cite{NDIS} have been proposed. In each iteration of these decentralized algorithms, all agents are involved in communication and computation, which means that each iteration will consume the PLF. Additionally, for the global nonsmooth term $nr(\vx)$, each agent is required to compute $\mathrm{prox}_{ r}(\cdot)$ per iteration, while RECAL implements the consistency process through the ``relay mechanism'' and computes the $\prox_{nr}(\cdot)$ only once for $\vx$ of lower dimensionality, which is economical when there is no analytic solution for the $\prox_{nr}(\cdot)$. Furthermore, different from algorithms \cite{Chang2015,Aybat2017,PGEXTR,NDIS}, whose stepsize condition are coordinated and dependent on the network, RECAL allows for uncoordinated and network-independent stepsizes $\alpha_i\in(0,\frac{2}{L_i+1})$. It implies that multiple communications are not required to achieve global consensus stepsizes. Therefore, the stepsizes of agents can be decided locally, which gives us a freedom in the choice of stepsize.

From the view point of communication, we compare RECAL with other related works. Different from I-ADMM \cite{Ye2020}, in RECAL the agents are activated in a predetermined Markov chain rather than a Hamiltonian circular, which would be suitable for a more general network topology. Different from RW-ADMM \cite{Shah2018}, iterations rely on the passing of the baton and only one agent is activated per iteration, which would remove any ``consensus''-type preprocessing and save communication costs. Different from gossip algorithms \cite{Boyd2006,Hendrikx2019}, the selected edges are continuous and some global coordination to ensure non-overlapping iterations is not necessary.

\subsection{Differential privacy}
\begin{figure}[!t]
\setlength{\abovecaptionskip}{-7pt}
  \centering
  \includegraphics[width=1\linewidth]{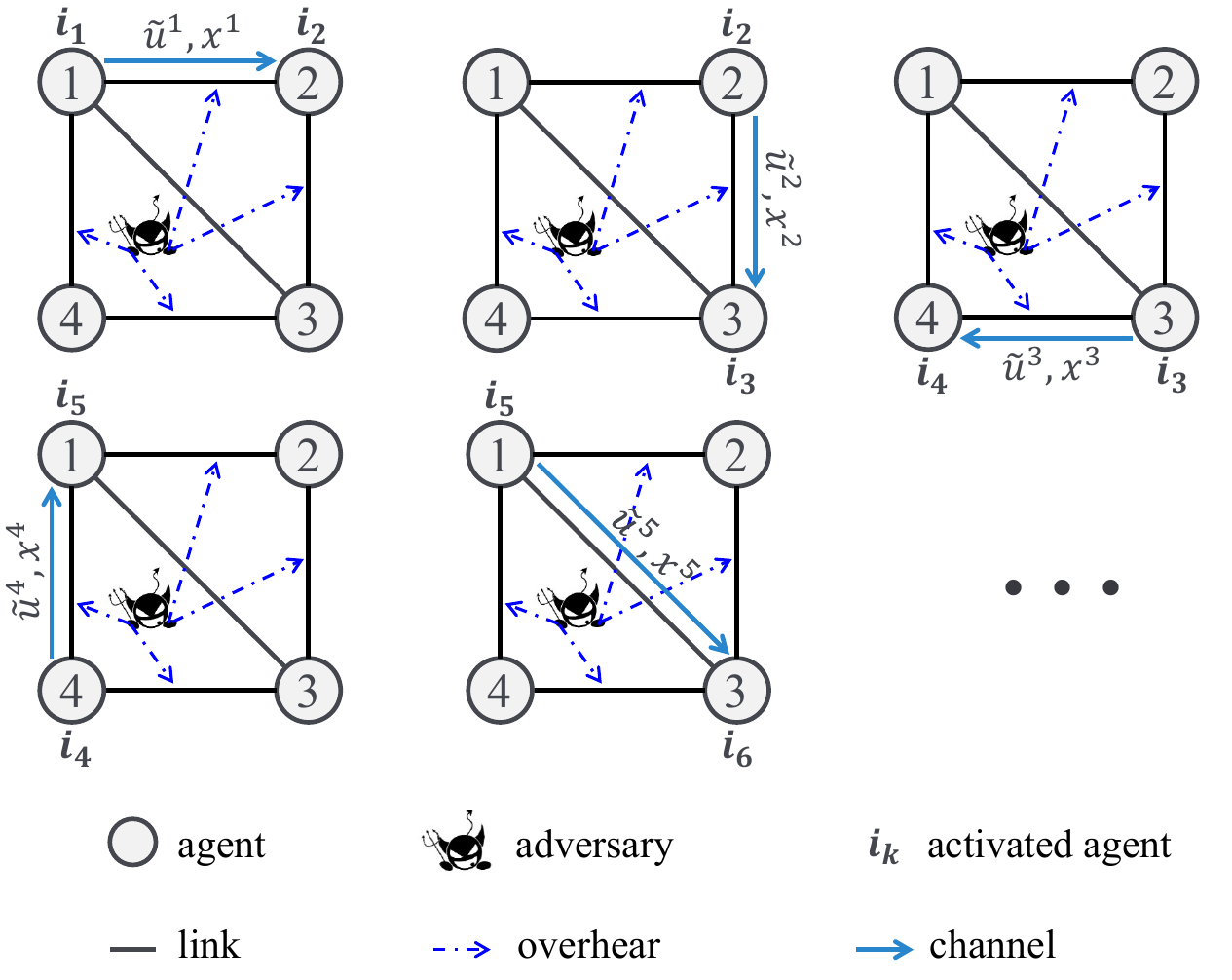}
  \caption{An example of DP-RECAL implementation with external adversaries}
  \label{fig:Alg1}
\end{figure}

\begin{algorithm}[h!]
  \caption{DP-RECAL} 
    \label{Alg2}
  \begin{algorithmic}[1]
    \Require
     Stepsizes $\alpha$, $\beta = \nicefrac{1}{2(n+1)}$, and matrix $\mP$.
    \State Initial $\vx^1$, $\vy^1$ and $\vlambda^{1}=0$.
    \For{$k=1,2,\cdots,\mathrm{K}$}
      \State Agent $i~(i= i_{k})$ computes
      \begin{align*}
      &\vlambda_i^{k+\frac{1}{2}}=\vlambda_{i}^{k}+\beta(\vx^{k}-\vy_i^{k}),\\
      &\vx^{k+1}=\prox_{nr}\left(\vx^{k}-(\tilde{\vu}^{k}+\vlambda_{i}^{k+\frac{1}{2}}-\vlambda_{i}^{k})\right) ,\\
      &\vy_i^{k+1}=\vy_i^k-\alpha_i(\nabla f_i(\vy_i^k)-\vlambda_i^{k+\frac{1}{2}}),\\
      &\vlambda_{i}^{k+1}=\vlambda_{i}^{k+\frac{1}{2}}+\beta[(\vx^{k+1}-\vx^{k})-(\vy_i^{k+1}-\vy_i^k)],\\
      &\vu^{k+1}=\vu^{k}+\vlambda_{i}^{k+1}-\vlambda_{i}^{k},\\
      &\tilde{\vu}^{k+1}=\vu^{k+1}-\vvarepsilon^k.
      \end{align*}
      \State Agent $i~(i\neq i_{k})$ remain the same.
      \State Agent $i_{k}$ pass $\tilde{\vu}^{k+1}$ and $\vx^{k+1}$ to the next chosen agent according to the transition probability
matrix $\mP$.
    \EndFor
    \Ensure
      $\vx^{\mathrm{K}+1}$.
  \end{algorithmic}
\end{algorithm}\
The gradient $\nabla f(\vx^k)$ of Recal can be easily inferred under the threat model of Section \ref{Attack model}, and the detailed procedure is described in the supplement. Such privacy concerns underscore the necessity for a robust privacy mechanism. Given the practical constraints associated with computational and communication demands, an approach that relies on the collaborative efforts of each agent becomes imperative to realize comprehensive privacy protection. Differential privacy mechanisms stand as a fitting choice to be incorporated into RECAL, thus aligning with the broader objective of achieving a collective privacy safeguarding effort.

To customize the proper differential privacy noise, we utilize the local counter $\tau_i(K)$ for $i \in \sV$, which tracks the number of times agent $i$ has been activated (participating in calculations and publishing messages) after $K$ iterations. When agent $i$ is activated in the $k$-iteration, we add the noise $\vvarepsilon^{k}$ drawn from the Gaussian distribution $\mathcal{N}(0,\sigma^2_{\tau_i(k)}\mI)$, where $\sigma^2_{\tau_i(k)}=R\sigma^2_{\tau_i(k)+1}$ and $R$ is the attenuation coefficient of noise, as a perturbation term to the passed $\vu^{k+1}$. Combined with \eqref{OurAlg2}, we finally obtain DP-RECAL. Algorithm \ref{Alg2} and Fig. \ref{fig:Alg1} illustrate the implementation process of DP-RECAL. It is obvious that Gaussian noise $\vvarepsilon^{k}$ is coupled into the update of $\vx^{k+1}$ when $\vu^{k+1}$ is perturbed, so no additional protection is required for the passing of $\vx^{k+1}$.

On the one hand, compared with \cite{Zhang2017}, we introduce $\sigma^2_{\tau_i(k)}$ rather than $\sigma^2_{k}$ to retard the decay of the noise variance, which is essential to reduce the overall privacy budget. See Section \ref{Privacy Analysis} for a detailed analysis. On the other hand, compared with \cite{Huang2020}, we do not prefer to fix the privacy budget of each iteration through decaying stepsizes, which inherently compromises the attainment of linear convergence.

\section{Convergence analysis}\label{Convergence analysis}
In this section, the convergence analysis of DP-RECAL is provided in virtue of the variational inequality context and the fixed point iteration theory.

\subsection{Concise formulation of  DP-RECAL}
In the subsection, we try to simplify the form of DP-RECAL to facilitate the analysis of its convergence.

The existence of both relay mechanism and differential privacy mechanism in DP-RECAL makes it difficult to elucidate the property of the sequence $\{\vw^k\}$. Thus, we firstly let $\hat{\vw}^{k+1}=\col\{\hat{\vx}^{k+1},\hat{\vy}^{k+1},\hat{\vlambda}^{k+1}\}$ denote the data processed by the operator $\mathbb{T}^{[\vvarepsilon^k]}$ from $\vw^{k}$ representing all of the agents are activated, in which
\begin{subequations}\label{Alg3}
\begin{align}
&\hat{\vlambda}^{k+\frac{1}{2}}=\vlambda^{k}+\beta(\mU\vx^k-\vy^k), \label{Alg3a}\\
&\hat{\vx}^{k+1}=\prox_{nr}\left(\vx^k-(\mU\tr\hat{\vlambda}^{k+\frac{1}{2}}-\vvarepsilon^k)\right), \label{Alg3b} \\
&\hat{\vy}^{k+1}=\vy^{k}-\mGamma(\nabla f(\vy^k)-\hat{\vlambda}^{k+\frac{1}{2}}),\label{Alg3c}\\
&\hat{\vlambda}^{k+1}=\hat{\vlambda}^{k+\frac{1}{2}}+\beta[\mU(\hat{\vx}^{k+1}-\vx^k)-(\hat{\vy}^{k+1}-\vy^k)].\label{Alg3d}
\end{align}
\end{subequations}

Then, we give the setup to activating agents for relay communication mechanism in the following definition.
\begin{defn}\label{Active}
(Activating settings)
\begin{enumerate}
  \item The activating indicator vector $\veta^{k}=[\eta^{k}_1,\eta^{k}_2,\cdots,\eta^{k}_n]\tr$ has only one non-zero element $\eta_{i}^k=1$ which implies that only agent $i$ is activated in the $k$-iteration.
  \item The probability of agent $i$ being activated in each iteration is $g_{i}$ defined as $g_i=\mathbf{P}\left(\eta_{i}^{k}=1\right)$, where $i \in \mathcal{V}$.

   \item The activating probability matrix $\tilde{\mPi}$ is defined as $\diag \{g_1,\cdots,g_n\}$ and the activating indicator matrix $\tilde{\mE}^{k}$ is $ \diag\{ \eta_{1}^{k},\cdots,\eta_{n}^{k}\}$.
\end{enumerate}
\end{defn}

In fact, probability vector $\vg=[g_1, g_2,\cdots, g_n]\tr$ can be derived from the transition probability matrix $\mP$. In the following lemma, we give an explicit expression for $\vg$.
\begin{lem}\label{lem:active}
Suppose that Assumption \ref{ass:choose probability} holds. We have
\begin{align*}
\vg=(\mD^{-1}\mL^2\mD^{-1}+\one\one\tr)^{-1}\cdot \one.
\end{align*}
 Moreover, it holds that $g_i>0$, $\forall i \in \mathcal{V}$, which implies the probability of each agent to be activated in the relay mechanism is greater than $0$.
\end{lem}
\begin{IEEEproof}
The proof is provided in Appendix \ref{APP0} in the supplemental materials.
\end{IEEEproof}

Following \cite{Li2022, Latafat2019}, we abbreviate the conditional expectation $\mathbb{E}\left[\cdot \mid \mathcal{F}_{k}\right]$ as $\mathbb{E}_{k}[\cdot]$, where $\mathcal{F}_k$ is the filtration generated by $\{\veta^1,\veta^2,\cdots,\veta^k\}$, and
similarly have two properties for the activating indicator matrix $\tilde{\mE}^{k}$:  $(\tilde{\mE}^{k})^{t}=\tilde{\mE}^{k}$ and $\mathbb{E}_{k}[\tilde{\mE}^{k+1}]=\tilde{\mPi}$. Let
\begin{align*}
\mE^k=\left(
\begin{array}{ccc}
\tilde{\mE}^{k} & 0 & 0 \\
0 & \mI &0\\
0 & 0 & \tilde{\mE}^{k}\\
\end{array}
\right),\quad\mPi=\left(
\begin{array}{ccc}
\tilde{\mPi}& 0 & 0 \\
0 & \mI  &0\\
0 & 0 & \tilde{\mPi}\\
\end{array}
\right),
\end{align*}
 and then we have the following lemma.
\begin{lem} \label{lem2}
 $\mE^{k}$ is an idempotent matrix, i.e., $\left(\mE^{k}\right)^{t}=\mE^{k}$ for any positive integer $t$.
  For the expectation of $\mE^{k+1}$, there exists $\mathbb{E}_{k}\left[\mE^{k+1}\right]=\mPi$.
\end{lem}

Therefore, Lemmas \ref{lem:active} and \ref{lem2} indicate that the relay mechanism can be taken as a special example of the randomized block-coordinate method \cite{Latafat2019,Binchi2016}. Then, we can develop the concise formulation of the proposed DP-RECAL as follows:
\begin{subequations}\label{OP1}
\begin{align}
&\hat{\vw}^{k+1}=\TP^{[\vvarepsilon^k]} \vw^{k}, \label{OP11}\\
&\vw^{k+1}=\vw^{k}+\mE^{k+1}(\hat{\vw}^{k+1}-\vw^{k}), \label{OP12}
\end{align}
\end{subequations}
which is crucial in obtaining Proposition \ref{prop1} in the interest to establish the convergence of DP-RECAL in Theorem \ref{thm1}.

\subsection{Global convergence}
In this subsection, we prove that DP-RECAL is convergent globally by establishing a quasi-Fej\'{e}r monotone sequence. Let
\begin{align*}
&\mH: (\vlambda,\vx,\vy)\mapsto(\frac{1}{\beta}\vlambda,\mU\tr \vlambda+\vx,-\vlambda+\mGamma^{-1}\vy),\\
&\mS:(\vlambda,\vx,\vy)\mapsto (\frac{1}{\beta}\vlambda,\vx,\mGamma^{-1}\vy).
\end{align*}
We have the following result.
\begin{lem} \label{prop1}
Suppose that Assumptions \ref{Ass1} and \ref{ass:choose probability} hold.
For $\forall~\vw^*\in\mathcal{W}^*$, if the stepsizes $\alpha_i\in(0,\frac{2}{L_i+1})$, the sequence $\{\hat{\vw}^{k+1}\}$ generated by \eqref{Alg3} satisfies
\begin{align}\label{Con5}
&\left\|\vw^k-\vw^*\right\|^2_{\mS}-\left\|\hat{\vw}^{k+1}-\vw^*\right\|^2_{\mS}-\left\|\hat{\vw}^{k+1}-\vw^k\right\|^2_{\widehat{\mS}}\nonumber\\
&+\left\langle \hat{\vw}^{k+\frac{1}{2}}-\vw^*,\mD^k \right\rangle\geq0,
\end{align}
where $\mD^{k}=[(\vvarepsilon^k)\tr,0,0]\tr$,  $\vw^{k+\frac{1}{2}}=\col\{\vlambda^{k+\frac{1}{2}},\vx^{k+1},\vy^{k+1}\}$ and
\begin{align*}
&\widehat{\mS} = \left(
      \begin{array}{ccc}
        \frac{1}{\beta}\mI& -\mU&\mI \\
        -\mU\tr&\mI&0\\
        \mI&0&\mGamma^{-1}-\frac{1}{2}\mQ\\
      \end{array}
    \right).
\end{align*}
\end{lem}
\begin{IEEEproof}
The proof is provided in Appendix \ref{APP1} in the supplemental materials.
\end{IEEEproof}
It follows from the proposition that, if ${\vw}^{k+1}=\mathbb{T}^{[0]}\vw^k$, i.e., for $\forall~\vw^*\in\mathcal{W}^*$, $\{\vw^k\}$ generated by \eqref{Alg1} satisfies
\begin{align*}
\left\|\vw^*-\vw^k\right\|^2_{\mS}-\left\|\vw^*-{\vw}^{k+1}\right\|^2_{\mS}-\left\|\vw^k-{\vw}^{k+1}\right\|^2_{\widehat{\mS}}\geq 0,
\end{align*}
which implies that $\{\vw^k\}$ is a Fej\'{e}r monotone sequence with respect to $\mathcal{W}^*$ in $\mS$-norm. Let $\vw^{\infty}$ be an accumulation point of $\{\vw^k\}$ and $\{\vw^{k_j}\}$ be a subsequence converging to $\vw^{\infty}$. It follows from $(\mathbb{T}^{[0]}-\mathbf{I})\vw^k\rightarrow 0$ that $(\mathbb{T}^{[0]}-\mathbf{I})\vw^{k_j}\rightarrow 0$. Since proximal and
linear mappings are continuous and $\nabla f$ is assumed continuous, the operator $\mathbb{T}^{[0]}$ is continuous. Taking $k_j\rightarrow\infty$, we have $(\mathbb{T}^{[0]}-\mathbf{I})\vw^{\infty}=0$, i.e., $\vw^{\infty} \in \text{Fix}\mathbb{T}^{[0]}$. Hence, $\|\vw^k-\vw^{\infty}\|$ monotonically decreases to $0$. In other words, the sequence $\{\vw^k\}$ generated by \eqref{Alg1} converges to the saddle point of the Lagrangian \eqref{Lagrangian}.

Next, based on Proposition \ref{prop1} and \cite{Combettes2015}, we establish the convergence of DP-RECAL from Lemma \ref{Quasi} to Theorem \ref{thm1}.
\begin{lem}\label{Quasi}
Suppose that Assumptions \ref{Ass1} and \ref{ass:choose probability} hold. If the stepsizes $\alpha_i\in(0,\frac{2}{L_i+1})$, there exist $M_1 \geq 0$ and $M_2 \geq 0$ such that $\{\vw^k\}$ generated by DP-RECAL satisfies
\begin{align}\label{ConPF3}
&\mathbb{E}\left[\left\|\vw^{k+1}-\vw^{*}\right\|_{\overline{\mS}}^{2}\right] \leq\left\|\vw^{k}-\vw^{*}\right\|_{\overline{\mS}}^{2}\nonumber\\
&\quad -\left\|\vw^{k}-\hat{\vw}^{k+1}\right\|_{\widehat{\mS}}^{2} +M_1M_2\cdot\|\vvarepsilon^k\|,
\end{align}
where $\overline{\mS}=\mPi^{-1}\mS$ for notational convenience.
\end{lem}
\begin{IEEEproof}
The proof is provided in Appendix \ref{APP2} in the supplemental materials.
\end{IEEEproof}

\begin{thm}\label{thm1}
Suppose that Assumptions \ref{Ass1} and \ref{ass:choose probability} hold. If the stepsizes $\alpha_i\in(0,\frac{2}{L_i+1})$, there exists $\vw^{\infty}\in \mathcal{W}^*$ such that
$$
\lim_{k\rightarrow\infty} \vw^{k+1}=\vw^{\infty}, \text{almost surely},
$$
i.e., $\vy^{\infty} = \one \otimes \vx^{\infty}$ and $\vx^{\infty}$ solves problem \eqref{pro1}.

\end{thm}
\begin{IEEEproof}
From Lemma \ref{Quasi} and \cite[Proposition 2.3]{Combettes2015}, it can be obtained.
\end{IEEEproof}
Due to the exiting differential privacy mechanism, extra difficulties have been imposed on the establishment of the convergence of DP-RECAL. We adopt a novel analysis method that employs ``transition" updates $\overline{\vw}^{k+1}=\mathbb{T}^{[\vvarepsilon^k=0]}\vw^k$ and $\breve{\vw}^{k+1}=\vw^{k}+\mE^{k+1}(\overline{\vw}^{k+1}-\vw^{k})$ as a reference in the $k$-iteration, which proves the boundedness of the sequence $\|\hat{\vw}^{k+1}-\vw^*\|$. In this way, \eqref{ConPF3} holds and $\{\vw^k\}$ generated by DP-RECAL is a quasi-Fej\'{e}r monotone sequence with respect to $\mathcal{W}^*$ in $\overline{\mS}$-norm, rather than a Fej\'{e}r monotone sequence in \cite{Li2022,Latafat2019,Guo2022}.
\subsection{Linear convergence rate}
In this subsection, we establish the global linear convergence of DP-RECAL under metric subregularity \cite{va2004}, i.e.,
\begin{defn}
A set-valued operator $\Psi:\mathbb{R}^{n}\rightrightarrows \mathbb{R}^{m}$ is metrically subregular at $(\bar{\vu},\bar{\vv})\in \mathrm{gph}(\Psi)$ if for some $\epsilon>0$ there exists $\tilde{\varsigma}\geq0$ such that
\begin{align*}
\dist(\vu,\Psi^{-1}(\bar{\vv}))\leq \tilde{\varsigma}\cdot \dist(\bar{\vv},\Psi(\vu)),\quad \forall \vu\in \mathcal{B}_{\epsilon}(\bar{\vu}),
\end{align*}
where $\mathrm{gph}(\Psi):=\{(\vu,\vv):\vv=\Psi(\vu)\}$, $\Psi^{-1}(\vv):=\{\vu\in\mathbb{R}^n:(\vu,\vv)\in\mathrm{gph}(\Psi)\}$ and $\dist_{\mA}(\vx,\mathcal{X}):=\|\vx-\mathcal{P}_{\mathcal{X}}^{\mA}(\vx)\|_{\mA}$.
\end{defn}

 From the conclusion and analysis in Theorem \ref{thm1}, we give the following theorem to show the linear convergence rate of the proposed algorithm.
\begin{thm}\label{thm2}
Suppose that Assumptions \ref{Ass1} and \ref{ass:choose probability} hold. If the stepsizes $\alpha_i\in(0,\frac{2}{L_i+1})$, and $\partial\mathcal{L}$ is metrically subregular at $(\vw^{\infty},0)$ with modulus $\tilde{\varsigma}$, there exits $\varsigma>0$ such that
\begin{align*}
&\EE_k\left[\dist^2_{\overline{\mS}}(\vw^{k+1},\mathcal{W^*})\right]\leq\frac{\phi}{\phi+1}\left[\dist^2_{\overline{\mS}}(\vw^{k},\mathcal{W^*})+V\big\|\vvarepsilon^k\big\|^2\right],
\end{align*}
where $\phi=(\varsigma+1)^2\|\overline{\mS}\|\|(\mH+\mA\tr)^{-1}\mS\|^2\|\widehat{\mS}^{-1}\|$ and $V=M_1M_2+\frac{\phi}{\phi+1}\mu^2$.

Let $\ell$ denote the maximum number of iterations between the two activations of any agent, i.e., for $\forall i \in \mathcal{V}$, $\tau_i(k+\ell)-\tau_i(k)\geq 1$. If $(\frac{\phi}{\phi+1})^\ell R>1$, the sequence $\{\vw^k\}$ generated by DP-RECAL satisfies
\begin{align}\label{Linear Eq}
&\mathbb{E}\left[\dist^2_{\overline{\mS}}(\vw^{k+1},\mathcal{W^*})\right]\nonumber\\
&\leq(\frac{\phi}{\phi+1})^k\big(\dist^2_{\overline{\mS}}(\vw^{1},\mathcal{W^*})+R\frac{(\frac{\phi}{\phi+1})^\ell V\sigma^2_1} {(\frac{\phi}{\phi+1})^\ell R-1}\big),
\end{align}
which implies that
$\{\vw^k\}$ converges to $\mathcal{W}^*$ R-linearly.
\end{thm}
\begin{IEEEproof}
The proof is provided in Appendix \ref{APP3} in the supplemental materials.
\end{IEEEproof}
Compared with optimization algorithms like  \cite{Shah2018} and \cite{Latafat2019} who need the global subregularity assumption or the strongly convex assumption, we apply a more general condition that $\partial\mathcal{L}$ is metrically subregular at $(\vw^{\infty},0)$. Some commonly used but not strongly convex loss functions in machine learning satisfying this metrically subregular condition. Some examples  from \cite{Xiaoming Yuan2020} can be seen.

Compared with several decentralized primal-dual proximal algorithms, the global linear convergence of DP-RECAL is established with a more liberal choice of parameters. In addition, we can also optimize the execution path of the algorithm in practice to reduce the overall communication distance and further improve the performance of our proposed algorithm. Algorithms such as Dijkstra may help. Specific numerical examples will show the advantages of DP-RECAL from various aspects in Section \ref{Performance evaluation}.
\section{Privacy Analysis}\label{Privacy Analysis}
In this section, we examine the privacy-preserving performance of DP-RECAL, i.e., the overall privacy budget $\epsilon$.

Zero concentrated differential privacy (zCDP) as a relaxed version of $(\epsilon,\delta)$-differential privacy is introduced as our specific privacy mechanism. It is proposed in \cite{Bun2016} aiming to establish better accuracy and lower privacy loss bounds compared with dynamic differential privacy \cite{Zhang2017}. The definition of zCDP is as follows:
\begin{defn} \cite{Bun2016}. \label{ZDP}
A randomised algorithm $\mathcal{M}$ is $\rho$-zCDP, if we have
\begin{align*}
\EE\left[e^{(\kappa-1)Z}\right] \leq e^{(\kappa-1)\kappa\rho}, \quad\forall \kappa \in (1,+\infty)
\end{align*}
for any two neighbouring datasets $D$ and $\hat{D}$ differing in at most one record, where $Z$ is the loss random variable of $\mathcal{M}$.
\end{defn}
Moreover, we illustrate the correlation between $\rho$-zCDP and $(\epsilon,\delta)$-differential privacy which can be illustrated in the following proposition so to describe the overall privacy budget $\epsilon$ of DP-RECAL.
\begin{prop}\label{z-DP}\cite{Bun2016}.
If an algorithm $\mathcal{M}$ is proven to be $\rho$-zCDP, then $\mathcal{M}$ satisfies $(\rho+2\sqrt{\rho  ln(1/\delta)},\delta)$-differential privacy for any $\delta \in (0,1)$.
\end{prop}

To carry out our analysis, we firstly focus on the private protection of the $k$-iteration occurring at activated agent $i$ and define $\rho_{\tau_i(k)}$-zCDP.
\begin{defn}\label{iZDP}
$\forall \kappa \in (1,+\infty)$, a randomised algorithm $\mathcal{M}_{\tau_i(k)}$ at activated agent $i$ in the $k$-iteration is $\rho_{\tau_i(k)}$-zCDP, if we have
$\EE\left[e^{(\kappa-1)Z_{\tau_i(k)}}\right] \leq e^{(\kappa-1)\kappa\rho_{\tau_i(k)}}$
for any two neighbouring datasets $D_{i}$ and $\hat{D}_{i}$ differing in at most one record, where $Z^{\tau_i(k)}$ is the privacy loss random variable of $\mathcal{M}_{\tau_i(k)}$.
\end{defn}
We introduce the following assumption widely used in the literature such as \cite{Zhang2017} and \cite{Huang2020}.
\begin{ass}\label{bound gradient}
For $\forall~i\in \mathcal{V}$, there exists a constant $c$ such that $\|\nabla f_i(\cdot)\| \leq c$.
\end{ass}
Subsequently, we estimate the sensitivity of DP-RECAL.
\begin{lem} \label{lem:sensitive}
Under Assumption \ref{bound gradient}, the sensitivity of DP-RECAL, denoted by $\Delta$, satisfies:
\begin{align*}
\Delta=\max_{D_{i}, \hat{D}_{i}}\|\vu^{k+1}_{D_{i}}-\vu^{k+1}_{\hat{D}_{i}}\|=4\alpha\beta c,
\end{align*}
where $\alpha=\max_{i \in \mathcal{V}}\alpha_i$.
\end{lem}
\begin{IEEEproof}
According to Algorithm \ref{Alg2}, we have
$\vu^{k+1}_{D_{i}}=\vu^{k}+\vlambda_{i,D_{i}} ^{k+1}-\vlambda_{i}^k$, and $\vu_{\hat{D}_{i}}^{k+1}=\vu^k+\vlambda_{i,\hat{D}_{i}}^{k+1}-\vlambda_{i}^k$.
Then, the sensitivity $\Delta$ is given by
\begin{align*}
 &\|\vu_{D_{i}}^{k+1}-\vu_{\hat{D}_{i}}^{k+1}\| =\|(\vlambda_{i,D_{i}}^{k+1}-\vlambda_{i}^k)-(\vlambda_{i,\hat{D}_{i}}^{k+1}-\vlambda_{i}^k)\| \\
 &=\beta\|-2\vy_{i,D_{i}}^{k+1}+2\vy_{i,\hat{D}_{i}}^{k+1}\|\leq2\beta \alpha\|\nabla f_{i,D_{i}}(\vy_i^k)-\nabla f_{i,\hat{D}_{i}}(\vy_i^k)\|\\
 &\leq 4\alpha\beta c.
\end{align*}
Thus, the sensitivity $\Delta$ equals to $4\alpha\beta  c$.
\end{IEEEproof}
The sensitivity serves to quantify the degree to which minor perturbations in the private dataset can lead to alterations in the algorithm's output within a single iteration. Essentially, higher sensitivity implies a higher likelihood for an attacker to deduce the private dataset through the algorithm's output changes. To counter this, we introduce differential privacy noise, which serves to obscure these output changes, and compute the exact privacy budget for agent $i$ under the introduced noise.
\begin{lem}\label{itZDP}
The agent $i$ DP-RECAL activated in the $k$-iteration satisfies $\rho_{\tau_i(k)}$-zCDP, where $\rho_{\tau_i(k)}=8\alpha^2\beta^2 c^2 /\sigma^2_{\tau_i(k)}$.
\end{lem}
\begin{IEEEproof}
The privacy loss random variable $Z_{\tau_i(k)}$ from $\vu^{k+1}$ among two neighbouring datasets $D_{i}$ and $\hat{D}_{i}$ is calculated by:
\begin{align}\label{PL}
Z_{\tau_i(k)}=\left| \ln\frac{\mathbf{P}(\vu^{k+1}|D_{i})}{\mathbf{P}(\vu^{k+1}|\hat{D}_{i})}\right|.
\end{align}
Since the probability distribution of $\hat{\vu}^{k+1}_{D_{i}}$  and $\hat{\vu}^{k+1}_{\hat{D}_{i}}$ are $\mathcal{N}\big(\vu_{D_{i}}^{k+1},\sigma^2_{\tau_i(k)}\mI\big)$ and $\mathcal{N}\big(\vu_{\hat{D}_{i}}^{k+1},\sigma^2_{\tau_i(k)}\mI\big)$ respectively, we follow \cite[Lemma 2.5]{Bun2016} and for $\forall \kappa \in (1,+\infty)$ rewrite \eqref{PL} as
\begin{align}\label{PL1}
Z_{\tau_i(k)}=\frac{\kappa\big\|\vu_{D_{i}}^{k+1}-\vu_{\hat{D}_{i}}^{k+1}\big\|^2}{2\sigma^2_{\tau_i(k)}}.
\end{align}
Then, applying \eqref{PL1} to Definition \ref{ZDP}, it holds that
\begin{align}
&\EE\left[e^{(\kappa-1)Z_{\tau_i(k)}}\right] \leq e^{(\kappa-1)\kappa \Delta /2\sigma^2_{\tau_i(k)}}\nonumber\\
&\leq e^{8(\kappa-1)\kappa\alpha^2\beta^2 c^2 /\sigma^2_{\tau_i(k)}}=e^{(\kappa-1)\kappa\rho_{\tau_i(k)}}.
\end{align}
Therefore, the agent $i$ activated in the $k$-iteration is shown to be $\rho_{\tau_i(k)}$-zCDP with $\rho_{\tau_i(k)}=8\alpha^2\beta^2 c^2 /\sigma^2_{\tau_i(k)}$.
\end{IEEEproof}
Finally, according to Proposition \ref{sumCDP}-\ref{sumCDP1} and \ref{z-DP}, we can obtain the overall privacy budget for DP-RECAL.
\begin{thm}\label{thm:privacy budget}
Under Assumption \ref{bound gradient}, Algorithm \ref{Alg2} achieves $(\epsilon, \delta)$-differential privacy after $K$ iterations, where
\begin{align*}
\epsilon = \frac{\rho_1(R^\xi - 1)}{R - 1} + 2\sqrt{\frac{\rho_1(R^\xi - 1)\ln(1/\delta)}{R - 1}}.
\end{align*}
Here, $\rho_1$ is $\frac{8\alpha^2\beta^2 c^2}{\sigma_1^2}$, and $\xi$ denotes the PLF of DP-RECAL.
\end{thm}
\begin{IEEEproof}
After $K$ iterations, suppose that there exits an agent $j$ satisfying $\tau_K(j)=\xi=\max_{i \in\mathcal{V}} \tau_i(K)$, which implies that the agent $j$ consumes the most privacy budget among agents. By Theorem \ref{LCI-pb}, DP-RECAL provides $\sum_{t=1}^{\tau_j(K)}\rho_{t}$-zCDP. Note that the noise variance is set as $\sigma^2_{t}=R\sigma^2_{t+1}$ and from Lemma \ref{itZDP}, it holds that $\rho_{t+1}=R\rho_{t}$. Thus, we have $\sum_{t=1}^{\xi}\rho_t=\frac{\rho_1(R^{\xi}-1)}{R-1}$, where $\xi$ is PLF. By Proposition \ref{z-DP}, DP-RECAL satisfies $(\epsilon,\delta)$-differential privacy, with
\begin{align*}
\epsilon=\frac{\rho_1(R^{\xi}-1)}{R-1} +2\sqrt\frac{\rho_1(R^{\xi}-1)\ln 1/\delta}{R-1},
\end{align*}
where $\rho^1=\frac{8\alpha^2\beta^2 c^2 }{\sigma^2_{1}}$, $R>1$ and $\delta\in(0,1)$.
\end{IEEEproof}
\begin{rmk}
Theorem \ref{thm:privacy budget} reveals the quantitative relationship between the privacy budget and PLF for DP-RECAL. When the PLF increases, the resultant overall privacy cost exhibits exponential growth. We observe that the choice of stepsize also has an effect on the privacy budget, but owing to our utilization of the maximum stepsize among agents, the result we derive does not represent a strict upper limit on the overall privacy budget. For individual agents with special security needs, it is feasible to adjust the local privacy budget to the desired level due to its independent ability to set local stepsizes.
\end{rmk}
From the above discussion, the privacy bound of DP-RECAL is manageable, i.e., DP-RECAL can be considered as privacy-preserving. Moreover, we can get that if the PLF of DP-RECAL is less than the number of iterations of other differentially private algorithms when converging to a fixed accuracy, then DP-RECAL consumes less privacy budget $\epsilon$ compared with existing algorithms. It is clear from \cite{Dwork2014} that DP-RECAL owns better privacy performance.

\section{Communication and privacy budget}\label{Communication and privacy budget}
In this section, we calculate the iteration numbers when RECAL converges to a fixed accuracy for the decentralize least squares problem as \eqref{pro:least squares}, and then compare the computation and communication complexity, and PLF between RECAL and other related algorithms.
\begin{equation}\label{pro:least squares}
\begin{aligned}
\min ~\frac{1}{n}\sum_{i=1}^{n}\frac{1}{2}\|\mB_i\vy_i-\vb_i\|^2,\quad
s.t.~\mU\vx-\vy=0.
\end{aligned}
\end{equation}
By Theorem \ref{thm2}, we establish the computation complexity of DP-RECAL as
\begin{align*}
\mathbb{E}\left[\mathrm{dist}_{\mS}(\vw^{k+1},\mathcal{W}^*)^2\right]\leq (\frac{\phi}{\phi+1})^{k}\mathrm{dist}_{\mS}(\vw^1,\mathcal{W}^*)^2\leq \epsilon,
\end{align*}
which implies
\begin{align}\label{itnumber}
k_{\epsilon}=\ln\left(\frac{\|\vw^1-\vw^*\|^2_{\mS}}{\epsilon}\right)/\ln\left(1+\frac{1}{\phi}\right),
\end{align}
where $\phi=(\varsigma+1)^2\|\overline{\mS}\|\|(\mH+\mA\tr)^{-1}\mS\|^2\|\widehat{\mS}\|^{-1}$.

 For the least squares problem, $\partial \mathcal{L}$ is global metric subregular, i.e., for any $\overline{\vw}^{k+\frac{1}{2}}$, $\mathrm{dist}(\overline{\vw}^{k+\frac{1}{2}},\mS)\leq \tilde{\varsigma}\left(\|\mQ\|+\|\mH-\mA\|\right) \|\overline{\vw}^{k+\frac{1}{2}}-\vw^{k}\|$. From \cite[Theorem 45 and Lemma 52]{Yuan2020}, we have $\tilde{\varsigma}=\mathcal{O}(n)$. In addition, combined with the fact that $\mQ$ and $\mH-\mA$ are the diagonal matrix and the upper triangular matrix, respectively, $\varsigma$ is presented to be of the same order as $\tilde{\varsigma}$, i.e., $\varsigma=\mathcal{O}(n)$. Similarly, $\|\overline{\mS}\|$ and $\|(\mH+\mA\tr)^{-1}\mS\|$ are calculated to be constants independent from the network scale $n$. As for $\|\widehat{\mS}\|$, with $\mL=(-\mU, \mI)$, the ``LDL'' decomposition of $\widehat{\mS}$ reveals that
\begin{align*}
\widehat{\mS}=\underbrace{\left[
\begin{array}{cc}
\mI & 0  \\
\beta\mL\tr& \mI \\
\end{array}
\right]}_{\mC\tr}\underbrace{\left[
\begin{array}{cc}
\frac{1}{\beta}\mI& 0  \\
0& \tilde{\mGamma}^{-1}-\frac{\tilde{\mQ}}{2}-\beta\mL\tr\mL\\
\end{array}
\right]}_{\mO}\underbrace{\left[
\begin{array}{cc}
\mI & \beta\mL \\
0& \mI \\
\end{array}
\right]}_{\mC},
\end{align*}
where $\tilde{\mGamma}=\diag\{\mI,\mGamma\}$ and $\tilde{\mQ}=\diag\{0,\mQ\}$.
It is obvious that $\|\mC\|=\|\mC\tr\|=1$ and $\|\mO\|=\mathcal{O}(n)$, and thus $\|\widehat{\mS}\|=\mathcal{O}(n)$. Therefore, it can be estimated that $\phi=\mathcal{O}(n)$, which is consistent with the actual situation. Since $\ln(1+\frac{1}{n}) \approx \frac{1}{n}$ for $\frac{1}{n}$ close to $0$, \eqref{itnumber} can be simplified to
$
k_{\epsilon}\sim\mathcal{O}\left(n\ln\nicefrac{1}{\epsilon}\right).
$

We start by comparing the computational complexity, communication complexity, and PLF of RECAL with other exsiting decentralized optimization algorithms like ADMM \cite{Shi2014}, EXTRA \cite{PGEXTR}, NIDS \cite{NDIS}, and RWADMM \cite{Shah2018} on two specific graph scenarios. Notably, differentially private algorithms \cite{Zhang2017}, \cite{Huang2020} and \cite{Zhang2018} are unable to precisely characterize their own computational and communication performance due to the presence of the differential privacy noise and decaying stepsizes. In light of this, we refer to the work in ADMM\cite{Shi2014} since the iterative sequences employed in \cite{Zhang2017} and \cite{Zhang2018} are rooted in the framework presented in \cite{Shi2014}.

\textbf{Example 1 (Complete graph)} In a complete graph, every agent connects with all the other agents. The number of edges $m=\mathcal{O}(n^2)$ and the degree of the agent $d_i=n$. Refer to Table \ref{Table2} for specific comparisons, where Comp. and Comm. stand for computation complexity and communication complexity.
\begin{table}[!htb]
\renewcommand\arraystretch{1.7}
\begin{center}
\caption{Computation, Communication and PLF Complexity of Related Algorithms in Complete Graph}
\scalebox{0.9}{
\begin{tabular}{c c c c}
\hline
\textbf{Algorithm} &\textbf{Comp.}&\textbf{PLF}&\textbf{Comm.}\\
\hline
ADMM \cite{Shi2014} &$\mathcal{O}\Big(\ln\nicefrac{1}{\epsilon}\Big)$ &$\mathcal{O}\Big(\ln\nicefrac{1}{\epsilon}\Big)$ & $\mathcal{O}\Big(n^2\ln\nicefrac{1}{\epsilon}\Big)$\\
EXTRA \cite{PGEXTR}& $\mathcal{O}\Big(\ln\nicefrac{1}{\epsilon}\Big)$&$\mathcal{O}\Big(\ln\nicefrac{1}{\epsilon}\Big)$ &$\mathcal{O}\Big(n^2\ln\nicefrac{1}{\epsilon}\Big)$\\
NIDS \cite{NDIS}& $\mathcal{O}\Big(\ln\nicefrac{1}{\epsilon}\Big)$&$\mathcal{O}\Big(\ln\nicefrac{1}{\epsilon}\Big)$& $\mathcal{O}\Big(n^2\ln\nicefrac{1}{\epsilon}\Big)$\\
RWADMM \cite{Shah2018}& $\mathcal{O}\Big(n^2\ln\nicefrac{1}{\epsilon}\Big)$&$\mathcal{O}\Big(n^2\ln\nicefrac{1}{\epsilon}\Big)$ &$\mathcal{O}\Big(n^3\ln\nicefrac{1}{\epsilon}\Big)$\\
\rowcolor{mygray}\textbf{RECAL}&$\mathcal{O}\Big(n\ln\nicefrac{1}{\epsilon}\Big)$&$\mathcal{O}\Big(\ln\nicefrac{1}{\epsilon}\Big)$ &$\mathcal{O}\Big(n\ln\nicefrac{1}{\epsilon}\Big)$\\
\hline
\end{tabular}}
\label{Table2}
\end{center}
\end{table}

\textbf{Example 2 (Circular graph)} Consider a circular, where each agent connects with its neighbors on both sides. We can verify that $m=n$, $d=2$ and $1-\sigma_m(\mP)=\mathcal{O}(\nicefrac{1}{n^2})$, where $\sigma_m(\cdot)$ represents the minimum nonzero eigenvalue. Refer to Table \ref{Table3} for specific comparisons.
\begin{table}[!htb]
\renewcommand\arraystretch{1.7}
\begin{center}
\caption{Computation, Communication and PLF Complexity of Related Algorithms in Circular Graph}
\scalebox{0.9}{
\begin{tabular}{cccc}
\hline
\textbf{Algorithm} &\textbf{Comp.}&\textbf{PLF}&\textbf{Comm.}\\
\hline
ADMM \cite{Shi2014}&$\mathcal{O}\Big(n^2\ln\nicefrac{1}{\epsilon}\Big)$ &$\mathcal{O}\Big(n^2\ln\nicefrac{1}{\epsilon}\Big)$ &$\mathcal{O}\Big(n^3\ln\nicefrac{1}{\epsilon}\Big)$\\
EXTRA \cite{PGEXTR}& $\mathcal{O}\Big(n^2\ln\nicefrac{1}{\epsilon}\Big)$&$\mathcal{O}\Big(n^2\ln\nicefrac{1}{\epsilon}\Big)$ &$\mathcal{O}\Big(n^3\ln\nicefrac{1}{\epsilon}\Big)$\\
NIDS \cite{NDIS}& $\mathcal{O}\Big(n^2\ln\nicefrac{1}{\epsilon}\Big)$&$\mathcal{O}\Big(n^2\ln\nicefrac{1}{\epsilon}\Big)$ &$\mathcal{O}\Big(n^3\ln\nicefrac{1}{\epsilon}\Big)$\\
RWADMM \cite{Shah2018}& $\mathcal{O}\Big(n^3\ln\nicefrac{1}{\epsilon}\Big)$&$\mathcal{O}\Big(n^2\ln\nicefrac{1}{\epsilon}\Big)$ &$\mathcal{O}\Big(n^3\ln(\nicefrac{1}{\epsilon})\Big)$\\
\rowcolor{mygray}\textbf{RECAL}&$\mathcal{O}\Big(n\ln\nicefrac{1}{\epsilon}\Big)$&$\mathcal{O}\Big(\ln\nicefrac{1}{\epsilon}\Big)$ &$\mathcal{O}\Big(n\ln\nicefrac{1}{\epsilon}\Big)$\\
\hline
\end{tabular}}
\label{Table3}
\end{center}
\end{table}

It is clear that RECAL imposes less communication burden on the optimization model than other decentralized algorithms listed in both Table \ref{Table2} and \ref{Table3}.
As to the complete graph, Table \ref{Table2} reveals that the PLF experienced by RECAL when converging to a fixed precision is less than RW-ADMM and similar as D-ADMM, EXTRA and Exact diffusion. It implies that, when other algorithms in Table \ref{Table2} introduce differential privacy mechanism, the privacy budget $\epsilon$ of DP-RECAL is lower than that of DP-RW-ADMM, and is similar as that of the remaining privacy-preserving algorithms. Likewise, concerning the circular graph, as shown in Table \ref{Table3}, the privacy budget $\epsilon$ of DP-RECAL is lower than that of all other  algorithms.

Next, we consider the general decentralized network topology, where $m\in [n,n(n-1)] $, $g_i \sim \mathcal{O}(\frac{1}{n}) $ and $1-\sigma_m(\mP)$ measures the connectivity of the network. Let $\chi=1-\sigma_m(\mP)$.
Refer to Table \ref{Table4} for comparisons.
\begin{table}[!htb]
\renewcommand\arraystretch{1.7}
\begin{center}
\caption{Computation, Communication and PLF Complexity of Related Algorithms in General Graph}
\scalebox{0.9}{
\begin{tabular}{cccc}
\hline
\textbf{Algorithm} &\textbf{Comp.}&\textbf{PLF}&\textbf{Comm.}\\
\hline
ADMM \cite{Shi2014}&$\mathcal{O}\Big(\nicefrac{1}{\chi}\ln\nicefrac{1}{\epsilon}\Big)$ &$\mathcal{O}\Big(\nicefrac{1}{\chi}\ln\nicefrac{1}{\epsilon}\Big)$ &$\mathcal{O}\Big(\nicefrac{m}{\chi}\ln\nicefrac{1}{\epsilon}\Big)$\\
EXTRA \cite{PGEXTR}& $\mathcal{O}\Big(\nicefrac{1}{\chi}\ln\nicefrac{1}{\epsilon}\Big)$ &$\mathcal{O}\Big(\nicefrac{1}{\chi}\ln\nicefrac{1}{\epsilon}\Big)$ &$\mathcal{O}\Big(\nicefrac{m}{\chi}\ln\nicefrac{1}{\epsilon}\Big)$\\
NIDS \cite{NDIS}& $\mathcal{O}\Big(\nicefrac{1}{\chi}\ln\nicefrac{1}{\epsilon}\Big)$&$\mathcal{O}\Big(\nicefrac{1}{\chi}\ln\nicefrac{1}{\epsilon} \Big)$ &$\mathcal{O}\Big( \nicefrac{m}{\chi}\ln\nicefrac{1}{\epsilon}\Big)$\\
RWADMM \cite{Shah2018}& $\mathcal{O}\Big(\nicefrac{m}{\chi}\ln\nicefrac{1}{\epsilon}\Big)$& $\mathcal{O}\Big(\nicefrac{m^2}{n^2\chi}\ln\nicefrac{1}{\epsilon}\Big)$& $\mathcal{O}\Big(\nicefrac{m^2}{n\chi}\ln\nicefrac{1}{\epsilon}\Big)$\\
\rowcolor{mygray}\textbf{RECAL}&$\mathcal{O}\Big(n\ln\nicefrac{1}{\epsilon} \Big)$&$\mathcal{O}\Big(\ln\nicefrac{1}{\epsilon}\Big)$& $\mathcal{O}\Big(n\ln\nicefrac{1}{\epsilon}\Big)$\\
\hline
\end{tabular}}
\label{Table4}
\end{center}
\end{table}

In Table \ref{Table4}, we can confirm that in whatever network topology, DP-RECAL exhibits superior communication efficiency and better security than any compared algorithms.
\section{Numerical experiments}\label{Performance evaluation}
In this section, we conduct two kinds of numerical experiments to validate the efficiency of the proposed algorithm. All the algorithms are performed on Matlab R2022a on a computer with 2.50 GHz 12th Gen Intel(R) Core(TM) i5-12500H and 16 GB memory.

\subsection{Decentralized classification problems}
We experimentally evaluate our proposed algorithm DP-RECAL by solving a decentralized linear regression and a decentralized logistic regression  as classic machine learning models based on real-world datasets\footnote{Some experimental details are given in the supplement.}.

We consider the logistic regression and linear regression problems on a circular graph of 8 agents. In view of the space constraint, the experimental results of the linear regression are shown below, and rest of results are shown in the supplement.
\begin{align*}
&\textbf{Linear Regression:}\\
&\min _{\vx}\frac{1}{n}\sum_{i=1}^{n}\frac{1}{m_i}\sum_{j=1}^{m_{i}} \frac{1}{2}\|\mB_{ij}\vx-\vb_{ij}\|^2+\frac{1}{2}\|\vx\|^2+\frac{1}{2}\|\vx\|_1,
\end{align*}
where $(\mB_{ij},\vb_{ij})$ are training data samples of agent $i$, and $\mB_{ij} \in \RR^{q}$ is a feature vector and the corresponding label is $\vb_{ij} \in \{+1,-1\}$ for $j=1,\cdots,m_i$.

We use four types of datasets: ijcnn1 datasets owning 49990 samples with 22 features; `0' and `1' datasets in MNIST owning 12665 samples with 784 features; `T-shirt' and `trouser' datasets in Fashion-MNIST owning 12000 samples with 784 features; `cat' and `dog' datasets in CIFAR-10 owning 10000 samples with 3072 features.

The following steps should be applied to preprocess the datasets. We firstly remove missing values and then normalize all numerical features so that the range of each value is $[0,1]$. Moreover, we transform labels by the one-hot encoding method into $\{+1,-1\}$. Ultimately, these processed data will be allocated equally to each local agent.

\begin{figure*}[!t]
\centering
\setlength{\abovecaptionskip}{-7pt}
\subfigure{
\includegraphics[width=0.25\linewidth]{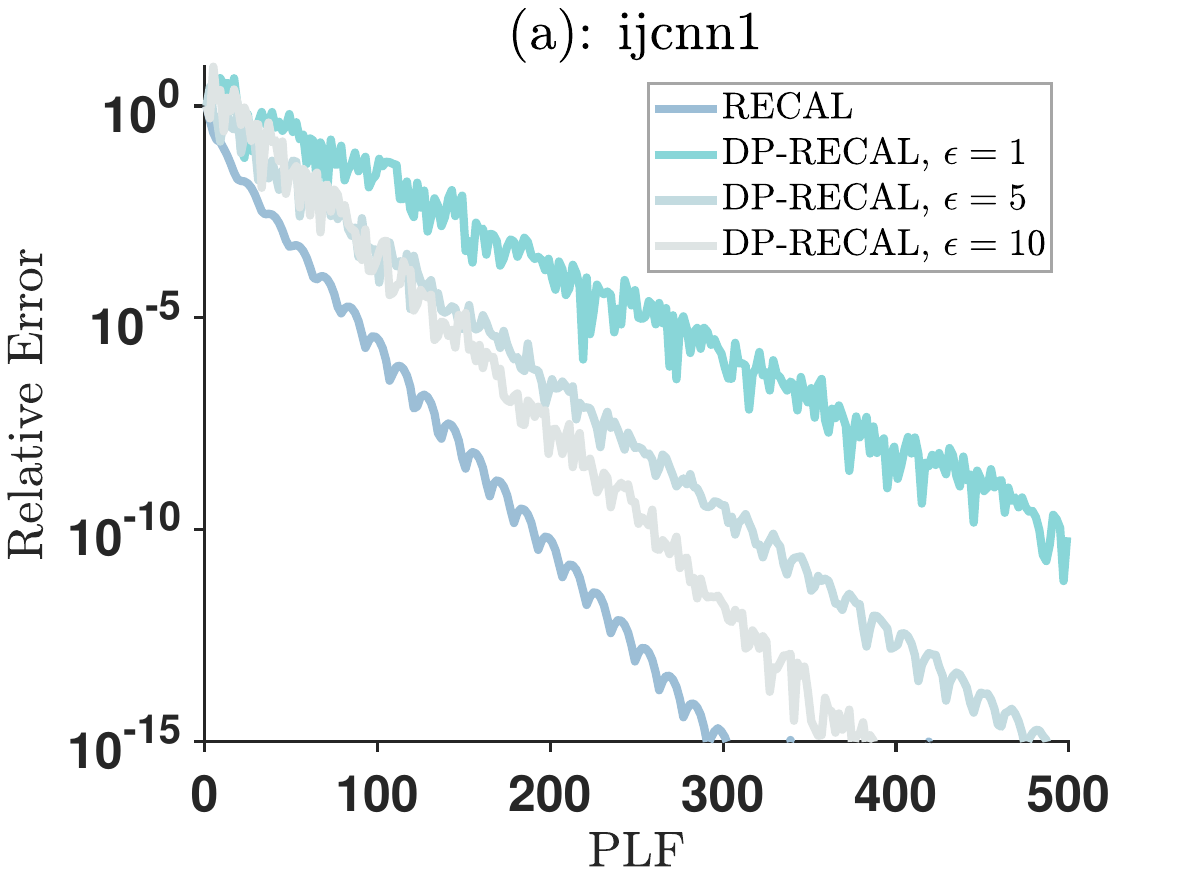}}\hspace{-2.85mm}
\subfigure{
\includegraphics[width=0.25\linewidth]{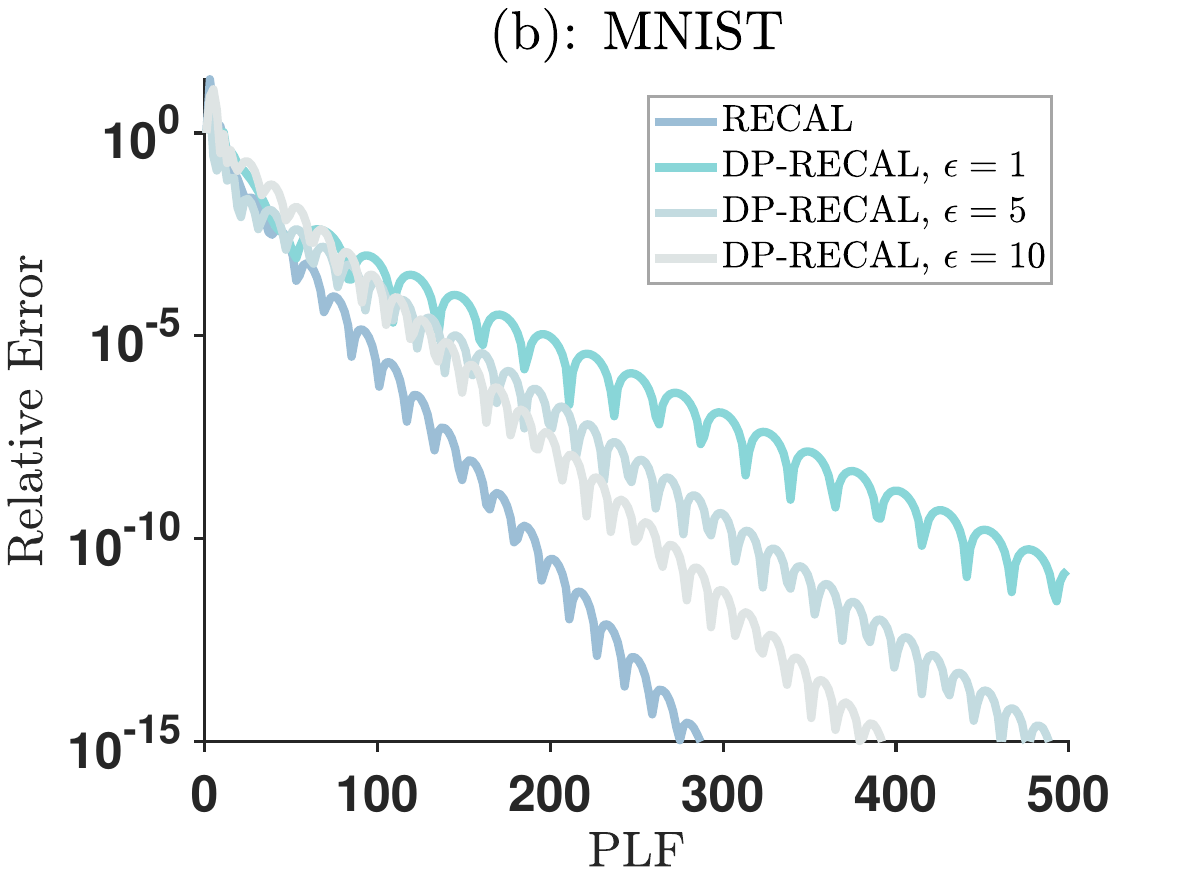}}\hspace{-2.85mm}
\subfigure{
\includegraphics[width=0.25\linewidth]{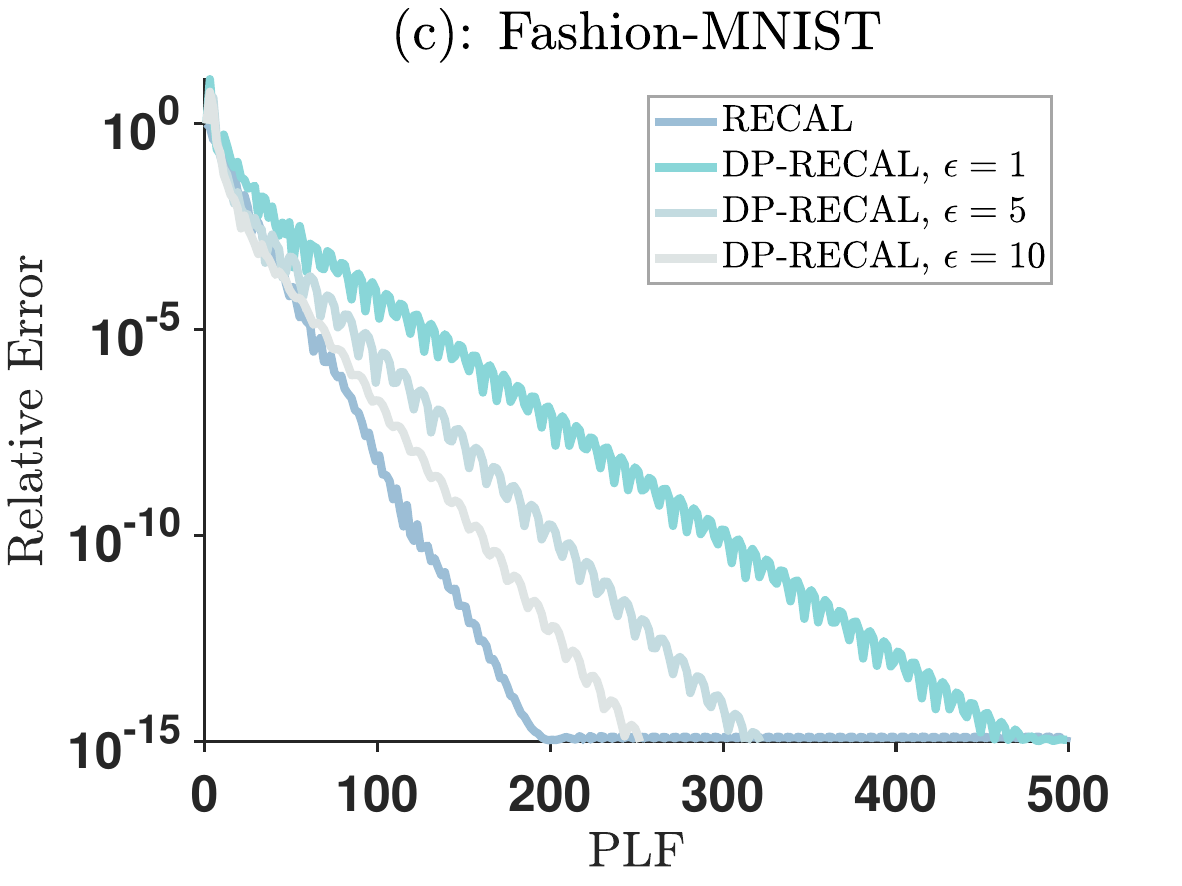}}\hspace{-2.85mm}
\subfigure{
\includegraphics[width=0.25\linewidth]{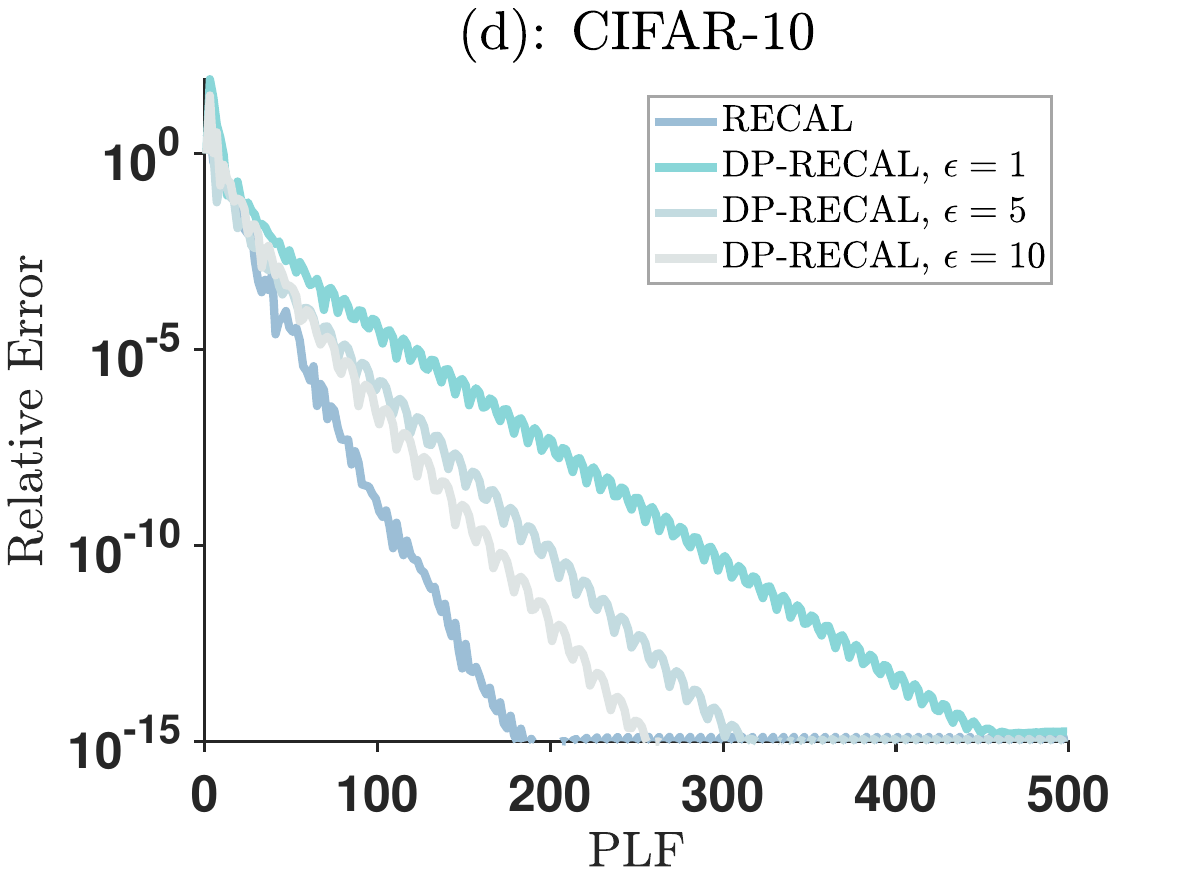}}
\caption{Convergence curves of DP-RECAL with various privacy budgets for the linear regression problem on four types of datasets.}
\label{fig:diffep}
\end{figure*}

\begin{figure*}[!t]
\centering
\setlength{\abovecaptionskip}{-7pt}
\subfigure{
\includegraphics[width=0.25\linewidth]{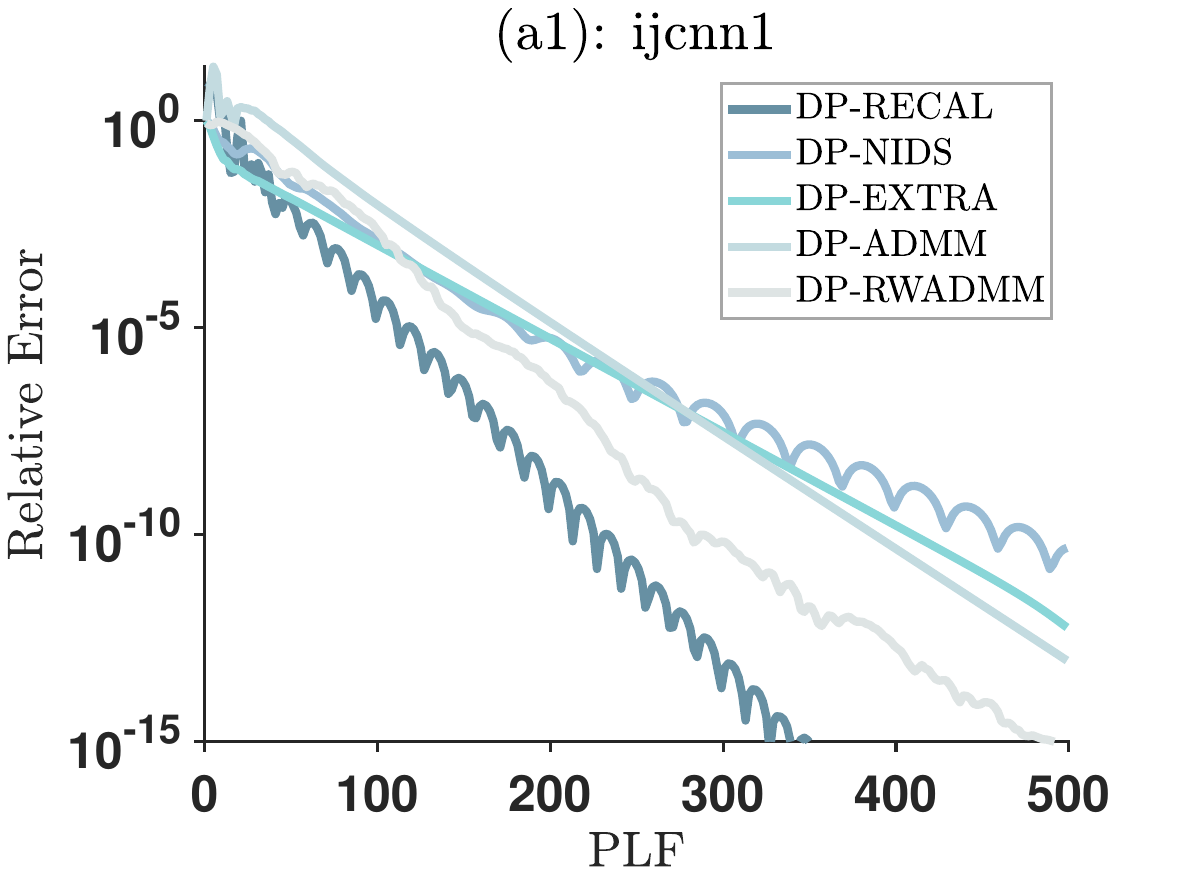}}\hspace{-2.85mm}
\subfigure{
\includegraphics[width=0.25\linewidth]{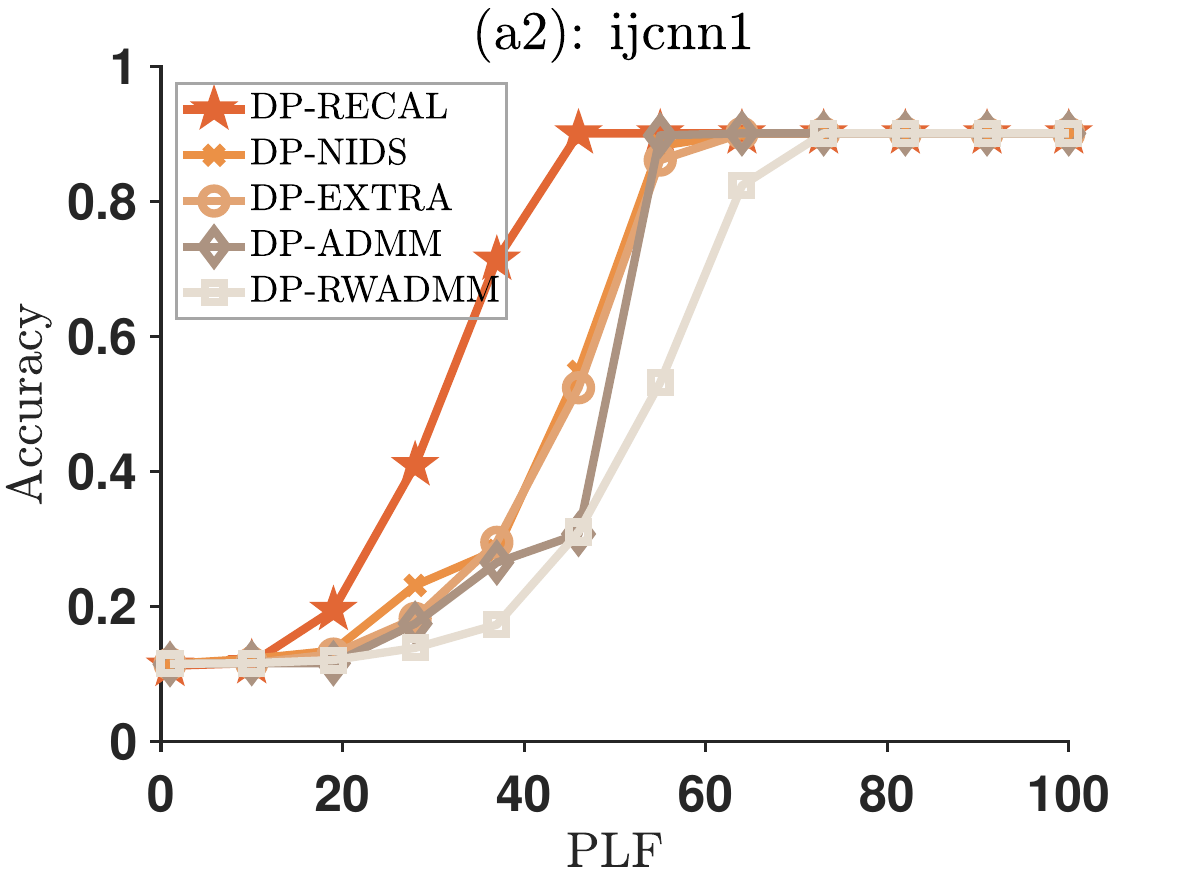}}\hspace{-2.85mm}
\subfigure{
\includegraphics[width=0.25\linewidth]{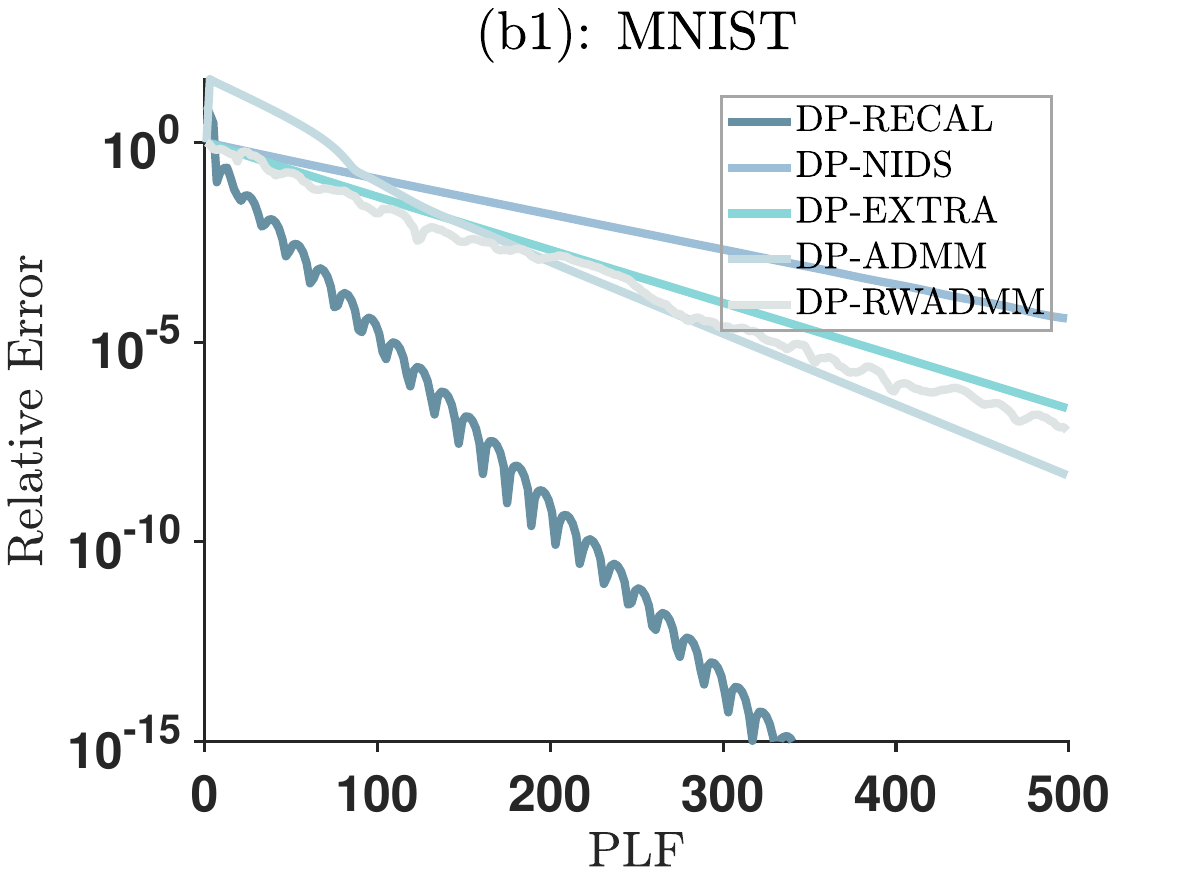}}\hspace{-2.85mm}
\subfigure{
\includegraphics[width=0.25\linewidth]{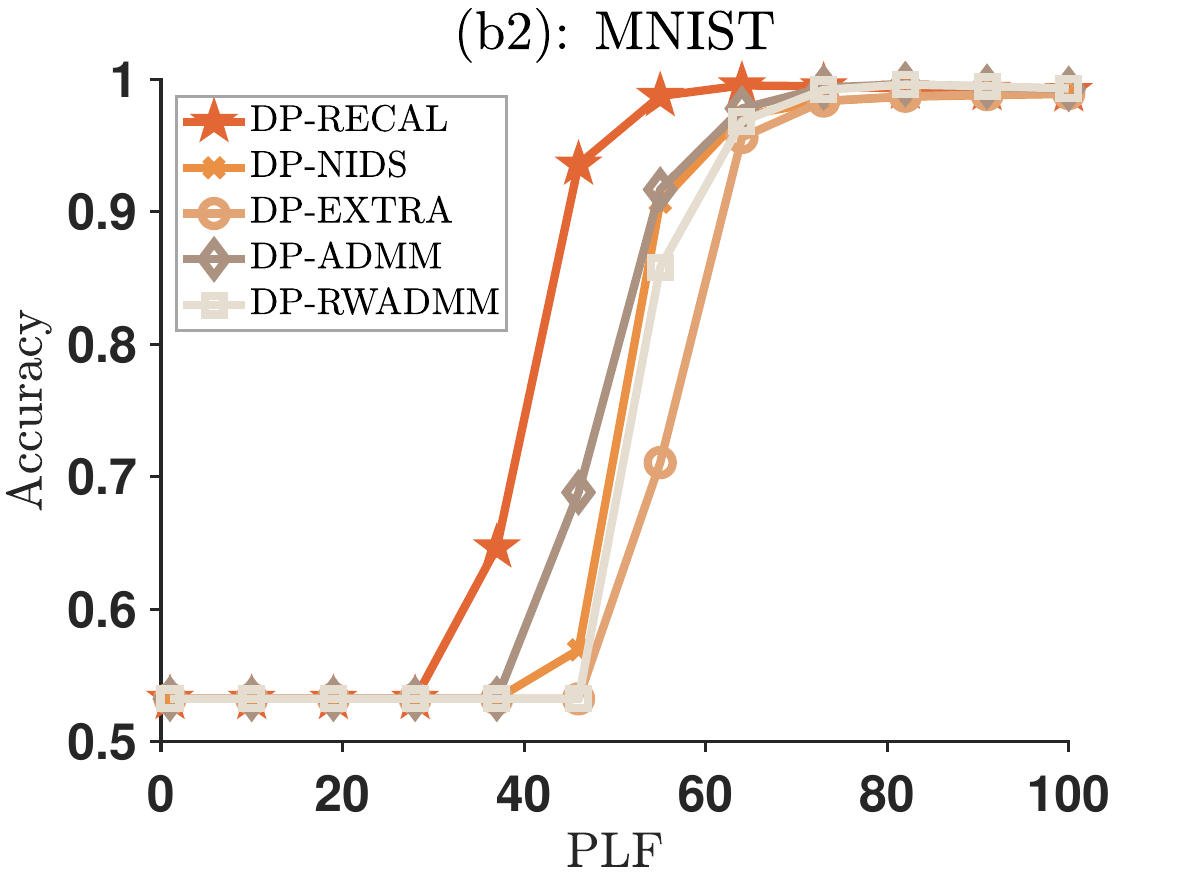}}
\subfigure{
\includegraphics[width=0.25\linewidth]{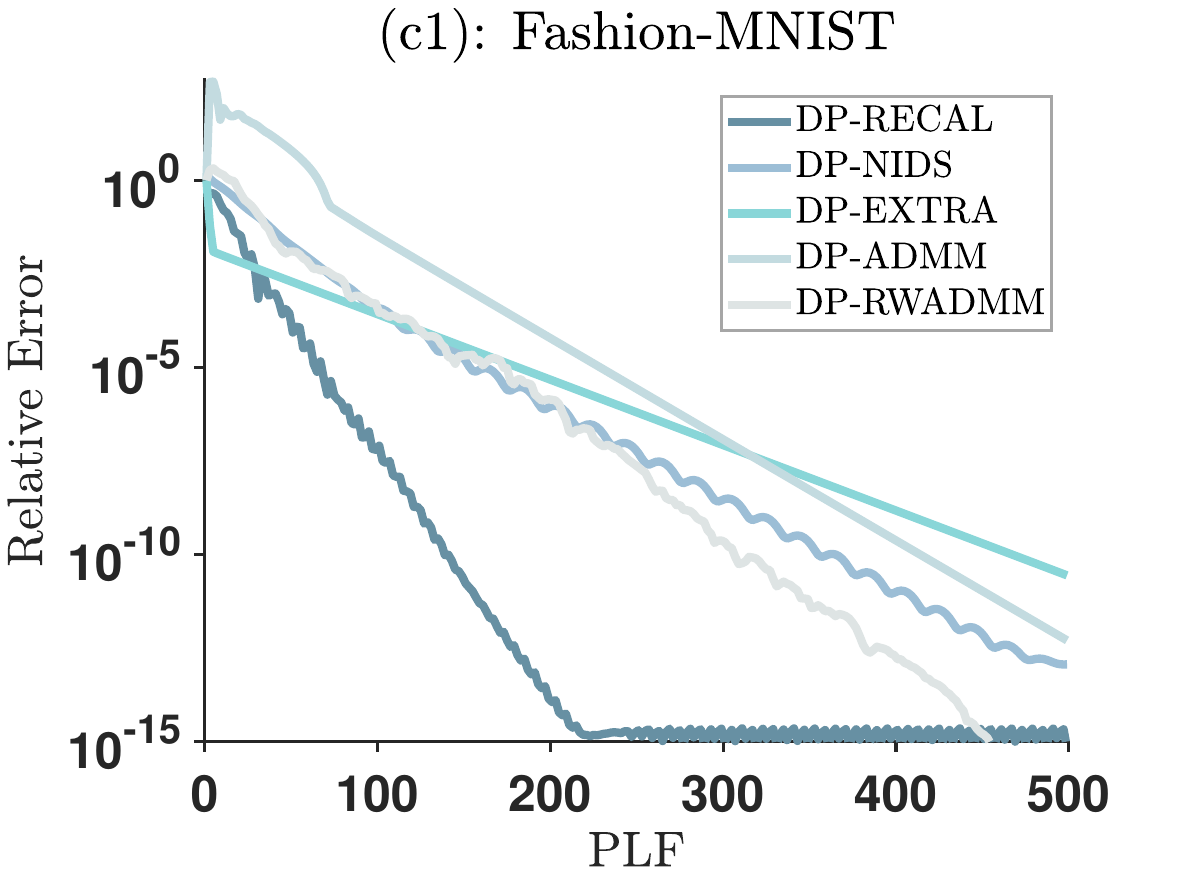}}\hspace{-2.85mm}
\subfigure{
\includegraphics[width=0.25\linewidth]{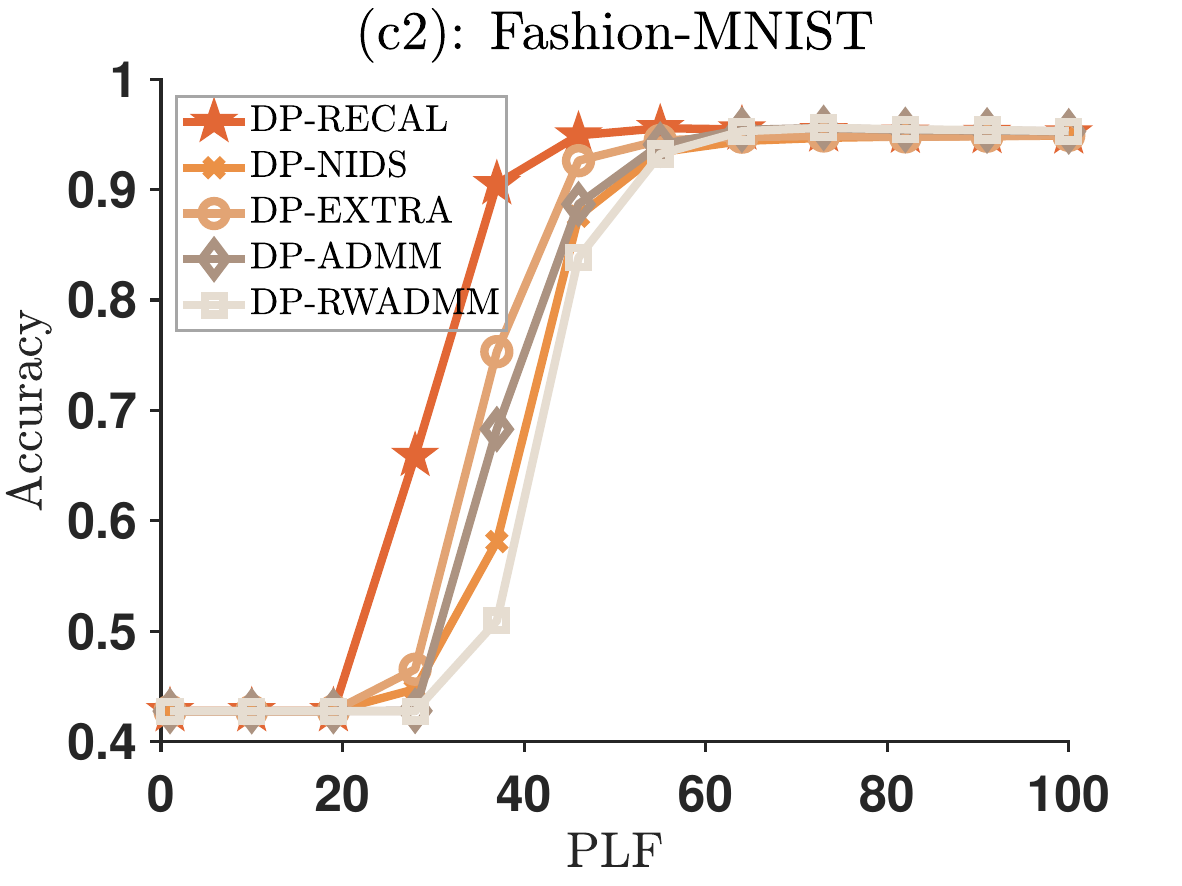}}\hspace{-2.85mm}
\subfigure{
\includegraphics[width=0.25\linewidth]{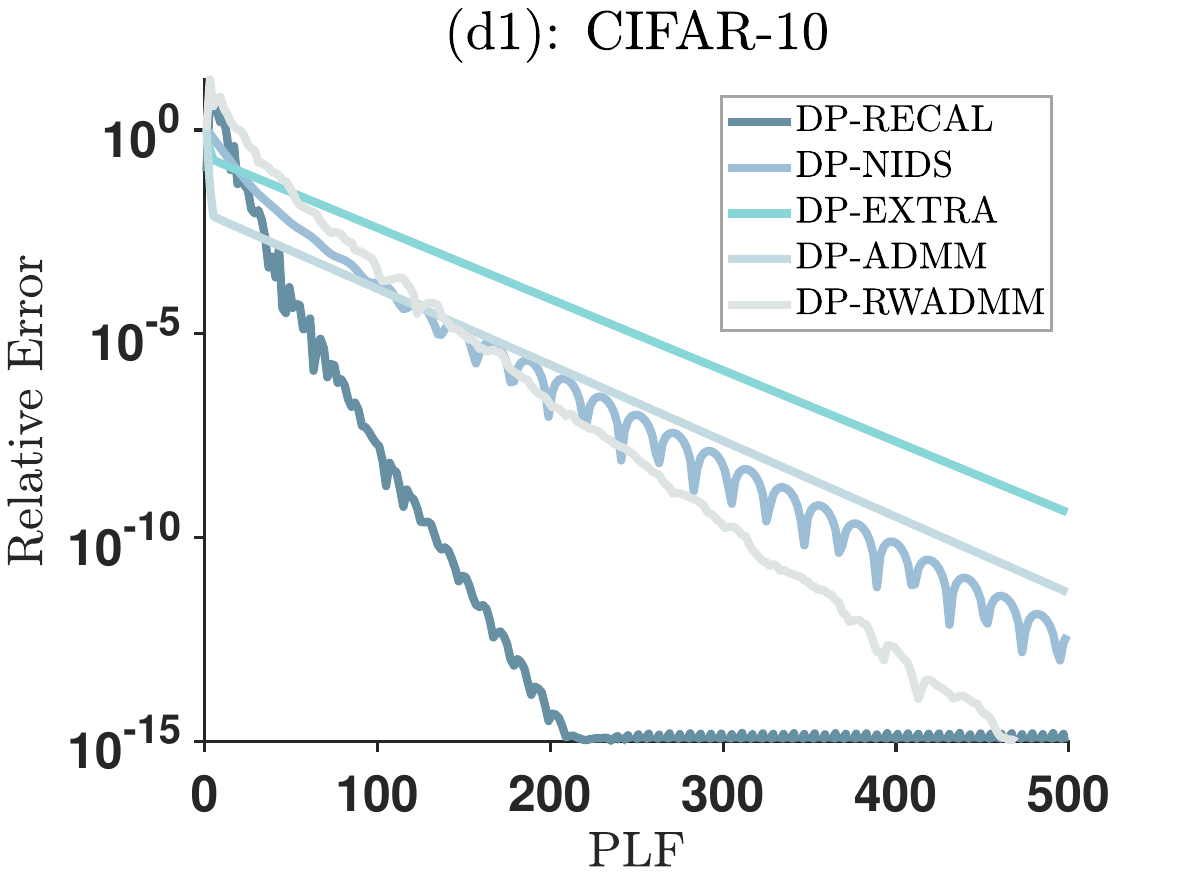}}\hspace{-2.85mm}
\subfigure{
\includegraphics[width=0.25\linewidth]{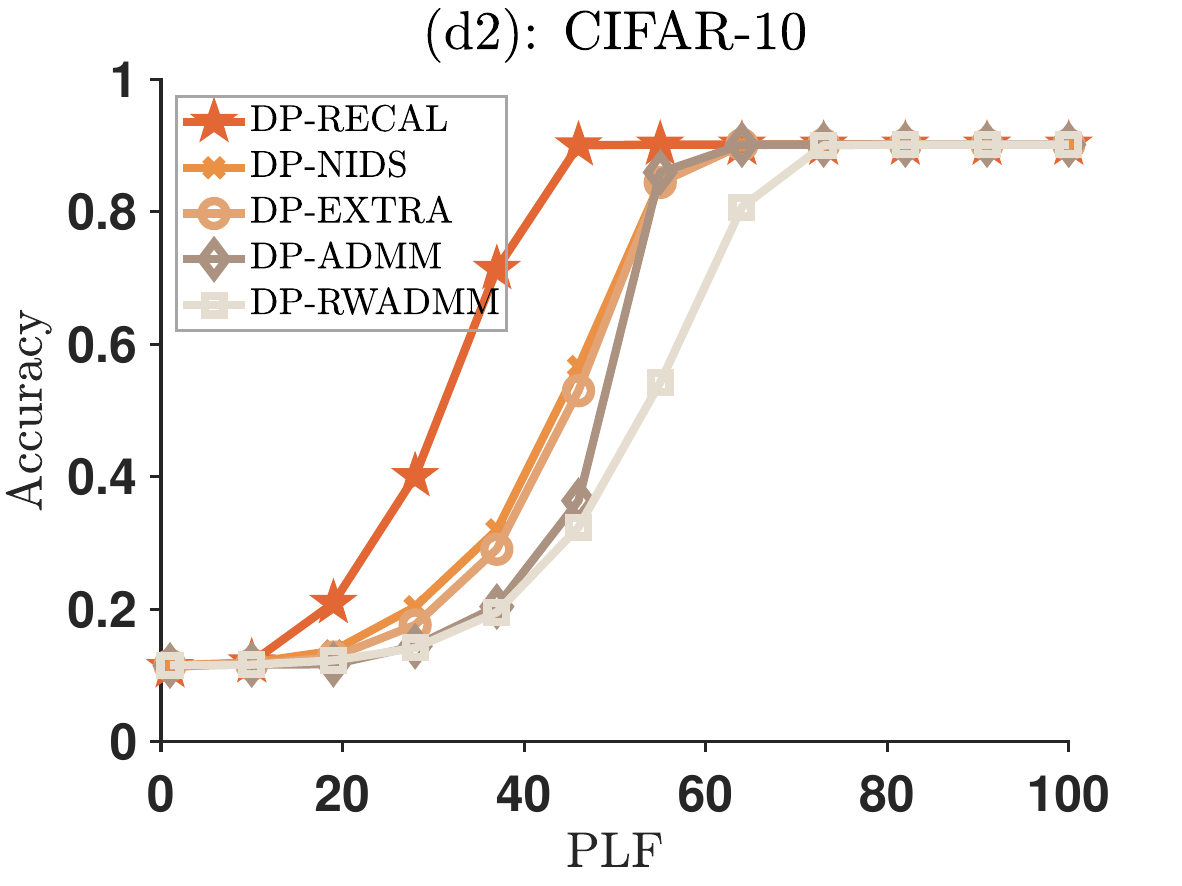}}
\caption{Convergence performance and classification accuracy of various differentially private algorithms with privacy budget $\epsilon = 12$ for the linear regression problems on four types of datasets.}
\label{fig:compare}
\end{figure*}

In this experiment, we firstly compare DP-RECAL with different privacy budgets and set RECAL, which is in non-private as a baseline, to evaluate the accuracy and the effectiveness of our method. Next, we compare various differentially private decentralized optimization algorithms with DP-RECAL to show its advantages. The latest relevant differentially private researches \cite{Zhang2017}, \cite{Huang2020} and \cite{Zhang2018} rooted in the ADMM algorithm should be fully considered. Note that \cite{Zhang2018} is pointed out by their authors to represent an improved version of \cite{Zhang2018_1}, thus we opt to incorporate \cite{Zhang2018} into our control group due to its better performance. NIDS \cite{NDIS}, EXTRA \cite{PGEXTR}, due to their wide-ranging applicability, are injected with DP noise in our experiments. Furthermore, we include a asynchronous gossip algorithm RWADMM \cite{Shah2018} in the comparison, which differs from the above algorithms.

For an internal comparison within DP-RECAL, we set different privacy budgets $\epsilon=\{1,5,10\}$ with $\delta=10^{-3}$ by using different initial values and consuming the same number of PLF. This allows us to understand the impact of different privacy budgets on the performance of DP-RECAL.

For a comprehensive comparison between algorithms, we maintain a consistent privacy budget of $\epsilon = 12$ across all algorithms. Subsequently, we evaluate their convergence performance and classification accuracy on datasets. Specifically, we assess DP-ADMM \cite{Zhang2017}, DP-NIDS \cite{NDIS}, DP-EXTRA \cite{PGEXTR}, DP-RWADMM \cite{Shah2018}, and our proposed approach DP-RECAL. These algorithms address regression problems featuring nonsmooth terms with fixed stepsizes, which allows us to construct Fig. \ref{fig:compare} for a visual representation of their distinctions while ensuring a uniform security level by maintaining a consistent PLF of 500. On the other hand, DP-ADMM \cite{Huang2020} is limited to decaying stepsizes, while DP-MRADMM \cite{Zhang2018} is not applicable to optimization problems involving nonsmooth terms. As these special cases are difficult to portray in images, Table \ref{Table-total} is additionally provided for a summarized overview of all relevant algorithms, including dimensions such as communication complexity, CPU time, convergence error, and classification accuracy. We set $\alpha_k = \alpha_0\cdot \gamma^k$ as well as PLF = 2000 for decaying stepsizes, $\alpha_k = \alpha_0$ as well as PLF = 300 for fixed stepsizes, and default the regular term $r=0$ in the regression problems for DP-MRADMM. Note that the convergence performance is expressed in Relative Errors ($\|\vx^k-\vx^*\|/\|\vx^0-\vx^*\|$), and classification accuracy is presented as Accuracy (number of correctly classified samples divided by total number of samples).

\begin{table*}[!t]
\renewcommand\arraystretch{1.5}
\begin{center}
\caption{Comparisons of various differentially private algorithms for linear regression problems on four types of datasets with the same privacy budget $\epsilon = 12$.}
\scalebox{1.05}{
\begin{tabular}{ccccccccccccc}
\hline
\multirow{2}{*}{$r(\vx)$}&\multirow{2}{*}{$\alpha$} &\multirow{2}{*}{PLF}&\multirow{2}{*}{ALG}&\multirow{2}{*}{Comm.}&\multicolumn{2}{c}{ijcnn1}&\multicolumn{2}{c}{MNIST}&\multicolumn{2}{c}{Fashion-MNIST}&\multicolumn{2}{c}{CIFAR-10}\\
\cline{6-13}
&&&&&{Time}&{RE.}&Time&RE.&Time&RE.&Time&RE.\\
\hline
\multirow{5}{*}{$\ell_1$}&$\mathcal{O}(\gamma^{-k})$&2000&DP-ADMM\cite{Huang2020}&6400&59.4&6.1e-8&158.3&1.9e-6&140.0&2.9e-8&1.0e+4&1.2e-9\\
\cline{2-13}
&\multirow{4}{*}{$\mathcal{O}(1)$}&\multirow{4}{*}{300}&DP-NIDS\cite{NDIS}&4800&12.0&1.5e-7&29.8&1.2e-3&26.4&2.9e-6&1835&5.8e-9\\
&&&DP-EXTRA\cite{PGEXTR}&9600&11.2&3.0e-8&25.5&9.6e-5&23.6&8.6e-7&1838&1.2e-6\\
&&&DP-ADMM\cite{Zhang2017}&9600&11.7&2.4e-8&25.2&1.6e-5&22.1&1.1e-7&1984&2.3e-8\\
&&&DP-RWADMM\cite{Shah2018}&9600&18.8&6.9e-9&41.8&2.4e-6&40.9&2.2e-10&3746&1.8e-10\\
\hline
\rowcolor{mygray}$\ell_1$&$\mathcal{O}(1)$&\textbf{300}&\textbf{DP-RECAL}&\textbf{2400}&\textbf{10.3}&\textbf{5.8e-15}&\textbf{26.2}&\textbf{4.7e-15}&\textbf{24.4}&\textbf{6.8e-15}&\textbf{1904}&\textbf{1.0e-15}\\
\hline
\hline
\multirow{2}{*}{-}&$\mathcal{O}(\gamma^{-k})$&2000&DP-MRADMM\cite{Zhang2018}&128000&1098&5.6e-10&3075&2.1e-7&2987&4.2e-10&3.6e+5&7.5e-11\\
\cline{2-13}
&$\mathcal{O}(1)$&300&DP-MRADMM\cite{Zhang2018}&19200&392&1.8e-10&1237&3.8e-9&1316&3.6e-11&1.7e+5&1.3e-11\\
\hline
\rowcolor{mygray}-&$\mathcal{O}(1)$&\textbf{300}&\textbf{DP-RECAL}&\textbf{2400}&\textbf{5.3}&\textbf{7.4e-15}&\textbf{23.4}&\textbf{6.0e-15}&\textbf{24.8}&\textbf{9.0e-16}&\textbf{1710}&\textbf{8.9e-16}\\
\hline
\end{tabular}}
\label{Table-total}
\end{center}
\end{table*}

\begin{figure*}[!t]
\setlength{\abovecaptionskip}{-7pt}
  \centering
  \includegraphics[width=1\linewidth]{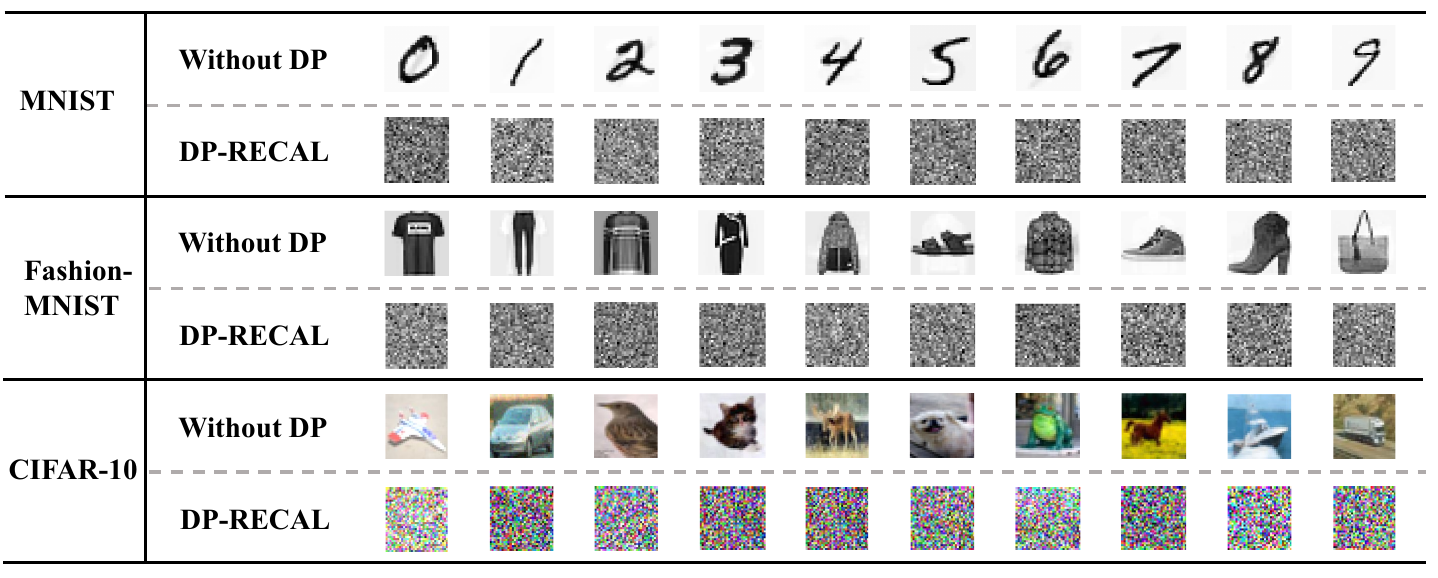}
  \caption{Comparisons of the results achieved by the adversary with DLG method \cite{Zhu2019} on classification of three datasets performed by DP-RECAL and RECAL.}
  \label{fig:dlg}
\end{figure*}

From Fig. \ref{fig:diffep}, DP-RECAL achieves faster convergence with a larger privacy budget $\epsilon$, which is consistent with the intuition, and when the privacy budget is large enough, the trained models of our approach is nearly as good as the one in non-private setting.  Fig. \ref{fig:compare} compares DP-RECAL with various differentially private algorithms and illustrates that our approach consistently presents the best convergence accuracy when consuming the same privacy budget for linear regression on different datasets. Moreover, Table \ref{Table-total}, where Comm. and RE. stand for communication complexity and relative error, lists various performance metrics of the algorithms based on both smooth and non-smooth objective functions, with DP-MRADMM being applicable to both fixed and decaying stepsizes. This table reveals that at the same security level, DP-RECAL holds not only better convergence accuracy but also less communication complexity compared to other algorithms. As far as CPU time is concerned, DP-RECAL doesn't seem to be inferior or even better, especially compared to those algorithms with decaying stepsizes and internal loops.
\subsection{Attack simulation}
To show that our method can indeed protect the privacy of participating agents, we simulate an adversary with the ability described in Section \ref{Attack model}. Firstly, this adversary extracts information (state variables) while classifying the MNIST, Fashion-MNIST, and CIFAR-10 datasets using DP-RECAL and non-private RECAL, respectively. Then, it deduces gradients and applies DLG method \cite{Zhu2019} aiming to recover datasets. Moreover, the latest attack methods \cite{Geiping2020} and \cite{Yin2021}, which have the similar core idea as that of DLG, are also used to test the security of our differentially private algorithm \footnote{Some experimental details are given in the supplement.}. In addition, ijcnn1 dataset is not shown here due to its inability to be presented as images, which can also be seen in the supplement for details.

Fig. \ref{fig:dlg} depicts the results achieved by the adversary's attacks on the classification process of three datasets performed by DP-RECAL and RECAL. It is evident that the adversary could effectively reconstruct the original image from a shared model update of the non-private RECAL. In contrast, DP-RECAL effectively thwarts the adversary's attempts to infer the original images, resulting in meaningless noise. Due to the same structure of DP-RECAL and RECAL, the experimental results also demonstrate that privacy-preserving performance is caused by differential privacy noise in our algorithm.

\section{Conclusion}\label{Conclusion}
We have considered a convex decentralized optimization problem on an undirected network where the loss function consists of local smooth terms and a  global nonsmooth term. For privacy concern of such problem, we proposed an efficient algorithm named DP-RECAL. Compared with preexisting privacy-preserving algorithms, DP-RECAL offers a distinctive advantage: it reduces the PLF with similar convergence accuracy by means of an AFBS-based framework and its relay communication technique. As substantiated through analysis, DP-RECAL demonstrates exceptional privacy preservation and commendable communication efficiency. Moreover, we have established DP-RECAL's convergence with uncoordinated network-independent stepsizes and its linear convergence under metric subregularity. However, it is pertinent to acknowledge that DP-RECAL's convergence rate may not be as competitive as popular decentralized optimization algorithms, due to the inherent sequential nature of agent computations. Thus, our subsequent attention could also be paid to find the optimal PLF to achieve optimal privacy protection for differentially private decentralized algorithms while maintaining parallel computation at the agents.


\vspace{-15 mm}
\begin{IEEEbiography}[{\includegraphics[width=1.3in,height=1.2in,clip,keepaspectratio]{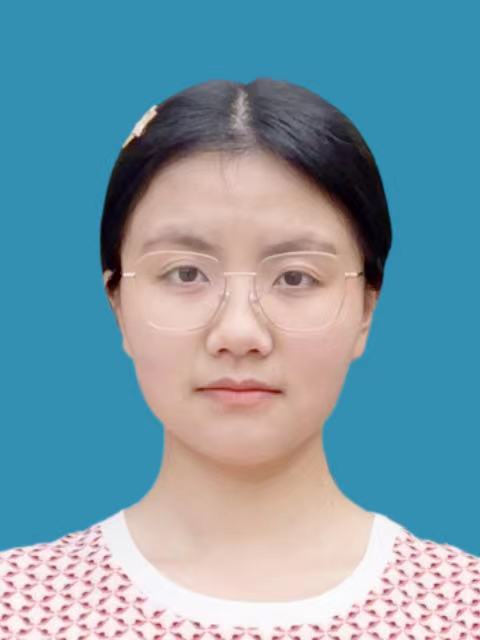}}]{Luqing Wang}
received the B.S. degree in Cyber Science and Engineering
from Southeast University, Nanjing, China, in 2023. She is currently pursuing the Ph.D.
degree in computer science and technology with the School of Computer Science and Engineering, Southeast University, Nanjing, China. Her current research focuses on distributed optimization, federal learning, and algorithmic security.
\end{IEEEbiography}

\vspace{-15 mm}
\begin{IEEEbiography}[{\includegraphics[width=1.3in,height=1.2in,clip,keepaspectratio]{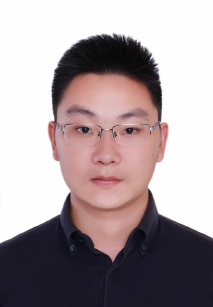}}]{Luyao Guo}
received the B.S. degree in information and computing science
from Shanxi University, Taiyuan, China, in 2020. He is currently pursuing the Ph.D.
degree in applied mathematics with the School of Mathematics, Southeast University, Nanjing, China. His current research focuses on distributed optimization and learning, federal learning, and decentralized learning.
\end{IEEEbiography}
\vspace{-15 mm}
\begin{IEEEbiography}[{\includegraphics[width=1.3in,height=1.2in,clip,keepaspectratio]{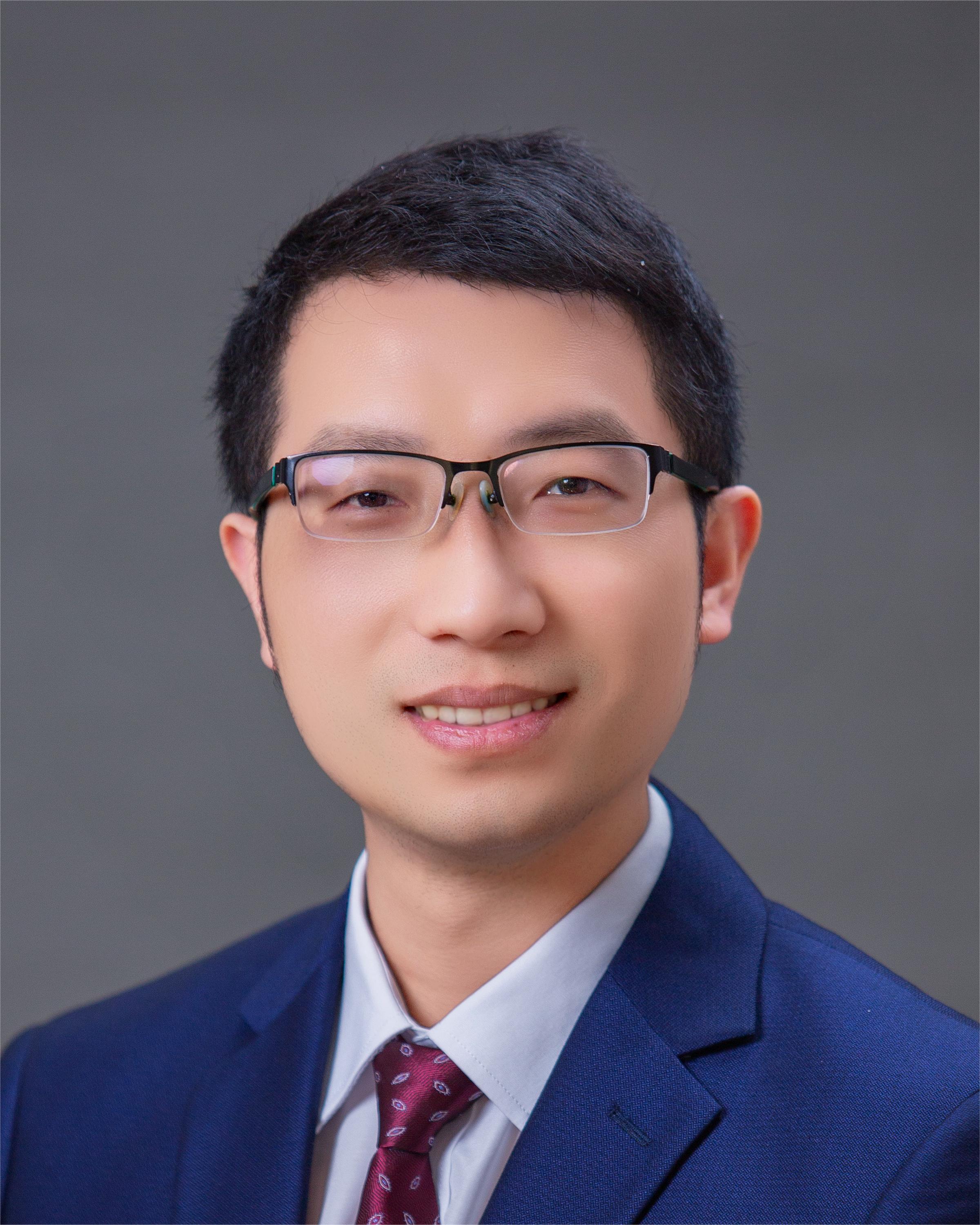}}]{Shaofu Yang (Member, IEEE)}
received the B.S. and M.S. degrees in applied mathematics from the Department of Mathematics, Southeast University, Nanjing, China, in 2010 and 2013, respectively, and the Ph.D. degree in engineering from the Department of Mechanical and Automation Engineering, The Chinese University of Hong Kong, Hong Kong, in 2016. He was a Post-Doctoral Fellow with the City University of Hong Kong, Hong Kong, in 2016. He is currently an Associate Professor with the School of Computer Science and Engineering, Southeast University. His current research interests include distributed optimization and learning, game theory, and their applications.
\end{IEEEbiography}
\vspace{-15 mm}
\begin{IEEEbiography}[{\includegraphics[width=1.3in,height=1.2in,clip,keepaspectratio]{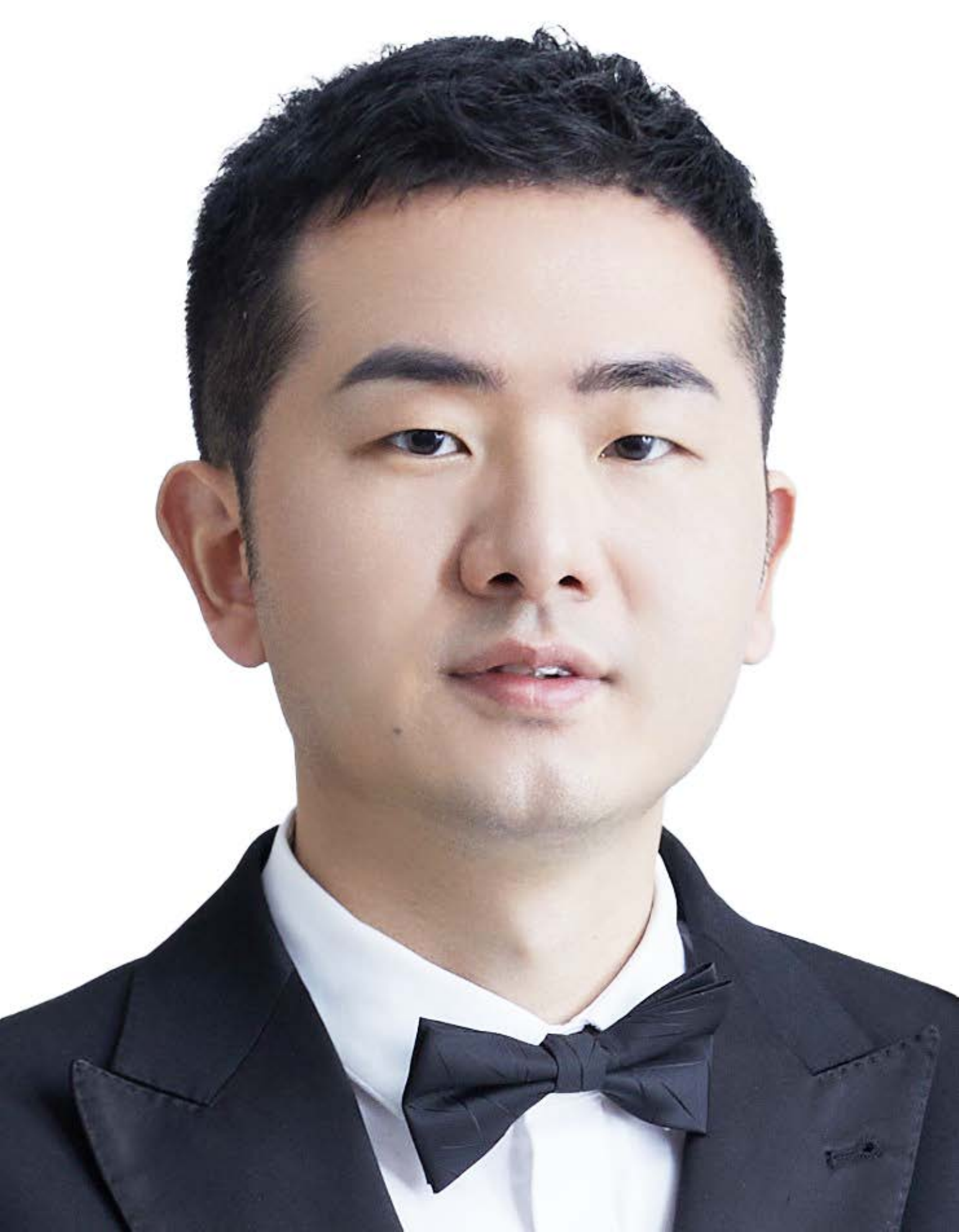}}]{Xinli Shi (Senior Member, IEEE)}
received the B.S. degree in software engineering, the M.S. degree in applied mathematics and the Ph.D. degree in control science and engineering from Southeast University, Nanjing, China, in 2013, 2016 and 2019, respectively. He was the recipient of the Outstanding Ph.D. Degree Thesis Award from Jiangsu Province, China. He is a recipient of the Australian Research Council Discovery Early Career Researcher Award. He is currently an associate professor at Southeast University.
His current research interests include distributed optimization, reinforcement learning, and network control systems.
\end{IEEEbiography}

\newpage

\setcounter{page}{1}
\section{SUPPLEMENTAL MATERIALS}
\appendices
\section{Gradient inference of DP-RECAL}\label{Gradient inference of DP-RECAL}

In RECAL \eqref{OurAlg2}, given the shared nature of $\{\vu^k\}$ and $\{\vx^k\}$ during iterations, an adversary could exploit them to infer agents' private cost functions through iteration sequences. In such a context, the adversary's inference process would resemble a mathematical induction. Assuming that variables $\vy_{i_k}^k$ and $\vlambda_{i_k}^k$, generated before the $k$-iteration, are known to the adversary, we can demonstrate how the adversary deduces $\vy_{i_k}^{k+1}$ and $\nabla f_{i_k}(\vy_{i_k}^k)$ based on the information disclosed by the activated agent $i_k$.

The adversary can deduce $\vlambda_{i_k}^{k+\frac{1}{2}}$ by utilizing \eqref{OurAlg2a}, and intercept $\vx^{k+1}$ as well as $\vu^{k+1}$ from the communication channel. Then, $\vlambda_{i_k}^{k+1}$, $\vy_{i_k}^{k+1}$ and $\nabla f_{i_k}(\vy_{i_k}^k)$ can be reduced from \eqref{OurAlg2e}, \eqref{OurAlg2d} and \eqref{OurAlg2c} in turn. Specifically, the derivation is as follows.
\begin{align*}
\vlambda_{i_k}^{k+\frac{1}{2}}&=\vlambda_{i_k}^{k}+\beta(\vx^{k}-\vy_{i_k}^{k}),\\
\vlambda_{i_k}^{k+1}&=\vu^{k+1}-\vu^k+\vlambda_{i_k}^k,\\
\vy_{i_k}^{k+1}&=\vy_{i_k}^k+\vx^{k+1}-\vx^k+\frac{1}{\beta}(\vlambda_{i_k}^{k+\frac{1}{2}}-\vlambda^{k+1}_{i_k}),\\
\nabla f_{i_k}(\vy^k_{i_k}) &= \frac{1}{\alpha_i}(\vy_{i_k}^k-\vy^{k+1}_{i_k})+\vlambda_{i_k}^{k+\frac{1}{2}}.
\end{align*}
So on for subsequent steps, the adversary recursively obtains this information, eventually inferring the private $\nabla f$.
\section{Proof of Lemma \ref{lem:active}}\label{APP0}
\begin{IEEEproof}
From the law of total probability, we have $g_i=\sum_{j\neq i}P_{ij}g_j$. Combining with $\sum_{i=1}^n g_i=1$, it holds that
\begin{align}\label{Pg}
\left(
\begin{array}{c}
 \mI-\mP\tr \\
 \one\tr \\
 \end{array}
 \right)\vg =\left(
\begin{array}{c}
 \bf{0} \\
 1 \\
 \end{array}
 \right).
\end{align}
Recalling that we have $P_{ij}=\frac{1}{d_i}$ if $(i,j)\in \mathcal{E} $; otherwise, $P_{ij}=0$.
Applying the Laplace matrix $\mL$ and the degree matrix $\mD$ of the graph $\mathcal{G}$, we can rewrite \eqref{Pg} as
\begin{align}\label{Geq1}
\left(
\begin{array}{c}
 \mL\mD^{-1} \\
 \one\tr \\
 \end{array}
 \right)\vg =\left(
\begin{array}{c}
 \bf{0} \\
 1 \\
 \end{array}
 \right).
\end{align}
To solve \eqref{Geq1}, we can consider its least squares form
\begin{align*}
 (\mD^{-1}\mL^2\mD^{-1}+\one\one\tr)\vg=\one.
\end{align*}
It is easy to verify that the matrix $\mL^2+\mD\one\one\tr\mD$ is positive definite, which implies that the matrix $\mD^{-1}\mL^2\mD^{-1}+\one\one\tr$ being positive definite. Therefore, the algebraic equation \eqref{Geq1} is equivalent to its least squares form, and the solution of \eqref{Geq1} can be expressed as
\begin{align}\label{valueg}
\vg=(\mD^{-1}\mL^2\mD^{-1}+\one\one\tr)^{-1}\cdot \one.
\end{align}
Furthermore, it follows from $\mL(\mD^{-1}\vg)=0$ that $\vg=c \mD\one$, where $c$ is a constant. Since $\one\tr\vg=1$, it holds that $c\cdot\one\tr\mD\one=1$. By the positive definiteness of $\mD$, we have $c>0$, which implies that $g_i>0$, $i\in \mathcal{V}$.

Thus, we complete the proof.
\end{IEEEproof}
\section{Proof of Lemma \ref{prop1}}\label{APP1}
\begin{IEEEproof}
It follows from \eqref{Alg3a} that
\begin{align}\label{dual1}
\big\langle \vlambda^*-\hat{\vlambda}^{k+\frac{1}{2}},-\mU\vx^k+\vy^k+\frac{1}{\beta}(\hat{\vlambda}^{k+\frac{1}{2}}-\vlambda^k)\big\rangle \geq0.
\end{align}
Rearranging \eqref{dual1} with $\mU\vx^*-\vy^*=0$, we obtain
\begin{align}\label{dual1-1}
&\big\langle \vlambda^*-\hat{\vlambda}^{k+\frac{1}{2}},-\mU\vx^*+\vy^*-\mU\vx^k+\vy^k\nonumber\\
&\quad +\frac{1}{\beta}(\hat{\vlambda}^{k+\frac{1}{2}}-\vlambda^k)\big\rangle \geq0.
\end{align}
It follows from \eqref{Alg3b} that
$$
\hat{\vx}^{k+1}\in (\mI+\cdot n \partial r)^{-1}\cdot \left(\vx^{k}-(\mU\tr\hat{\vlambda}^{k+\frac{1}{2}}-\vvarepsilon^k)\right),
$$
which implies
$
n\partial r(\hat{\vx}^{k+1}) \in -\mU\tr\hat{\vlambda}^{k+\frac{1}{2}}-\hat{\vx}^{k+1}+\vx^k+\vvarepsilon^{k}.
$
From the monotonicity of $\partial r$, it holds that
\begin{align}\label{primal1}
&\big\langle \vx^*-\hat{\vx}^{k+1},-\mU\tr\vlambda^*+\mU\tr\hat{\vlambda}^{k+\frac{1}{2}}\nonumber\\
&\quad+\hat{\vx}^{k+1}-\vx^k-\vvarepsilon^{k}\big\rangle\geq 0,
\end{align}
which is equivalent to
\begin{align}\label{primal1-1}
&\big\langle \vx^*-\hat{\vx}^{k+1},-\mU\tr\vlambda^*+\mU\tr\vlambda^k+\mU\tr\hat{\vlambda}^{k+\frac{1}{2}}-\mU\tr\vlambda^k\nonumber\\
&\quad+\hat{\vx}^{k+1}-\vx^k-\vvarepsilon^{k}\big\rangle\geq 0.
\end{align}
By \eqref{Alg3c}, we have
\begin{align*}
&\big\langle \vy^*-\hat{\vy}^{k+1}, \nabla f(\vy^k)-\hat{\vlambda}^{k+\frac{1}{2}}+\mGamma^{-1}(\hat{\vy}^{k+1}-\vy^k) \big\rangle \geq 0.
\end{align*}
Similarly, rearranging this inequality with $\vlambda^*-\nabla f(\vy^*)=0$, we obtain
\begin{align}\label{primal2-1}
&\big\langle \vy^*-\hat{\vy}^{k+1}, \vlambda^*- \nabla f(\vy^*) -\vlambda^k+ \nabla f(\vy^k) \nonumber\\
&\quad -\hat{\vlambda}^{k+\frac{1}{2}}+\vlambda^k+\mGamma^{-1}(\hat{\vy}^{k+1}-\vy^k) \big\rangle \geq 0.
\end{align}
Remembering the structure of $\mA$, $\mH$ and $\FP$, we combine \eqref{dual1-1} \eqref{primal1-1} and \eqref{primal2-1} and then have
\begin{align}\label{THEQ1}
&\big\langle \vw^*-\hat{\vw}^{k+\frac{1}{2}}, -\mA\vw^*-\FP(\vw^*)+\mA\vw^k+\FP(\vw^k)\nonumber\\
&\quad+\mH\hat{\vw}^{k+\frac{1}{2}}-\mH\vw^k-\mD^k\big\rangle\geq0,
\end{align}
where $\vw^{k+\frac{1}{2}}=\col\{\vlambda^{k+\frac{1}{2}},\vx^{k+1},\vy^{k+1}\}$, $\mD^{k}=[(\vvarepsilon^k)\tr,0,0]\tr$.
Since $f_i(\vy_i)$ is $L_i$-smooth and convex, it holds that
\begin{align}
&\big\langle \vy_i^*-\hat{\vy_i}^{k+\frac{1}{2}},\nabla f_i(\vy_i^k)-\nabla f_i(\vy_i^*)\big\rangle\nonumber\\
&=\big\langle \vy_i^k-\hat{\vy_i}^{k+\frac{1}{2}}+\vy_i^*-\vy_i^k, \nabla f_i(\vy_i^k)-\nabla f_i(\vy_i^*)\big\rangle \nonumber\\
&\leq \big\langle \vy_i^k-\hat{\vy_i}^{k+\frac{1}{2}}, \nabla f_i(\vy_i^k)-\nabla f_i(\vy_i^*) \big\rangle \nonumber\\
&\quad-\frac{1}{L_i} \big\|\nabla f(\vy_i^k)-\nabla f(\vy_i^*)\big\|^2\nonumber\\
&=-\big\|\frac{\sqrt{L_i}}{2}(\hat{\vy_i}^{k+\frac{1}{2}}-\vy_i^k)+\frac{1}{\sqrt{L_i}}(\nabla f_i(\vy_i^k)-\nabla f(\vy_i^*))\big\|^2\nonumber\\
&\quad+ \frac{L_i}{4} \big\|\hat{\vy_i}^{k+\frac{1}{2}}-\vy_i^k\big\|^2 \leq \frac{L_i}{4} \big\|\hat{\vy_i}^{k+\frac{1}{2}}-\vy_i^k\big\|^2,
\end{align}
which implies that
\begin{align}\label{LY}
&\big\langle \vw^*-\hat{\vw}^{k+\frac{1}{2}}, \FP(\vw^k)- \FP(\vw^*)\big\rangle\nonumber\\
&=\sum_{i=1}^n \big\langle \vy_i^*-\hat{\vy_i}^{k+\frac{1}{2}},\nabla f_i(\vy_i^k)-\nabla f_i(\vy_i^*)\big\rangle\nonumber\\
&\leq \sum_{i=1}^n \frac{L_i}{4} \big\|\hat{\vy_i}^{k+\frac{1}{2}}-\vy_i^k\big\|^2= \frac{1}{4} \big\|\hat{\vy}^{k+\frac{1}{2}}-\vy^k\big\|^2_{\mQ},
\end{align}
where $\mQ=\diag\{L_1\mI,\cdots,L_n\mI\}$.

Further we define two operators
\begin{align*}
{\mJ}=\left(
               \begin{array}{ccc}
                 \frac{1}{\beta}\mI & \frac{\mU}{2} & -\frac{\mI}{2} \\
                 \frac{\mU\tr }{2}& \mI & 0 \\
                 -\frac{\mI}{2}& 0 & \mGamma^{-1} \\
               \end{array}
             \right), {\mK}=\left(
               \begin{array}{ccc}
                 0 & -\frac{\mU}{2} & \frac{\mI}{2} \\
                 \frac{\mU\tr}{2} & 0 & 0 \\
                 -\frac{\mI}{2}& 0 & 0 \\
               \end{array}
             \right),
\end{align*}
where $\mJ+\mK=\mH$.
Note that the operator $\mK$ and $\mA$ is affine with a skew-symmetric matrix. We have
\begin{align*}
&\big\langle \vw^{*}-\hat{\vw}^{k+\frac{1}{2}},\mK\vw^{*}-\mK\hat{\vw}^{k+\frac{1}{2}}\big\rangle\equiv 0.\\
&\big\langle \vw^{*}-\vw^{k},(\mA-\mK)\vw^{*}-(\mA-\mK)\vw^{k}\big\rangle\equiv 0.
\end{align*}
Combining \eqref{THEQ1} and \eqref{LY}, we have
\begin{align}\label{THEQ3}
&\big\langle \vw^*-\hat{\vw}^{k+\frac{1}{2}}, -\mA\vw^*-\FP(\vw^*)+\mA\vw^k+\FP(\vw^k)\nonumber\\
&\quad+\mH\hat{\vw}^{k+\frac{1}{2}}-\mH\vw^k-\mD^k\big\rangle\nonumber\\
&\leq \big\langle \vw^*-\hat{\vw}^{k+\frac{1}{2}},(\mA-\mK)(\vw^k-\vw^*)+\mJ(\hat{\vw}^{k+\frac{1}{2}}-\vw^k) \big\rangle \nonumber\\
&\quad + \frac{1}{4}\big\|\hat{\vy}^{k+\frac{1}{2}}-\vy^k\big\|^2_{\mQ}-\big\langle  \vw^*-\hat{\vw}^{k+\frac{1}{2}}, \mD^k\big\rangle\nonumber\\
&=\big\langle \vw^*-\vw^{k},(\mA-\mK)(\vw^k-\vw^*)+\mJ(\hat{\vw}^{k+\frac{1}{2}}-\vw^k) \big\rangle \nonumber\\
&\quad +\big\langle \vw^k-\hat{\vw}^{k+\frac{1}{2}},(\mA-\mK)(\vw^k-\vw^*)+\mJ(\hat{\vw}^{k+\frac{1}{2}}-\vw^k) \big\rangle \nonumber\\
&\quad + \frac{1}{4} \big\|\hat{\vy}^{k+\frac{1}{2}}-\vy^k\big\|^2_{\mQ}-\big\langle  \vw^*-\hat{\vw}^{k+\frac{1}{2}}, \mD^k\big\rangle\nonumber\\
&=\big\langle \vw^*-\vw^{k},\mJ(\hat{\vw}^{k+\frac{1}{2}}-\vw^k) \big\rangle -\|\hat{\vw}^{k+\frac{1}{2}}-\vw^k\|^2_{\mJ}\nonumber\\
&\quad+\big\langle \vw^k-\hat{\vw}^{k+\frac{1}{2}},(\mA-\mK)(\vw^k-\vw^*)\big\rangle\nonumber\\
&\quad + \frac{1}{4} \big\|\hat{\vy}^{k+\frac{1}{2}}-\vy^k\big\|^2_{\mQ}-\big\langle  \vw^*-\hat{\vw}^{k+\frac{1}{2}}, \mD^k\big\rangle\nonumber\\
&=\big\langle \vw^{k}-\vw^*,(\mH+\mA\tr)(\vw^k-\hat{\vw}^{k+\frac{1}{2}})\big \rangle -\big\|\hat{\vw}^{k+\frac{1}{2}}-\vw^k\big\|^2_{\mJ}\nonumber\\
&\quad + \frac{1}{4} \big\|\hat{\vy}^{k+\frac{1}{2}}-\vy^k\big\|^2_{\mQ}+\big\langle  \hat{\vw}^{k+\frac{1}{2}}-\vw^*, \mD^k\big\rangle.
\end{align}
Due to $\mS^{-1}(\mH+\mA\tr)(\vw^k-\hat{\vw}^{k+\frac{1}{2}})=\vw^k-\hat{\vw}^{k+1}$, we have
\begin{align*}
&\big\langle \vw^{k}-\vw^*,(\mH+\mA\tr)(\vw^k-\hat{\vw}^{k+\frac{1}{2}})\big\rangle\\
&=\big\langle \vw^{k}-\vw^*,\vw^k-\hat{\vw}^{k+1}\big\rangle_{\mS}.
\end{align*}
In addition, we conclude that
\begin{align*}
\frac{1}{4} \big\|\hat{\vy}^{k+\frac{1}{2}}-\vy^k\big\|^2_{\mQ}-\big\|\hat{\vw}^{k+\frac{1}{2}}-\vw^k\big\|^2_{\mJ}=-\big\|\hat{\vw}^{k+1}-\vw^k\big\|^2_{\widehat{\mJ}},
\end{align*}
where
$$
\widehat{\mJ}=\left(
               \begin{array}{ccc}
                 \frac{1}{\beta}\mI & -\frac{\mU}{2}& \frac{\mI}{2} \\
                 -\frac{\mU\tr }{2}& \mI& 0 \\
                 \frac{\mI}{2}& 0 & \mGamma^{-1}-\frac{\mQ}{4} \\
               \end{array}
             \right).
$$
Then, from \eqref{THEQ3}, we obtain
\begin{align}
&\big\langle \vw^{k}-\vw^*,\vw^k-\hat{\vw}^{k+1}\big\rangle_{\mS}\nonumber\\
&\geq \left\|\hat{\vw}^{k+1}-\vw^k\right\|^2_{\widehat{\mJ}}-\big\langle \hat{\vw}^{k+\frac{1}{2}}-\vw^*, \mD^k\big\rangle.
\end{align}
Therefore, we have
\begin{align*}
&\big\|\hat{\vw}^{k+1}-\vw^*\big\|^2_{\mS}=\big\|\hat{\vw}^{k+1}-\vw^k+\vw^k-\vw^*\big\|^2_{\mS}\\
&=\big\|\vw^k-\vw^*\big\|^2_{\mS}+\big\|\hat{\vw}^{k+1}-\vw^k\big\|^2_{\mS}\\
&\quad+2\big\langle \vw^k-\vw^*,\hat{\vw}^{k+1}-\vw^k\big\rangle_{\mS}\\
&\leq\big\|\vw^k-\vw^*\big\|^2_{\mS}-\big\|\hat{\vw}^{k+1}-\vw^k\big\|^2_{\widehat{\mS}}+\big\langle  \hat{\vw}^{k+\frac{1}{2}}-\vw^*, \mD^k\big\rangle,
\end{align*}
where $\widehat{\mS}=2\widehat{\mJ}-\mS$. The stepsize condition $\alpha_i\in(0,\frac{2}{L_i+1})$ can ensure the positive definiteness of $\widehat{\mS}$.
\end{IEEEproof}

\section{Proof of Theorem \ref{thm1}}\label{APP2}
\begin{IEEEproof}
It follows from
$\vw^{k+1}=\vw^{k}+\mE_{k+1}(\hat{\vw}^{k+1}-\vw^{k}) $ that, for any $\vw^*\in\mathcal{W}^*$,
\begin{align}\label{ConPF1}
&\mathbb{E} {\left[\left\|\vw^{k+1}-\vw^{*}\right\|_{\overline{\mS}}^{2}\right] }\nonumber \\
&= \mathbb{E}\left[\begin{array}{l}
\left\|\vw^{k}-\vw^{*}\right\|_{\overline{\mS}}^{2}+\left\|\mE_{k+1}\left(\hat{\vw}^{k+1}-\vw^{k}\right)\right\|_{\overline{\mS}}^{2}\nonumber \\
+2\left(\vw^{k}-\vw^{*}\right)\tr \overline{\mS}\mE_{k+1}\left(\hat{\vw}^{k+1}-\vw^{k}\right)
\end{array}\right] \nonumber \\
&=\left\|\vw^{k}-\vw^{*}\right\|_{\overline{\mS}}^{2}+\mathbb{E}\left[\left\|\hat{\vw}^{k+1}-\vw^{k}\right\|_{\mE_{k+1} \mPi^{-1}\mS \mE_{k+1}}^{2}\right] \nonumber\\
&\quad+\mathbb{E}\left[2\left(\vw^{k}-\vw^{*}\right)^{\top} \mPi^{-1} \mS \mE_{k+1}\left(\hat{\vw}^{k+1}-\vw^{k}\right)\right] \nonumber\\
&=\left\|\vw^{k}-\vw^{*}\right\|_{\overline{\mS}}^{2}+\left\|\hat{\vw}^{k+1}-\vw^{k}\right\|_{\mS}^{2} \nonumber\\
&\quad+2\left(\vw^{k}-\vw^{*}\right)^{\top} \mS\left(\hat{\vw}^{k+1}-\vw^{k}\right),
\end{align}
where we use the results in Lemma \ref{lem2} at the last step. Note that for the last two terms in the above inequality, there exists an upper bound from \eqref{Con5}:
\begin{align}\label{ConPF2}
&\left\|\hat{\vw}^{k+1}-\vw^{k}\right\|_{\mS}^{2}+2\left(\vw^{k}-\vw^{*}\right)\tr\mS\left(\hat{\vw}^{k+1}-\vw^{k}\right) \nonumber\\
&=\left\|\hat{\vw}^{k+1}-\vw^{k}\right\|_{\mS}^{2}-\left\|\vw^{k}-\vw^{*}\right\|_{\mS}^{2}+\left\|\hat{\vw}^{k+1}-\vw^{*}\right\|_{\mS}^{2}\nonumber\\
&\quad - \left\|\hat{\vw}^{k+1}-\vw^{k}\right\|_{\mS}^{2}=\left\|\hat{\vw}^{k+1}-\vw^{*}\right\|_{\mS}^{2}-\left\|\vw^{k}-\vw^{*}\right\|_{\mS}^{2} \nonumber\\
& \leq-\left\|\hat{\vw}^{k+1}-\vw^{k}\right\|_{\widehat{\mS}}^{2}+\big\langle \hat{\vw}^{k+\frac{1}{2}}-\vw^*,\mD^k \big\rangle.
\end{align}
Substituting \eqref{ConPF2} into \eqref{ConPF1} we obtain
\begin{align}\label{ABC1}
\mathbb{E}\left[\left\|\vw^{k+1}-\vw^{*}\right\|_{\overline{\mS}}^{2}\right]\leq &\left\|\vw^{k}-\vw^{*}\right\|_{\overline{\mS}}^{2}-\left\|\vw^{k}-\hat{\vw}^{k+1}\right\|_{\widehat{\mS}}^{2}\nonumber\\ &\quad+\big\langle \hat{\vw}^{k+\frac{1}{2}}-\vw^*,\mD^k \big\rangle.
\end{align}
When $\mD^k=0$, we denote $\vw^{k+1}$ as $\breve{\vw}^{k+1}$.
Thus, by \eqref{ABC1}, we have $
\mathbb{E}[\|\breve{\vw}^{k+1}-\vw^{*}\|_{\overline{\mS}}^{2}] \leq\mathbb{E}[\|\vw^{k}-\vw^{*}\|_{\overline{\mS}}^{2}]$. Then, it holds
\begin{align}\label{Boundpf1}
&\mathbb{E}\left[\left\|\vw^{k+1}-\vw^{*}\right\|_{\overline{\mS}}^{2}\right]\nonumber\\
&\leq\mathbb{E}\left[\left\|\breve{\vw}^{k+1}-\vw^{*}\right\|_{\overline{\mS}}^{2}\right]+\mathbb{E}\left[\left\|\breve{\vw}^{k+1}-\vw^{k+1}\right\|_{\overline{\mS}}^{2}\right]\nonumber\\
&\leq\mathbb{E}\left\|\vw^{k}-\vw^{*}\right\|_{\overline{\mS}}^{2}+\mathbb{E}\left[\left\|\breve{\vw}^{k+1}-\vw^{k+1}\right\|_{\overline{\mS}}^{2}\right].
\end{align}
Let $\overline{\vw}^{k+1}=\mathbb{T}^{[\vvarepsilon^k=0]}\vw^k$ and we have he fact that $\hat{\vw}^{k+1}=\overline{\vw}^{k+1}$ with $\vvarepsilon^k$ setting to zero. Recall from \eqref{Alg3b} that
\begin{align*}
&\hat{\vx}^{k+1}=\prox_{nr}\left(\vx^k-(\mU\tr\vlambda^{k+\frac{1}{2}}-\vvarepsilon^k)\right),\\
&\overline{\vx}^{k+1}=\prox_{nr}\left(\vx^k-\mU\tr\vlambda^{k+\frac{1}{2}}\right).
\end{align*}
By the Lipschitz continuity of $\mathrm{prox}_{nr}(\cdot)$ one can readily get
\begin{align*}
\left\|\hat{\vx}^{k+1}-\overline{\vx}^{k+1}\right\|^2&\leq \left\langle\vvarepsilon^k,\hat{\vx}^{k+1}-\overline{\vx}^{k+1} \right\rangle\\
&\leq  \left\|\vvarepsilon^k\right\|\cdot\left\|\hat{\vx}^{k+1}-\overline{\vx}^{k+1}\right\|.
\end{align*}
Thus, we have $\|\hat{\vx}^{k+1}-\overline{\vx}^{k+1}\|\leq \|\vvarepsilon^k\|$. On the other hand,
\begin{align*}
\left\|\hat{\vlambda}^{k+1}-\overline{\vlambda}^{k+1}\right\|&\leq \beta\left\|\mU(\hat{\vx}^{k+1}-\overline{\vx}^{k+1})-(\hat{\vy}^{k+1}-\overline{\vy}^{k+1})\right\|\\
&\leq \beta \big\|\mU \big\|\cdot \left\|\vvarepsilon^k\right\|.
\end{align*}
Therefore, it holds that
\begin{align}\label{Bound1}
\left\|\hat{\vw}^{k+1}-\overline{\vw}^{k+1}\right\|_{\overline{\mS}}^2 \leq \underbrace{\left(1+ \beta^2\big\|\mU\big\|^2\right)\big\|\overline{\mS}\big\|}_{:=\mu^2>0}\cdot\big\|\vvarepsilon^k\big\|^2.
\end{align}
It follows from \eqref{OP1} that
\begin{align*}
&\vw^{k+1}=\vw^{k}+\mE^{k+1}(\hat{\vw}^{k+1}-\vw^{k}),\\
&\breve{\vw}^{k+1}=\vw^{k}+\mE^{k+1}(\overline{\vw}^{k+1}-\vw^{k}).
\end{align*}
Since $\mE^{k+1}$ is an idempotent matrix, it holds that
\begin{align*}
&\left\|\vw^{k+1}-\breve{\vw}^{k+1}\right\|= \left\|\mE^{k+1}(\hat{\vw}^{k+1}-\overline{\vw}^{k+1})\right\|\\
&\leq\left \|\mE^{k+1}\right\| \left\|\hat{\vw}^{k+1}-\overline{\vw}^{k+1}\right\| \leq \left\|\hat{\vw}^{k+1}-\overline{\vw}^{k+1}\right\|.
\end{align*}
Therefore, it deduces that
$\|\vw^{k+1}-\breve{\vw}^{k+1}\|_{\overline{\mS}}^{2}\leq \eta(\overline{\mS})\|\hat{\vw}^{k+1}-\overline{\vw}^{k+1}\|_{\overline{\mS}}^2
$.
Then, by \eqref{Boundpf1}, we have
\begin{align}
&\mathbb{E}\left[\left\|\vw^{k+1}-\vw^{*}\right\|_{\overline{\mS}}^{2}\right]\nonumber\\
&\leq\mathbb{E}\left\|\vw^{k}-\vw^{*}\right\|_{\overline{\mS}}^{2}+\eta(\overline{\mS})\mathbb{E}\left[\left\|\hat{\vw}^{k+1}-\overline{\vw}^{k+1}\right\|_{\overline{\mS}}^2\right]\nonumber\\
&\leq\mathbb{E}\left\|\vw^{k}-\vw^{*}\right\|_{\overline{\mS}}^{2} + \eta(\overline{\mS})\mu^2\left\|\vvarepsilon^k\right\|^2.
\end{align}
Summing the above inequality over $k=1,\cdots,\mathrm{K}$, one has for any $\mathrm{K}\geq 1$,
$
\sum_{k=1}^{\mathrm{K}}(\mathbb{E}\left\|\vw^{k+1}-\vw^{*}\right\|_{\overline{\mS}}^{2}-\mathbb{E}\left\|\vw^{k}-\vw^{*}\right\|_{\overline{\mS}}^{2})
\leq \sum_{k=1}^{\mathrm{K}} \eta(\overline{\mS})\mu^2\left\|\vvarepsilon^k\right\|^2
$,
which implies that for $\forall \mathrm{K}\geq 1$,
\begin{align*}
\mathbb{E}\left\|\vw^{\mathrm{K}}-\vw^{*}\right\|_{\overline{\mS}}^{2}\leq \mathbb{E}\left\|\vw^{1}-\vw^{*}\right\|_{\overline{\mS}}^{2}+ \sum_{k=1}^{\mathrm{K}} \eta(\overline{\mS})\mu^2\left\|\vvarepsilon^k\right\|^2.
\end{align*}

Since $\overline{\mS}$ is positive definite and $\{\vvarepsilon^k\}$ is summable, it deduces that for any $w^*\in \mathcal{W}^*$, $\{\|\vw^{k}-\vw^{*}\|\}$ is bounded, i.e., there exists $M_1>0$ such that $\|\vw^{k}-\vw^{*}\|\leq M_1$, for $\forall k\geq 0$.
Since proximal mapping and $\nabla f_i$ are Lipschitz continuous, and $\sum_{k=0}^{\infty}\|\vvarepsilon^k\|<\infty$, there exists $M_2>0$ satisfying
\begin{align}\label{ABC3}
\big\|\hat{\vw}^{k+\frac{1}{2}}-\vw^*\big\|\leq M_2\left\|\vw^{k}-\vw^*\right\|\leq M_1M_2, ~\forall k\geq0.
\end{align}
Therefore, from \eqref{ABC1} and \eqref{ABC3}, it holds that
\begin{align*}
&\mathbb{E}\left[\left\|\vw^{k+1}-\vw^{*}\right\|_{\overline{\mS}}^{2}\right] \leq\left\|\vw^{k}-\vw^{*}\right\|_{\overline{\mS}}^{2}\\
&\quad -\left\|\vw^{k}-\hat{\vw}^{k+1}\right\|_{\widehat{\mS}}^{2} +M_1M_2\cdot\big\|\vvarepsilon^k\big\|,
\end{align*}
which implies \eqref{ConPF3} holds. From it, the sequence $\{\vw^k\}$ is Quasi-Fej\'{e}r monotone with respect to $\mathcal{W}^*$ in $\overline{\mS}$-norm.
\end{IEEEproof}

\section{Proof of Theorem \ref{thm2}}\label{APP3}
\begin{IEEEproof}
It follows from \eqref{ConPF3} that
\begin{align}\label{KN1-1}
&\EE_k\left[\dist^2_{\overline{\mS}}(\vw^{k+1},\mathcal{W^*})\right]\leq\EE_k\left[\big\|\vw^{k+1}-\mathcal{P}_{\mathcal{W^*}}^{\overline{\mS}}(\vw^k)\big\|^2_{\overline{\mS}}\right]\nonumber\\
&\leq\big\|\vw^k-\mathcal{P}_{\mathcal{W^*}}^{\overline{\mS}}(\vw^k)\big\|^2_{\overline{\mS}}-\big\|\hat{\vw}^{k+1}-\vw^k\big\|_{\widehat{\mS}}^2+M_1M_2\cdot\big\|\vvarepsilon^k\big\|\nonumber\\
&=\dist^2_{\overline{\mS}}(\vw^{k},\mathcal{W^*})-\big\|\hat{\vw}^{k+1}-\vw^k\big\|_{\widehat{\mS}}^2+M_1M_2\big\|\vvarepsilon^k\big\|.
\end{align}
By \eqref{AFBS2}, we obtain
\begin{align}\label{Lind1}
&\big\|\overline{\vw}^{k+\frac{1}{2}}-\vw^k\big\|^2=\big\|(\mH+\mA\tr)^{-1}\mS(\overline{\vw}^{k+1}-\vw^k)\big\|^2\nonumber\\
&\leq\big\|(\mH+\mA\tr)^{-1}\mS\big\|^2\big\|\widehat{\mS}^{-1}\big\|\big\|\overline{\vw}^{k+1}-\vw^k)\big\|^2_{\widehat{\mS}}.
\end{align}
Considering the projection of $\hat{\vw}^{k+\frac{1}{2}}$ onto $\mathcal{W^*}$, we have
\begin{align}\label{Lind2}
&\dist^2_{\overline{\mS}}(\vw^{k},\mathcal{W^*})\leq\big\|\vw^{k}-\mathcal{P}_{\mathcal{W^*}}(\overline{\vw}^{k+\frac{1}{2}})\big\|^2_{\overline{\mS}}\nonumber\\
&\leq\big\|\overline{\mS}\big\| \big\|\vw^k-\mathcal{P}_{\mathcal{W^*}}(\overline{\vw}^{k+\frac{1}{2}})\big\|^2\nonumber\\
&\leq\big\|\overline{\mS}\big\|\left(\big\|\overline{\vw}^{k+\frac{1}{2}}-\mathcal{P}_{\mathcal{W^*}}(\overline{\vw}^{k+\frac{1}{2}})\big\|+\big\|\overline{\vw}^{k+\frac{1}{2}}-\vw^k\big\|\right)^2\nonumber\\
&=\big\|\overline{\mS}\big\|\left(\dist(\overline{\vw}^{k+\frac{1}{2}},\mathcal{W^*})+\big\|\overline{\vw}^{k+\frac{1}{2}}-\vw^k\big\|\right)^2.
\end{align}
Since $\overline{\vw}^{k+1}=\mathbb{T}^{[\vvarepsilon^k=0]}\vw^k$, i.e., $\overline{\vw}^{k+1/2}=(\mH+\mathbf{G})^{-1}(\mH-\mA-\FP)\vw^k$, it holds that
\begin{align*}
(\mH-\mA)(\vw^{k}-\overline{\vw}^{k+\frac{1}{2}})+\FP(\overline{\vw}^{k+\frac{1}{2}})-\FP(\vw^k)  \in \partial \mathcal{L}(\overline{\vw}^{k+\frac{1}{2}}).
\end{align*}
We can get that
\begin{align*}
&\mathrm{dist}(0,\partial \mathcal{L}(\overline{\vw}^{k+\frac{1}{2}}))\\
&\leq\big\|(\mH-\mA)(\vw^{k}-\overline{\vw}^{k+\frac{1}{2}})+\FP(\overline{\vw}^{k+\frac{1}{2}})-\FP(\vw^k)\big\|\\
&\leq\big\|\nabla f(\overline{\vy}^{k+\frac{1}{2}})-\nabla f(\vy^{k})\big\|+\left(\big\|\mH-\mA\big\|\right)\big\|\vw^k-\overline{\vw}^{k+\frac{1}{2}}\big\|\\
&\leq\left(\big\|\mQ\big\|+\big\|\mH-\mA\big\|\right)\cdot\big\|\overline{\vw}^{k+\frac{1}{2}}-\vw^k\big\|.
\end{align*}
Since the mapping $\partial \mathcal{L}$ is metrically subregular at $(\vw^{\infty},0)$ with modulus $\tilde{\varsigma}$, there exists $\epsilon>0$ such that
\begin{align*}
&\mathrm{dist}(\overline{\vw}^{k+\frac{1}{2}},(\partial \mathcal{L})^{-1}(0))\leq\tilde{\varsigma}\cdot \mathrm{dist}(0,\partial \mathcal{L}(\overline{\vw}^{k+\frac{1}{2}}))\\
&\leq\tilde{\varsigma}\cdot\left(\big\|\mQ\big\|+\big\|\mH-\mA\big\|\right)\cdot\big\|\overline{\vw}^{k+\frac{1}{2}}-\vw^k\big\|, \overline{\vw}^{k+\frac{1}{2}}\in \mathcal{B}_{\epsilon}(\vw^{\infty}).
\end{align*}
It follows from $(\partial \mathcal{L})^{-1}(0)=\mathcal{W}^*$ that if $\overline{\vw}^{k+\frac{1}{2}}\in \mathcal{B}_{\epsilon}(\vw^{\infty})$,
\begin{align}\label{EB1}
&\mathrm{dist}(\overline{\vw}^{k+\frac{1}{2}},\mathcal{W}^*)\leq\mathrm{dist}(\overline{\vw}^{k+\frac{1}{2}},(\partial \mathbf{L})^{-1}(0))\nonumber\\
&\leq\tilde{\varsigma}\cdot \left(\big\|\mQ\big\|+\big\|\mH-\mA\big\|\right)\big\|\overline{\vw}^{k+\frac{1}{2}}-\vw^k\big\|.
\end{align}
Since $\lim_{k\rightarrow \infty}\vw^{k}=\vw^{\infty}, \text{a.s.}$ and $\lim_{k\rightarrow \infty}\|\vvarepsilon^k\|=0$, it follows from \eqref{Bound1} that $\{\overline{\vw}^k\}$ also converges to $\vw^{\infty}$ as $k\rightarrow \infty$. From \eqref{AFBS2}, we have $\|\overline{\vw}^{k+1}-\overline{\vw}^{k+\frac{1}{2}}\| \leq (1+\|(\mH+\mA\tr)^{-1}\mS\|)\|\overline{\vw}^{k+1}-\vw^k\|$. Since $\lim_{k\rightarrow \infty}\|\overline{\vw}^{k+1}-\vw^k\|=0$, for the $\epsilon>0$, there exists an integer $\tilde{K}>0$ such that $\overline{\vw}^{k+\frac{1}{2}}\in \mathcal{B}_{\epsilon}(\vw^{\infty}),\forall k\geq \tilde{K}$. Let $\tilde{\epsilon}=\min_{0\leq k<\tilde{K}}\{\|\overline{\vw}^{k+\frac{1}{2}}-\vw^k\|\}$. It holds that $\|\overline{\vw}^{k+\frac{1}{2}}-\vw^k\|<\tilde{\epsilon}\Rightarrow k\geq \tilde{K}\Rightarrow\overline{\vw}^{k+\frac{1}{2}}\in \mathcal{B}_{\epsilon}(\vw^{\infty})$, i.e., for the $\epsilon>0$, there exists $\tilde{\epsilon}>0$ such that
$\|\overline{\vw}^{k+\frac{1}{2}}-\vw^k\|<\tilde{\epsilon}$ implies that $\overline{\vw}^{k+\frac{1}{2}}\in \mathcal{B}_{\epsilon}(\vw^{\infty})$. Therefore, by \eqref{EB1}, when $\|\overline{\vw}^{k+\frac{1}{2}}-\vw^k\|<\tilde{\epsilon}$, it gives that
\begin{align*}
&\mathrm{dist}(\overline{\vw}^{k+\frac{1}{2}},\mathcal{W}^*)\leq \tilde{\varsigma} \cdot\left(\big\|\mQ\big\|+\big\|\mH-\mA\big\|\right)\big\|\overline{\vw}^{k+\frac{1}{2}}-\vw^k\big\|.
\end{align*}
By \eqref{ConPF3}, we have $\|\overline{\vw}^{k+\frac{1}{2}}-\vw^k\|\rightarrow 0$, which implies that  $\|\overline{\vw}^{k+\frac{1}{2}}-\vw^k\|<\tilde{\epsilon}$ when $k$ is sufficient large. Thus, we only need to consider indices $k<\tilde{K}$ (without loss of generality, we still use $\tilde{K}$ to denote a sufficiently large positive integer) such that $\|\overline{\vw}^{k+\frac{1}{2}}-\vw^k\|\geq\tilde{\epsilon}$. According to the boundness of $\|\overline{\vw}^{k+\frac{1}{2}}-\vw^{\infty}\|$, it holds that there exists a constant $\tilde{M}>0$ such that $\|\overline{\vw}^{k+\frac{1}{2}}-\vw^{\infty}\|\leq \tilde{M}, \forall k\geq0$. It immediately has
\begin{align*}
\mathrm{dist}(\overline{\vw}^{k+\frac{1}{2}},\mathcal{W}^*)&\leq \big\|\overline{\vw}^{k+\frac{1}{2}}-\vw^{\infty}\big\|\leq \tilde{M}/\tilde{\epsilon}\cdot \big\|\overline{\vw}^{k+\frac{1}{2}}-\vw^k\big\|.
\end{align*}
Let $\varsigma=\max\{\tilde{\varsigma} \cdot (\|\mQ\|+\|\mH-\mA\|),\tilde{M}/\tilde{\epsilon}\}$. It follows that
\begin{align}\label{LinIn1}
\mathrm{dist}(\overline{\vw}^{k+\frac{1}{2}},\mathcal{W^*})\leq {\varsigma} \cdot\big\|\overline{\vw}^{k+\frac{1}{2}}-\vw^{k}\big\|, \quad \forall k\geq 0.
\end{align}
Combine \eqref{Lind1}, \eqref{Lind2} and \eqref{LinIn1} to deduce
$
\dist^2_{\overline{\mS}}(\vw^{k},\mathcal{W^*})\leq \phi \|\overline{\vw}^{k+1}-\vw^k)\|^2_{\widehat{\mS}}
$,
where $\phi=(\varsigma+1)^2\|\overline{\mS}\|\|(\mH+\mA\tr)^{-1}\mS\|^2\|\widehat{\mS}^{-1}\|$. With \eqref{Bound1}, we get
\begin{align}\label{Lind3}
&\dist^2_{\overline{\mS}}(\vw^{k},\mathcal{W^*})\leq \phi \cdot \big\|\overline{\vw}^{k+1}-\hat{\vw}^{k+1}+\hat{\vw}^{k+1}-\vw^k)\big\|^2_{\widehat{\mS}}\nonumber\\
&\leq \phi \left(\big\|\overline{\vw}^{k+1}-\hat{\vw}^{k+1}\big\|^2_{\widehat{\mS}} +\big\|\hat{\vw}^{k+1}-\vw^k\big\| ^2_{\widehat{\mS}}\right) \nonumber\\
&\leq \phi\left(\big\|\hat{\vw}^{k+1}-\vw^k\big\|^2_{\widehat{\mS}}+\mu^2\big\|\vvarepsilon^k\big\|^2\right).
\end{align}
By \eqref{KN1-1} and \eqref{Lind3} we have
\begin{align}
&\EE_k\left[\dist^2_{\overline{\mS}}(\vw^{k+1},\mathcal{W^*})\right]\nonumber\\
&\leq (1-\frac{1}{\phi})\dist^2_{\overline{\mS}}(\vw^{k},\mathcal{W^*})+(M_1M_2+\mu^2)\big\|\vvarepsilon^k\big\|^2.
\end{align}
From \eqref{ConPF3}, it holds that
\begin{align*}
\EE_k\left[\dist^2_{\overline{\mS}}(\vw^{k},\mathcal{W^*})\right]\geq\dist^2_{\overline{\mS}}(\vw^{k+1},\mathcal{W^*})-M_1M_2\big\|\vvarepsilon^k\big\|^2.
\end{align*}
Therefore, we have
\begin{align*}
&\EE_k\left[\dist^2_{\overline{\mS}}(\vw^{k+1},\mathcal{W^*})\right]\leq\frac{\phi}{\phi+1}\left[\dist^2_{\overline{\mS}}(\vw^{k},\mathcal{W^*})+V\|\vvarepsilon^k\|^2\right]
\end{align*}
where $V=M_1M_2+\frac{\phi}{\phi+1}\mu^2$.
Note that two updates of any agent are separated by at most $\ell$ iterations. Since $\sigma^2_{\tau_i(k)}=R\sigma^2_{\tau_i(k)+1}$, and $(\frac{\phi}{\phi+1})^\ell R>1$, we have
\begin{align*}
&\mathbb{E}\left[\dist^2_{\overline{\mS}}(\vw^{k+1},\mathcal{W^*})\right]\nonumber\\
&\leq (\frac{\phi}{\phi+1})^k\left(\dist^2_{\overline{\mS}}(\vw^{1},\mathcal{W^*})+V\sum_{s=1}^{k}(\frac{\phi}{\phi+1})^{-s}\big\|\vvarepsilon^s\big\|^2\right)\nonumber\\
&\leq (\frac{\phi}{\phi+1})^k\Bigg(\dist^2_{\overline{\mS}}(\vw^{1},\mathcal{W^*})+V\sum_{s=1}^{\lceil\frac{k}{\ell}\rceil}[(\frac{\phi}{\phi+1})^\ell]^{-s}\frac{(\sigma_1)^2}{R^s}\Bigg)\nonumber\\
&\leq(\frac{\phi}{\phi+1})^k\Bigg(\dist^2_{\overline{\mS}}(\vw^{1},\mathcal{W^*})+(\frac{\phi}{\phi+1})^\ell R\frac{V(\sigma_1)^2}{(\frac{\phi}{\phi+1})^\ell R-1}\Bigg).
\end{align*}
Finally, it is easy to see that $0<\frac{\phi}{\phi+1}<1$.
\end{IEEEproof}

\section{Experimental supplementation}
\subsection{Decentralized classification problems on logistic regression}
We classify each of the four datasets on the following logistic regression model in a decentralized manner.
\begin{align*}
&\textbf{Logistic Regression:}\\
&\min _{\vx}\frac{1}{n}\sum_{i=1}^{n}\frac{1}{m_i}\sum_{j=1}^{m_{i}} \ln(1+e^{-(\mB_{ij}\vx)\vb_{ij}})+\frac{1}{2}\|\vx\|^2+\frac{1}{2}\|\vx\|_1.
\end{align*}

The details and conclustions of this experiment are the same as those in Section VII-A of the main text.

\begin{figure*}[!h]
\centering
\subfigure{
\includegraphics[width=0.25\linewidth]{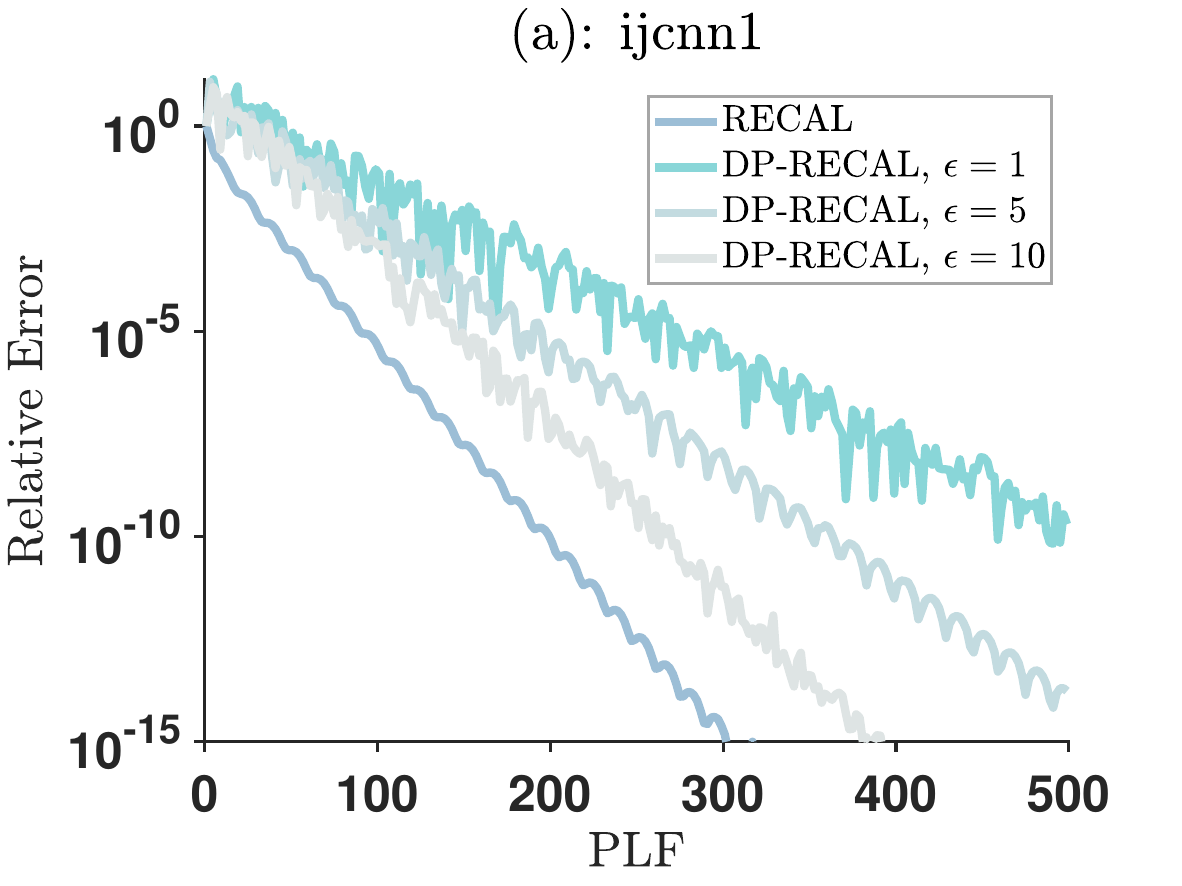}}\hspace{-2.85mm}
\subfigure{
\includegraphics[width=0.25\linewidth]{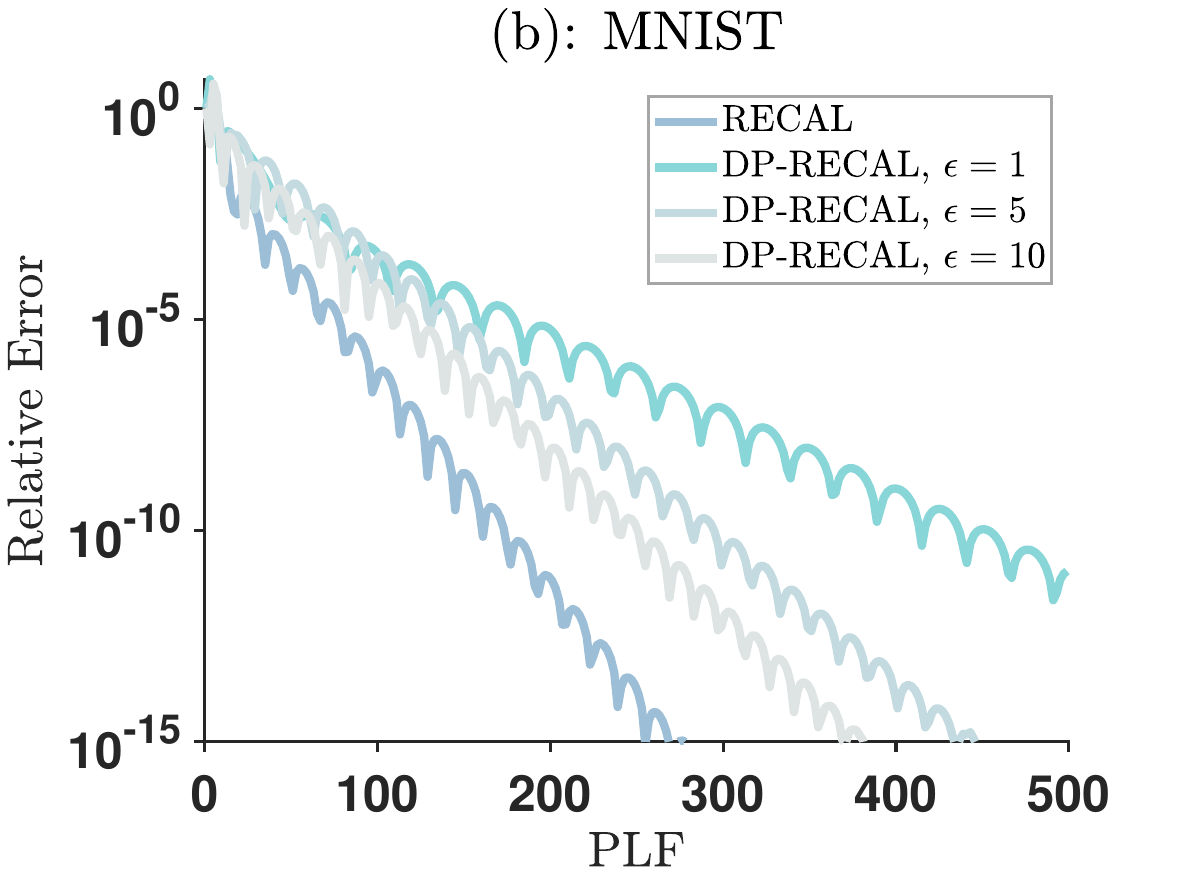}}\hspace{-2.85mm}
\subfigure{
\includegraphics[width=0.25\linewidth]{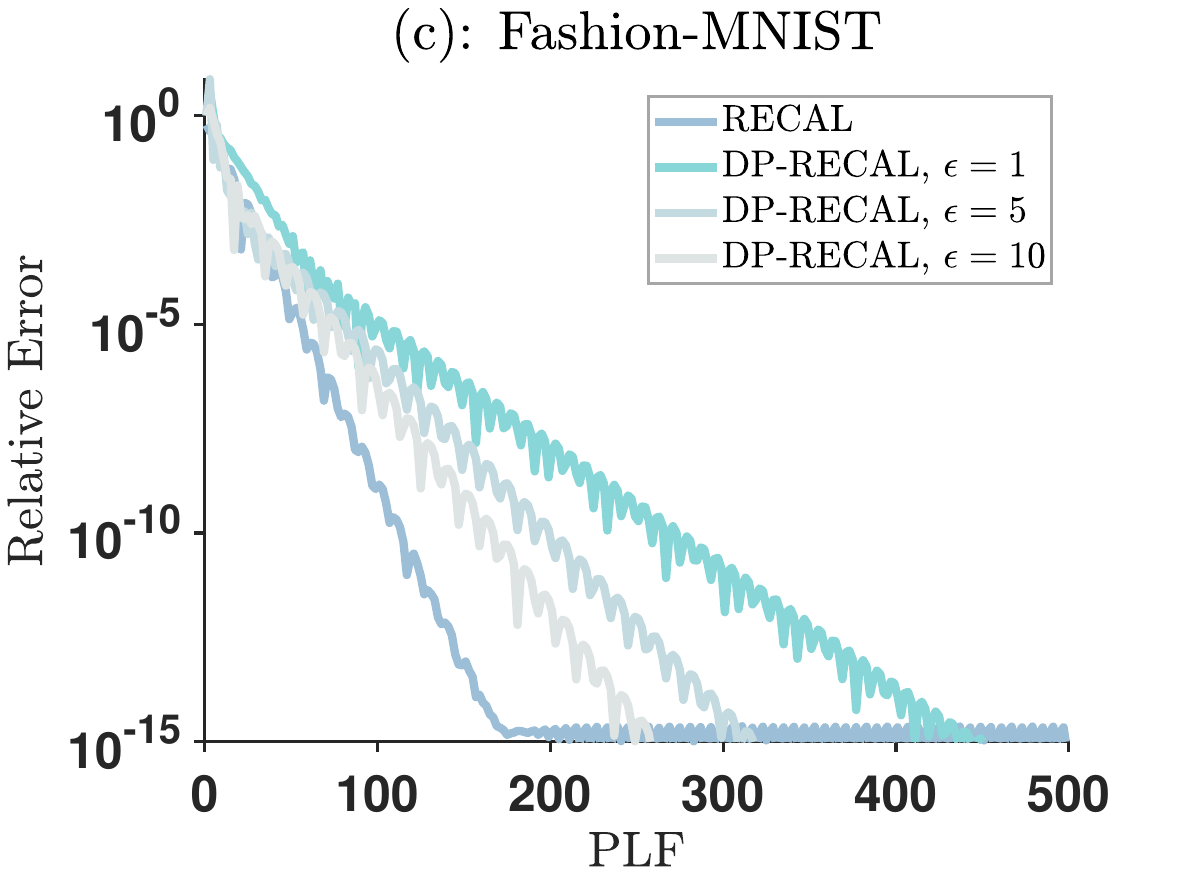}}\hspace{-2.85mm}
\subfigure{
\includegraphics[width=0.25\linewidth]{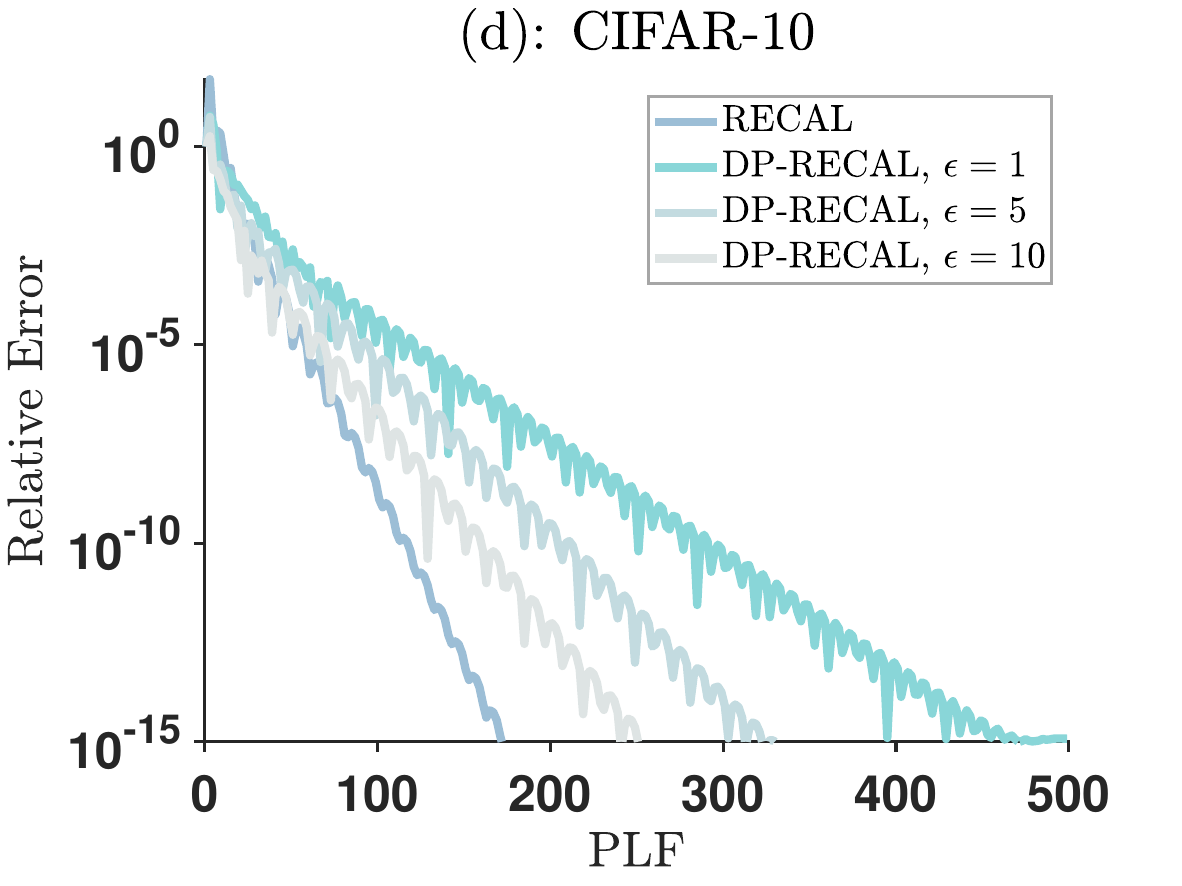}}
\caption{Convergence curves of DP-RECAL with various privacy budgets for the logistic regression problem on four types of datasets.}
\label{fig:diffep_LG}
\end{figure*}

\begin{figure*}[!t]
\centering
\subfigure{
\includegraphics[width=0.25\linewidth]{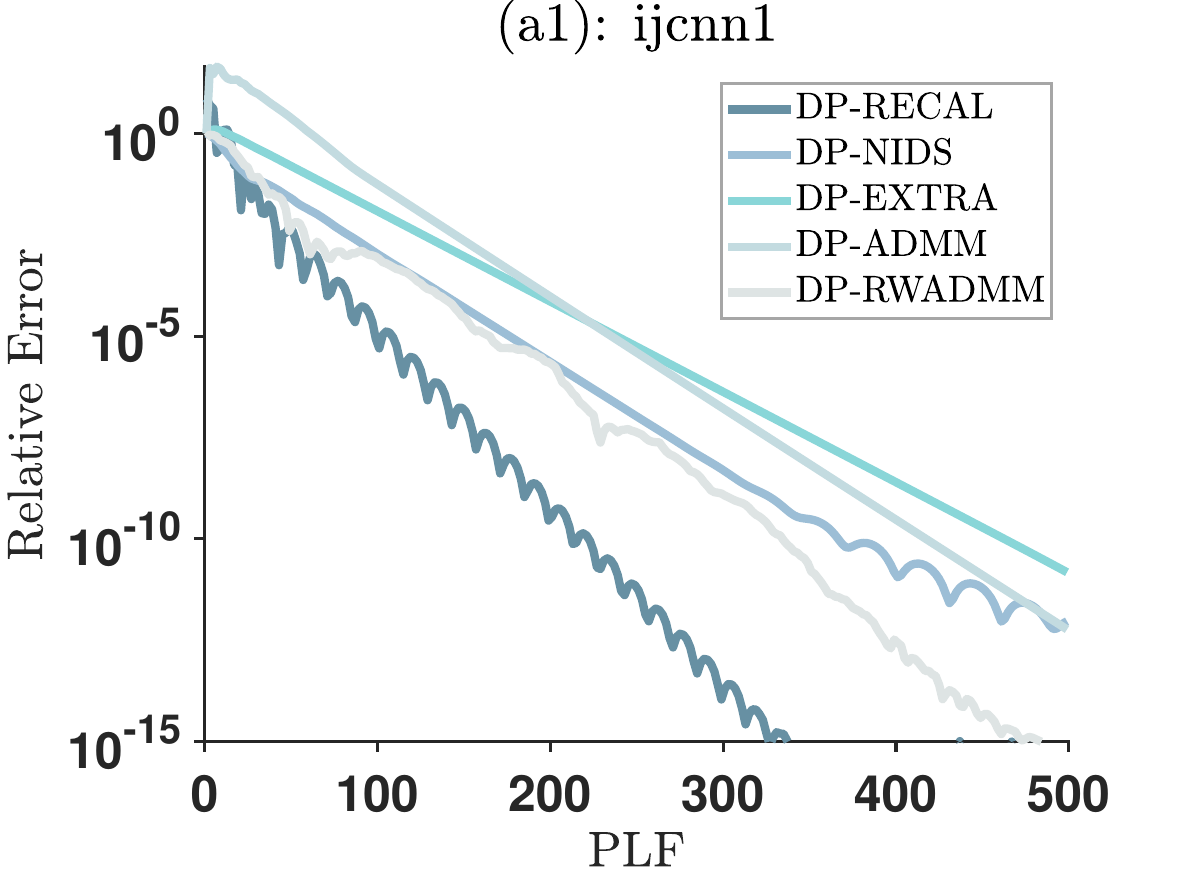}}\hspace{-2.85mm}
\subfigure{
\includegraphics[width=0.25\linewidth]{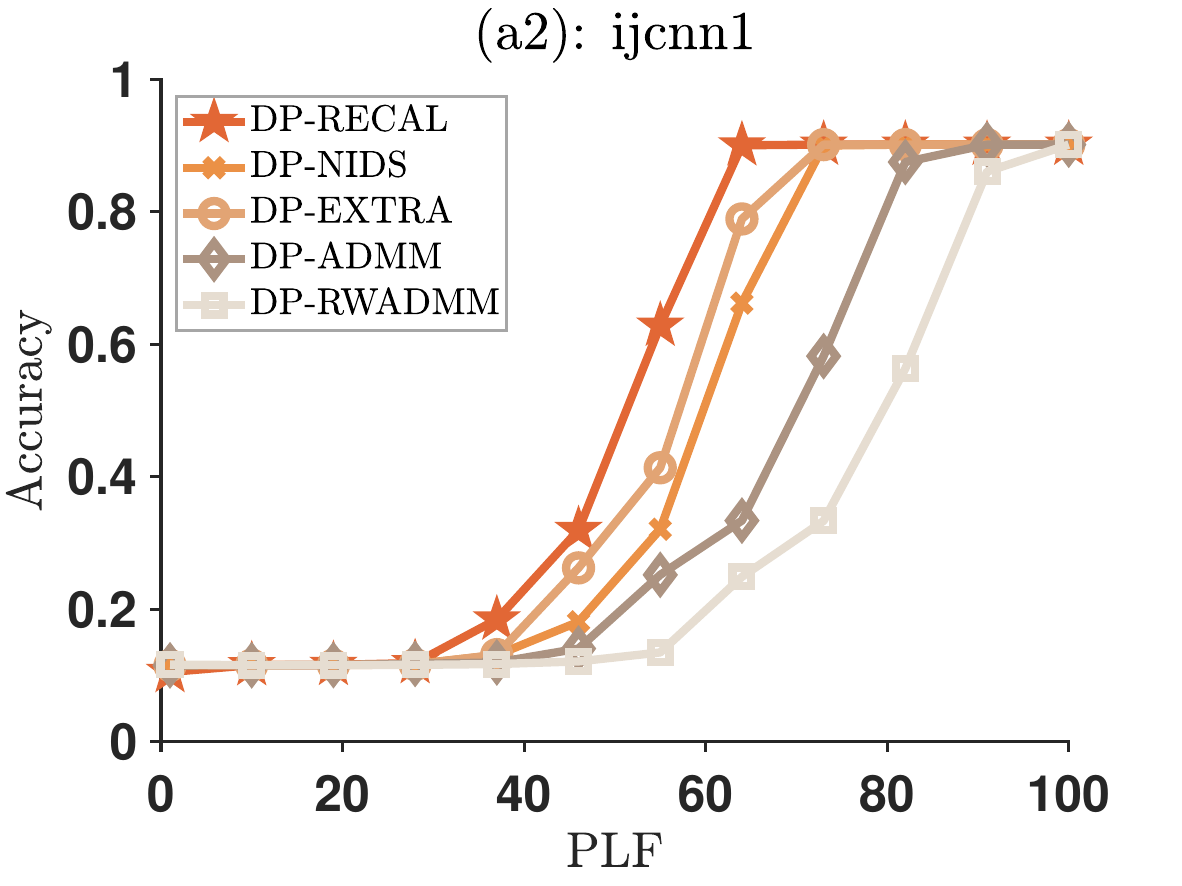}}\hspace{-2.85mm}
\subfigure{
\includegraphics[width=0.25\linewidth]{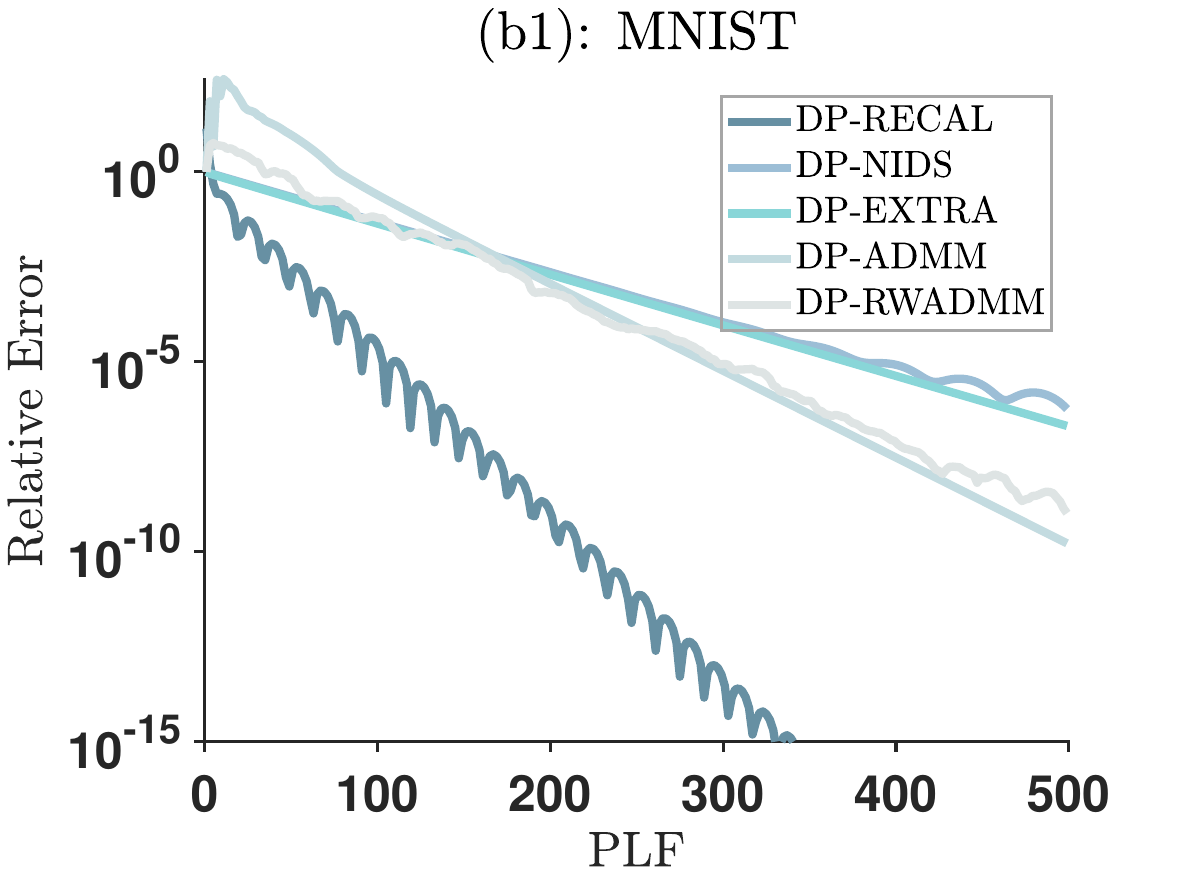}}\hspace{-2.85mm}
\subfigure{
\includegraphics[width=0.25\linewidth]{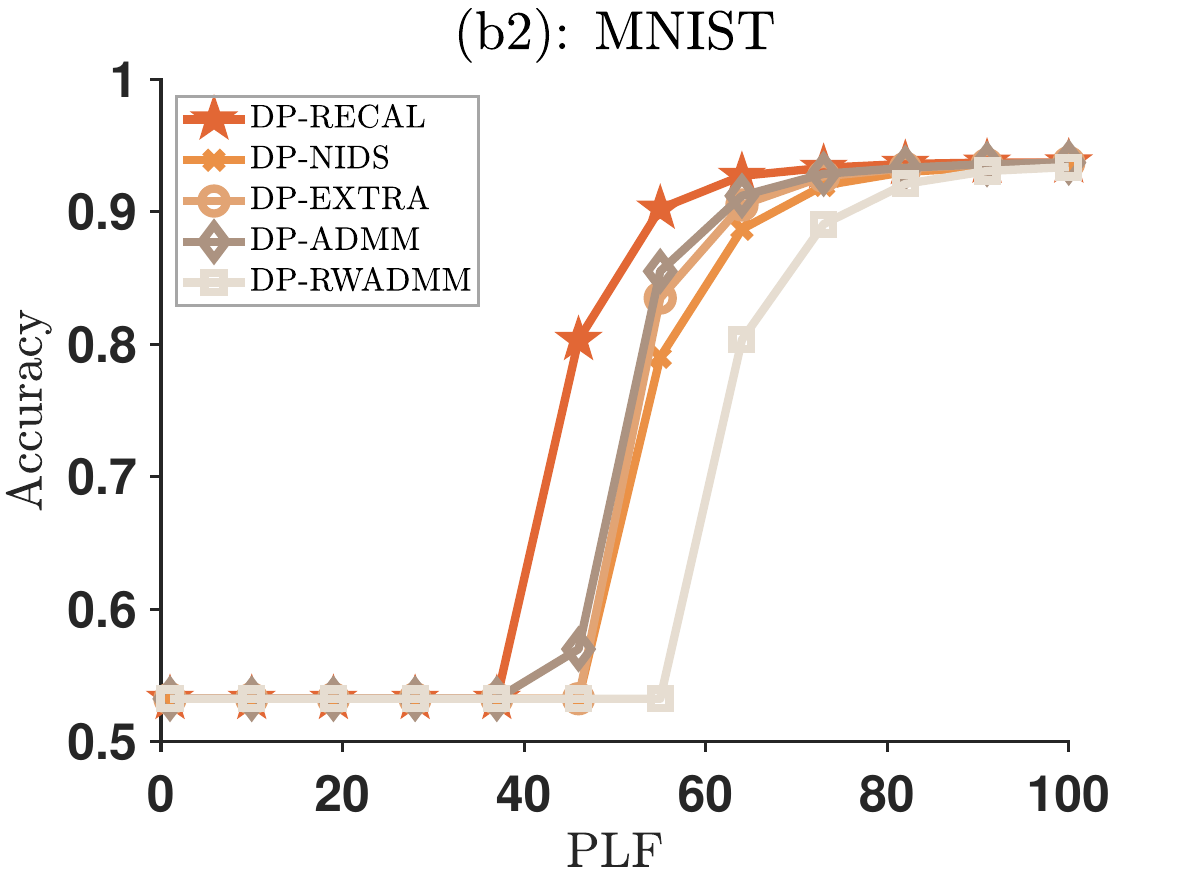}}
\subfigure{
\includegraphics[width=0.25\linewidth]{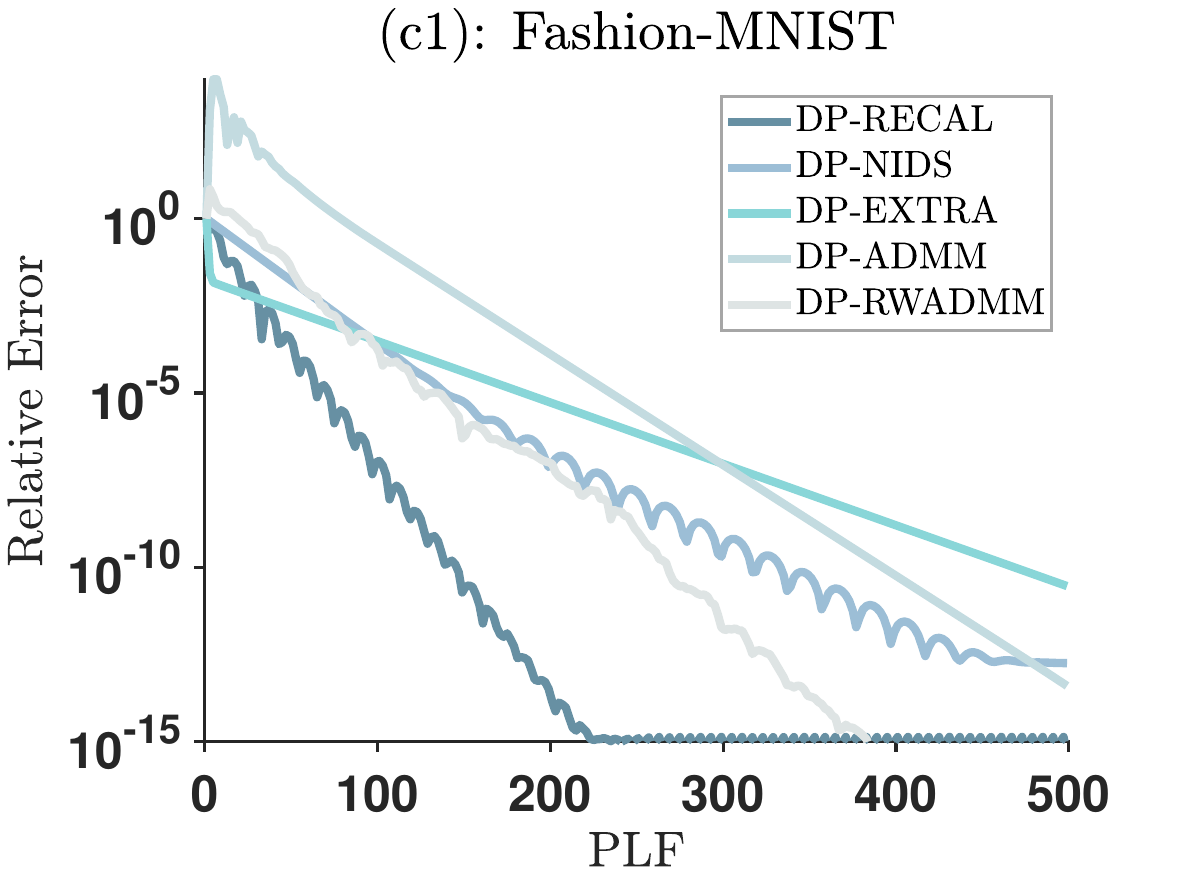}}\hspace{-2.85mm}
\subfigure{
\includegraphics[width=0.25\linewidth]{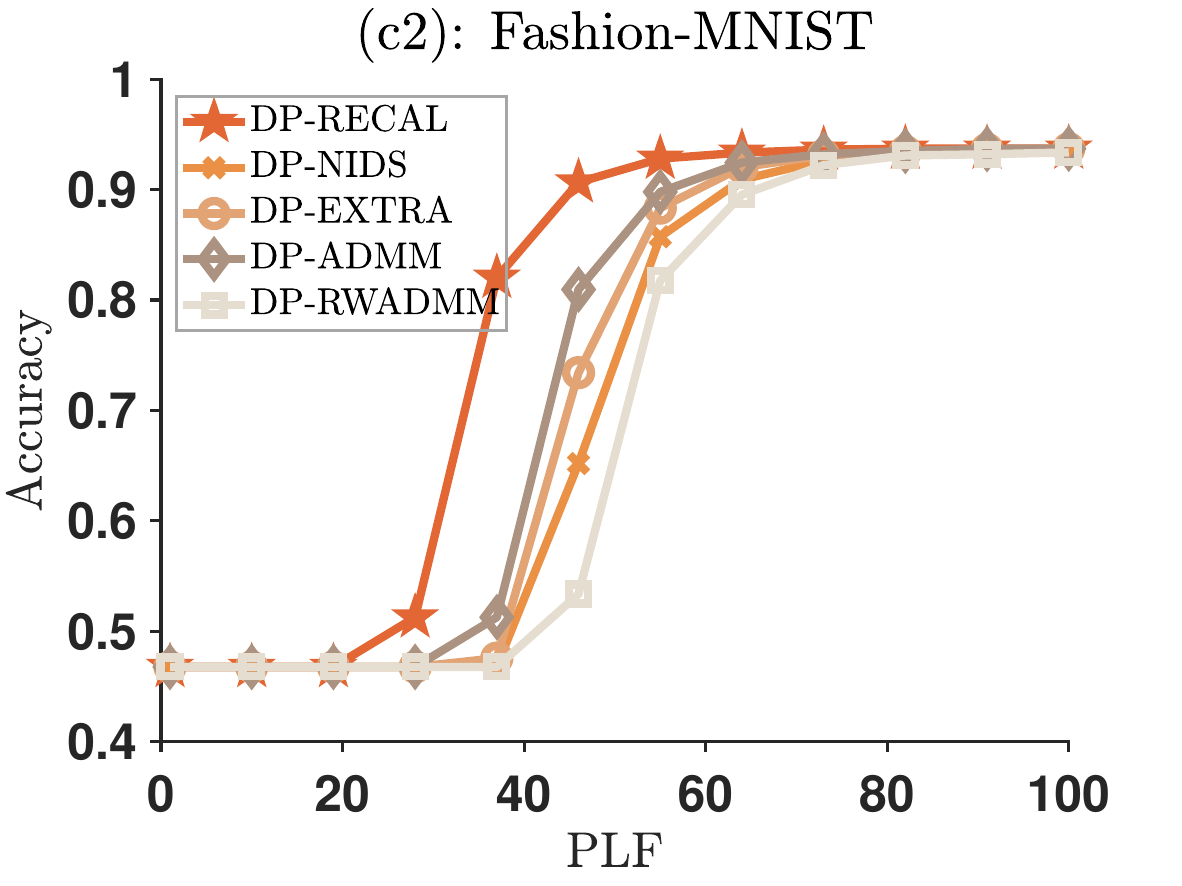}}\hspace{-2.85mm}
\subfigure{
\includegraphics[width=0.25\linewidth]{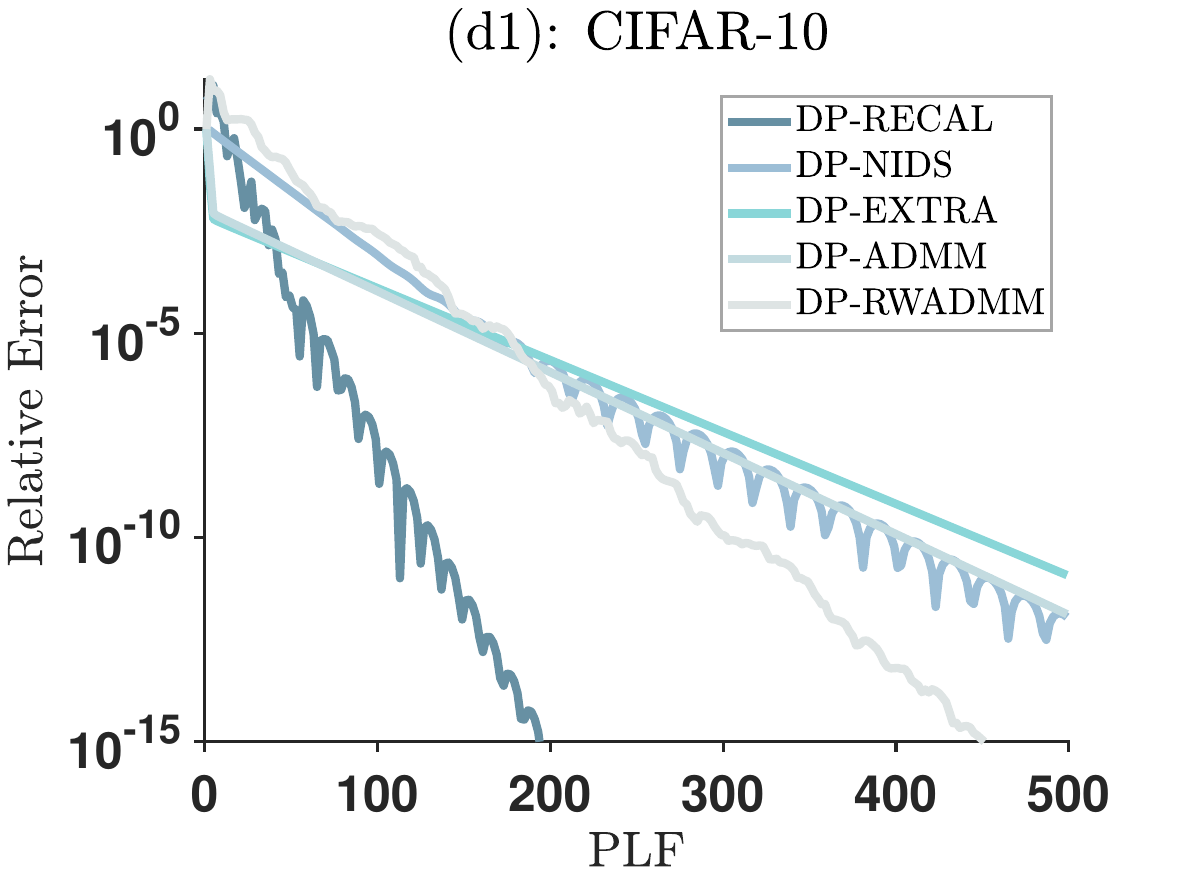}}\hspace{-2.85mm}
\subfigure{
\includegraphics[width=0.25\linewidth]{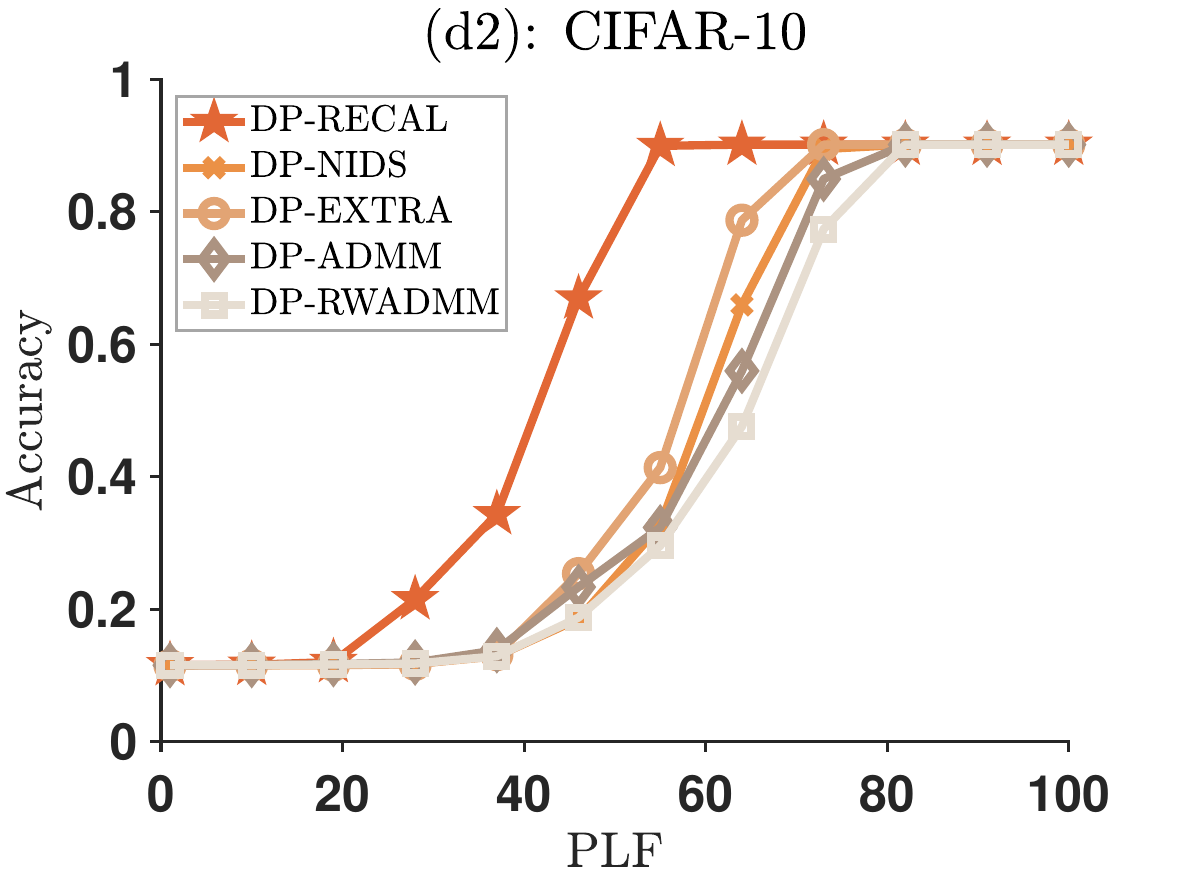}}
\caption{Convergence performance and classification accuracy of various differentially private algorithms with privacy budget $\epsilon = 12$ for the logistic regression problems on four types of datasets.}
\label{fig:compare_LG}
\end{figure*}

From Fig. \ref{fig:diffep_LG}, we can see that DP-RECAL achieves faster convergence with a larger privacy budget $\epsilon$. Fig. \ref{fig:compare_LG} illustrates that DP-RECAL consistently presents the best convergence accuracy when consuming the same privacy budget for logistic regression on different datasets. Moreover, Table \ref{Table-total_LG}, where Comm. and RE. stand for communication complexity and relative error, reveals that at the same level of security, DP-RECAL holds not only better convergence accuracy but also less communication complexity compared to other algorithms. As far as CPU time is concerned, DP-RECAL doesn't seem to be inferior or even better, especially compared to those algorithms with decaying stepsizes and internal loops.

\begin{table*}[!h]
\renewcommand\arraystretch{1.5}
\begin{center}
\caption{Comparisons of various differentially private algorithms for logistic regression problems on four types of datasets with the same privacy budget $\epsilon = 12$.}
\scalebox{1.05}{
\begin{tabular}{ccccccccccccc}
\hline
\multirow{2}{*}{$r(\vx)$}&\multirow{2}{*}{$\alpha$} &\multirow{2}{*}{PLF}&\multirow{2}{*}{ALG}&\multirow{2}{*}{Comm.}&\multicolumn{2}{c}{ijcnn1}&\multicolumn{2}{c}{MNIST}&\multicolumn{2}{c}{Fashion-MNIST}&\multicolumn{2}{c}{CIFAR-10}\\
\cline{6-13}
&&&&&{Time}&{RE.}&Time&RE.&Time&RE.&Time&RE.\\
\hline
\multirow{5}{*}{$\ell_1$}&$\mathcal{O}(\gamma^{-k})$&2000&DP-ADMM\cite{Huang2020}&6400&67.2&6.7e-8&173.6&8.5e-5&93.4&2.9e-8&8134&1.7e-8\\
\cline{2-13}
&\multirow{4}{*}{$\mathcal{O}(1)$}&\multirow{4}{*}{300}&DP-NIDS\cite{NDIS}&4800&13.4&1.1e-7&16.8&1.1e-4&17.3&3.0e-9&1328&8.5e-9\\
&&&DP-EXTRA\cite{PGEXTR}&9600&12.3&4.2e-6&13.9&8.1e-5&18.9&9.2e-8&1413&3.8e-8\\
&&&DP-ADMM\cite{Zhang2017}&9600&11.7&2.4e-8&14.3&5.6e-6&18.2&8.6e-8&1582&1.1e-8\\
&&&DP-RWADMM\cite{Shah2018}&9600&24.1&1.2e-9&18.9&8.3e-6&26.7&1.9e-12&3011&1.0e-10\\
\hline
\rowcolor{mygray}$\ell_1$&$\mathcal{O}(1)$&\textbf{300}&\textbf{DP-RECAL}&\textbf{2400}&\textbf{11.8}&\textbf{5.8e-15}&\textbf{12.6}&\textbf{4.0e-15}&\textbf{18.6}&\textbf{1.4e-15}&\textbf{1375}&\textbf{4.8e-16}\\
\hline
\hline
\multirow{2}{*}{-}&$\mathcal{O}(\gamma^{-k})$&2000&DP-MRADMM\cite{Zhang2018}&128000&1376&7.3e-9&2513&4.4e-7&2360&3.4e-11&2.8e+5&3.3e-11\\
\cline{2-13}
&$\mathcal{O}(1)$&300&DP-MRADMM\cite{Zhang2018}&19200&392&1.8e-10&988&2.3e-9&1000&5.8e-12&1.2e+5&9.4e-12\\
\hline
\rowcolor{mygray}-&$\mathcal{O}(1)$&\textbf{300}&\textbf{DP-RECAL}&\textbf{2400}&\textbf{6.1}&\textbf{8.2e-15}&\textbf{10.9}&\textbf{4.8e-15}&\textbf{20.2}&\textbf{4.6e-16}&\textbf{985}&\textbf{2.6e-15}\\
\hline
\end{tabular}}
\label{Table-total_LG}
\end{center}
\end{table*}

\subsection{Attack simulation with two other attack methods}
\begin{figure*}[!t]
\setlength{\abovecaptionskip}{-7pt}
  \centering
  \includegraphics[width=1\linewidth]{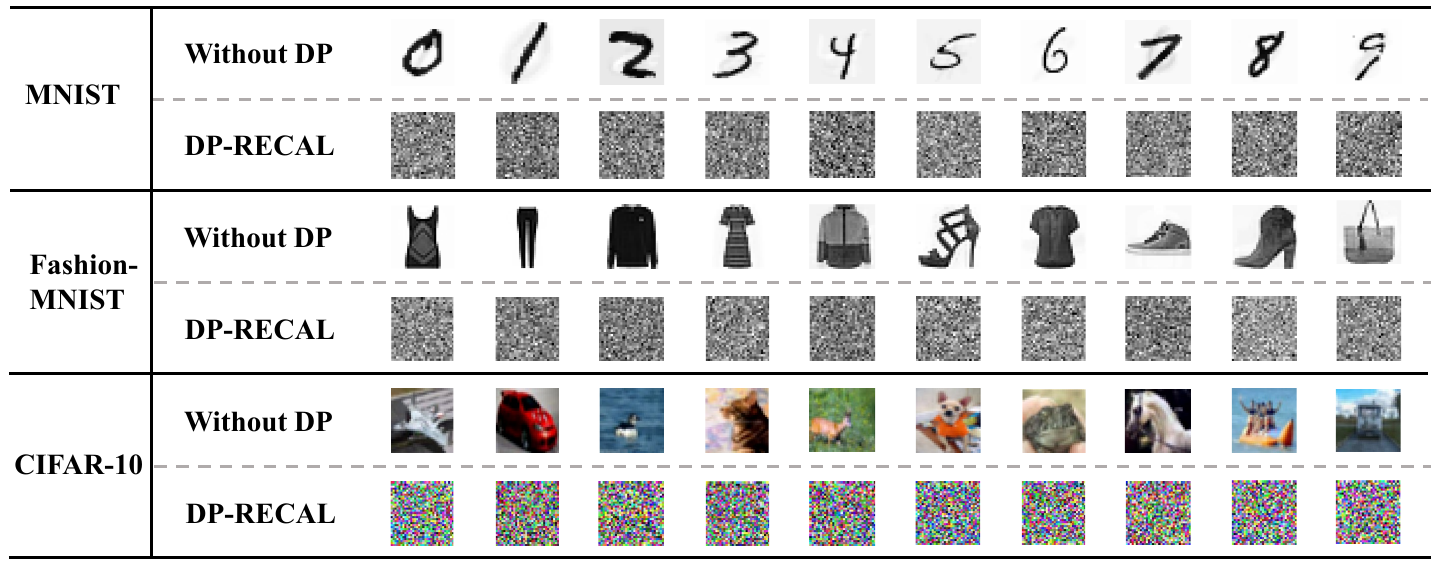}
  \caption{Comparisons of the results achieved by the adversary with attack method \cite{Geiping2020} on classification of three datasets performed by DP-RECAL and RECAL.}
  \label{fig:dlg1}
\end{figure*}

\begin{figure*}[!t]
\setlength{\abovecaptionskip}{-7pt}
  \centering
  \includegraphics[width=1\linewidth]{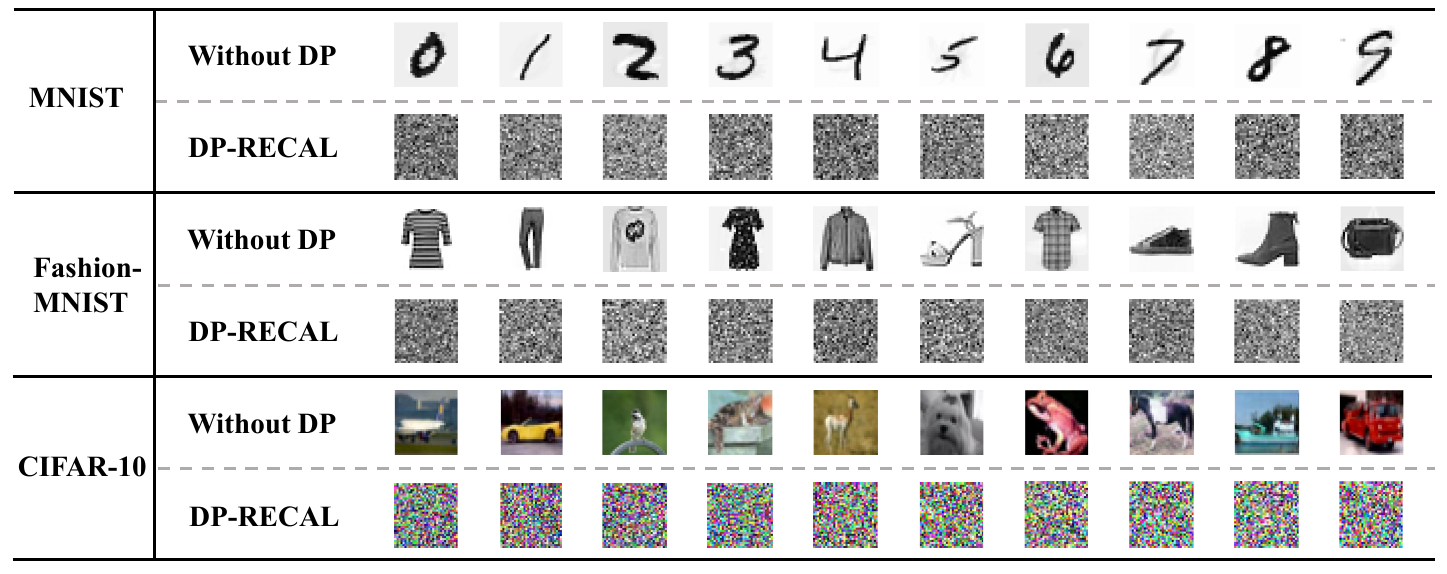}
  \caption{Comparisons of the results achieved by the adversary with attack method \cite{Yin2021} on classification of three datasets performed by DP-RECAL and RECAL.}
  \label{fig:dlg2}
\end{figure*}

We depict Fig. \ref{fig:dlg1} and Fig. \ref{fig:dlg2} to examine the security performance of our algorithm DP-RECAL against the latest gradient leakage attack methods \cite{Geiping2020} and \cite{Yin2021} through the classification process of three datasets. Both Fig. \ref{fig:dlg1} and Fig. \ref{fig:dlg2} show that DP-RECAL effectively thwarts the adversary’s attempts to infer the original images, while RECAL fails to defend against gradient leakage attasks. This also underscores the need for differential privacy mechanisms to protect sensitive data.

\subsection{Attack simulation on ijcnn1 datasets}
\begin{figure*}[!h]
  \centering
  \includegraphics[width=0.4\linewidth]{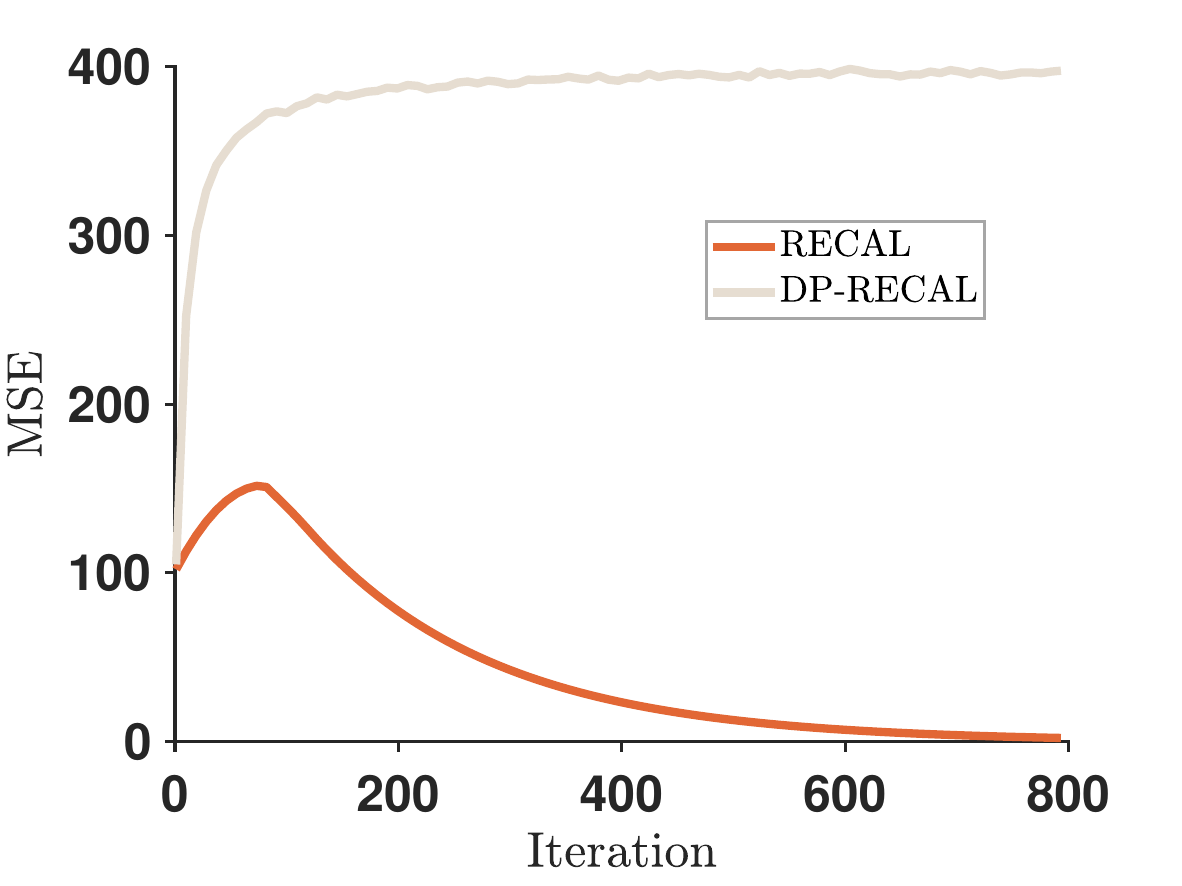}
  \caption{Mean-square error between the original and inferred data under RECAL and DP-RECAL.}
  \label{fig:dlg_ijcnn}
\end{figure*}

We continue to use the adversary in Section VII-B of the main text to attack the process of classifying the ijcnn1 dataset by the RECAL and DP-RECAL algorithms, respectively. Since ijcnn1 dataset cannot be presented as images, we use the mean square error between the original data and the inferred data during the iterations of the DLG algorithm to show the effect of the adversary's attacks. As illustrated in Fig. \ref{fig:dlg_ijcnn}, the attacker managed to accurately infer the unaltered image in the case of RECAL, as indicated by the convergence of estimation error to zero. Conversely, DP-RECAL effectively countered the attacker, as the estimation error remained consistently large. Combined with the results in the main text on other datasets, we can conclude that DP-RECAL is privacy-preserving.

\end{document}